\documentclass{imsart}

\RequirePackage{amsthm,amsmath,amsfonts,amssymb,subfig,bm}
\RequirePackage[numbers,sort&compress]{natbib}
\RequirePackage[colorlinks,citecolor=blue,urlcolor=blue]{hyperref}
\RequirePackage{graphicx}

\theoremstyle{plain}

\newtheorem{corollary}{Corollary}
\newtheorem{proposition}{Proposition}

\newtheorem{theorem}{Theorem}[section]
\newtheorem{lemma}[theorem]{Lemma}

\theoremstyle{remark}
\newtheorem{definition}[theorem]{Definition}

\newcommand{\E}{\mathbb{E}}
\newcommand{\Var}{\mathrm{Var}}

\newcommand{\Cov}{\mathrm{Cov}}
\renewcommand{\P}{\mathbb{P}}

\begin{document}

\begin{frontmatter}
\title{A minimum Wasserstein distance approach to Fisher's combination of independent discrete p-values}

\begin{aug}
\author[A]{\fnms{Gonzalo}~\snm{Contador}\ead[label=e1]{gonzalo.contador@usm.cl}} \and
\author[B]{\fnms{Zheyang}~\snm{Wu}\ead[label=e2]{zheyangwu@wpi.edu}}

\address[A]{Department of Mathematics, Universidad T\'ecnica Federico Santa Mar\'ia, Santiago, Chile\printead[presep={,\ }]{e1}}

\address[B]{Department of Mathematical Sciences, Worcester Polytechnic Institute, Worcester, MA, U.S.A.\printead[presep={,\ }]{e2}}
\end{aug}

\begin{abstract}
This paper introduces a comprehensive framework to adjust a discrete test statistic for improving its hypothesis testing procedure. The adjustment minimizes the Wasserstein distance to a null-approximating continuous distribution, tackling some fundamental challenges inherent in combining statistical significances derived from discrete distributions. The related theory justifies Lancaster's mid-p and mean-value chi-squared statistics for Fisher's combination as special cases. However, in order to counter the conservative nature of Lancaster's testing procedures, we propose an updated null-approximating distribution. It is achieved by further minimizing the Wasserstein distance to the adjusted statistics within a proper distribution family. Specifically, in the context of Fisher's combination, we propose an optimal gamma distribution as a substitute for the traditionally used chi-squared distribution. This new approach yields an asymptotically consistent test that significantly improves type I error control and enhances statistical power.
\end{abstract}



\end{frontmatter}

\section{Introduction}\label{sec:Intro}

\subsection{Challenges in testing discrete statistics and combining discrete p-values}

Combining the results of multiple hypothesis tests is a critical approach for global hypothesis testing, which finds applications in a wide range of fields including meta-analysis, signal detection, and other data-integrative studies. Many studies in the literature assume that the underlying statistics are continuous due to their homogeneous and straightforward mathematical structure. However, in reality, data and its corresponding statistics are often discrete. For example, genetic mutations can be rare, resulting in discrete genetic association test statistics \citep{Neale2011}.

The process of combining discrete tests and their corresponding p-values presents additional theoretical and computational challenges when compared to their continuous counterparts. One fundamental challenge lies in the fact that the exact test, based on the exact distribution of a discrete statistic, always results in a conservative test \citep{Berry}. To understand this point, it is helpful to recall that the p-value of any continuously distributed test statistic is always uniformly distributed between 0 and 1 under the null hypothesis. This uniform distribution endows the test with two preferred attributes: (i) the mean of the p-value is 0.5, which represents a neutral point between 0 and 1, and (ii) the type I error rate can be precisely controlled at any desired {\it size} $\alpha\in (0, 1)$, due to the continuity of the p-value \cite[Chapter 8]{Casella2002Statistical-inf}. However, these beneficial properties are not present in a discrete scenario. 

Specifically, consider a discrete test statistic $X\in \{x_i; i \in \mathbb{N}\}$ with a probability mass function $p_i = \P(X=x_i)$ and a cumulative distribution function (CDF) $F_i=\P(X\leq x_i) = \sum_{k=1}^{i}p_k$, $0=F_0<F_1< \cdots < F_i \cdots$, under the null hypothesis. Without loss of generality, consider the left-tailed p-value $P$, a random variable such that
\begin{equation}
	\label{eq:disc_P_left}
	P \in \{F_i, i \in \mathbb{N}\}, \text{ with } \P(P =F_i)=\P(X=x_i) =p_i = F_i - F_{i-1}. 
\end{equation}
The distribution of $P$ depends on the distribution of $X$ and is no longer uniformly distributed (not necessarily discrete uniform either). Moreover, under the null, its expected value $\E(P) = \frac{1}{2}\left(1+\sum_{i \in \mathbb{N}}p_i^2\right)> 0.5$, indicating a conservative bias towards 1 \citep{Berry}. Furthermore, the discreteness of $P$ means that the type I error rate cannot be precisely controlled at an exact $\alpha$ value (unless an $F_i$ value happens to be $\alpha$), but must be strictly less than $\alpha$. This constraint results in the {\it level}-$\alpha$ tests that are often stricter than desired, leading to conservative tests that sacrifice statistical power.

Additional challenges arise when combining discrete tests for global hypothesis testing. Suppose $F^{(j)}$ represents the null distribution of the $j$th individual statistic $X_j$. Consider the global null hypothesis
\begin{equation}
	\label{eq:globalNull}
	H_0 = \bigcap_{j=1}^n H_{0j}, \text{ composed of individual nulls } H_{0j}: X_j \sim F^{(j)}\quad j =1 \ldots n.
\end{equation}
To test $H_0$, a typical strategy is to combine the corresponding individual p-values $P_j$. For example, Fisher's combination test \citep{fisher1925statistical} applies the statistic
\begin{equation}
	\label{eq:fisherchi}
	T_n=-2\sum_{j=1}^n \log P_j.
\end{equation}
When the individual test statistics $X_j$'s are independent and continuous, $-2\log P_j \sim \chi^2_2$, implying that $T_n \sim \chi^2_{2n}$ under $H_0$. However, if the $X_j$'s are discrete, deriving the exact discrete distribution of $T_n$ requires an $n$-fold convolution, which becomes computationally demanding when $n$ is large \citep{mielke2004combining}. Resampling-based methods, such as permutation or simulation, are alternative approaches to estimating the distribution of $T_n$. However, these methods also demand heavy computation, especially in large-scale data analysis.

Various strategies have been put forth to counter the challenges mentioned above. One such strategy is Lancaster's mid-p adjustment, which modifies the discrete p-value in (\ref{eq:disc_P_left}) to be:
\begin{equation}
	\label{eq:midp_P}
	\tilde{P} = 
	\frac{F_i+F_{i-1}}{2} \text{ whenever } P=F_i. 
\end{equation}
Lancaster's premise rests on the fact that $\P(X< x_i)$ and $\P(X\leq x_i)$ are identical under continuity. Therefore, it's justifiable to define $\tilde{P}$ as the median of these two values in the discrete scenario. Moreover, \cite{barnard} (page 1477) suggests that any reasonable measure of significance should fall between $\P(X< x_i)$ and $\P(X\leq x_i)$.
Despite the fact that $\tilde{P}$ is unbiased with respect to Uniform(0,1) (since $\E(\tilde{P})=0.5$), its variance is strictly less than that of Uniform(0,1). This issue leads to concerns about its appropriateness for continuous approximations \citep{routledge1994practicing} and for p-value combination tests \citep{adjmidp}.

When combining discrete p-values by Fisher's statistic $T_n$ in \eqref{eq:fisherchi}, \cite{lancaster1949combination} proposed two adjustments for $-2\log (P)$ to bring it closer to $\chi^2_2$, i.e., its null distribution under continuity. Specifically, whenever $X=x_i$, or equivalently $P=F_i$, the mean-value-$\chi^2$ statistic (denoted by $\ddot{Z}$) and the median-value-$\chi^2$ statistic ($\tilde{Z}$) take the following values, respectively:
\begin{align}
	\ddot{Z}= \frac{\int_{F_{i-1}}^{F_i} (-2 \log (p)) dp}{F_i-F_{i-1}} = 2 - 2\frac{F_i\log (F_i)-F_{i-1}\log (F_{i-1})}{F_i-F_{i-1}},
	\label{eq:lmod_lexplicit}
\end{align}
\begin{equation}
	\tilde{Z} = -2\log(\tilde{P}) = - 2\log \left(\frac{F_i+F_{i-1}}{2}\right).
	\label{eq:medianchi}
\end{equation}
Put simply, $\ddot{Z}$ is the conditional average of the log-transformed p-value over intervals $[F_{i-1}, F_i]$, and $\tilde{Z}$ is the log-transformation of the mid-p. For combining $P_j$, $j=1, ..., n$, the $T_n$ statistic in \eqref{eq:fisherchi} is then adjusted to be, respectively,
\begin{equation}
	\label{eq:sum12}
	S_n = \sum_{j=1}^n \ddot{Z}_j,\quad \tilde{S}_n = \sum_{j=1}^n \tilde{Z}_j.
\end{equation}

\begin{figure}
	\begin{centering}
		\includegraphics[width=3.5in]{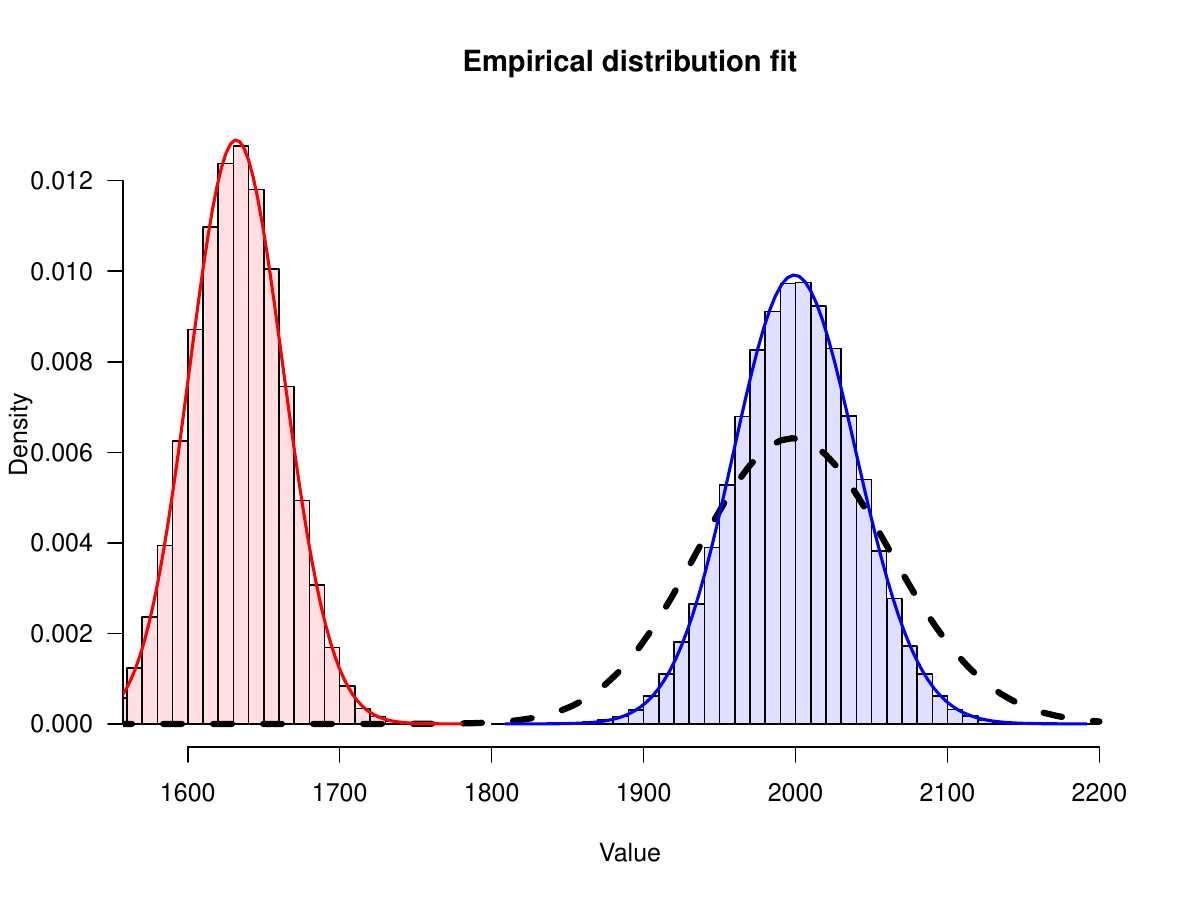}
		\caption{Densities of $\chi_{2n}^2$ (dashed curve) and optimal gamma distributions (solid curves) are compared with histograms of simulated statistics $S_n$ (right) and $\tilde{S}_n$ (left), which combine $n=1000$ independent and identically distributed p-values from $X_j \sim Binomial(5, 0.1)$. The optimal gamma distributions outperform $\chi_{2n}^2$ as approximations for the null distributions of $S_n$ and $\tilde{S}_n$.
		}
		\label{fig:histlanc_midp}
	\end{centering}
\end{figure}

\cite{lancaster1949combination} argued that both statistics in \eqref{eq:sum12} are closer to $\chi^2_{2n}$ than $T_n$ is, and that $\chi^2_{2n}$ can be used as a convenient null distribution for testing $H_0$. However, as we will prove, the mean of $\tilde{S}_n$ and the variances of $S_n$ and  $\tilde{S}_n$ are strictly smaller than the counterparts of $\chi^2_{2n}$, rendering $\chi^2_{2n}$ inappropriate for approximating their null distributions. Figure \ref{fig:histlanc_midp} displays the significant discrepancies between $\chi^2_{2n}$ (dashed curve) and the null distributions of $S_n$ and $\tilde{S}_n$ (histograms) in an illustrative example.

\subsection{Contributions of this work}

This work contributes to the theory and methodology of hypothesis testing for discrete statistics, especially in the context of Fisher's combination of independent discrete p-values. Our first contribution is the introduction of the Wasserstein metric into the development of optimal continuous approximations for hypothesis testing procedures involving discrete statistics.

Secondly, we establish a framework capable of adjusting any discrete statistic according to a given approximating continuous distribution by minimizing the Wasserstein distance. Specifically, given a discrete statistic $X$ and a target continuous distribution $Y$, we calculate a modified discrete statistic $Z$. Whenever $X = x_i$, we derive $Z = z_i$, ensuring that the Wasserstein distance between $Z$ and $Y$ is minimized. As the Wasserstein metric well describes the inherent geometry of the probability space \citep{Kantorovich2006}, $Z$ maintains the structure of $X$ while aligning as closely as possible with the distribution of $Y$.
In specific cases, if $X = P$ is a discrete p-value and its approximating continuous distribution is given by $Y \sim \text{Uniform}(0, 1)$, then the adjusted statistic is $Z=\tilde{P}$, the mid-p in \eqref{eq:midp_P}. Alternatively, if $X = -2\log(P)$ and $Y \sim \chi^2_2$, then $Z=\ddot{Z}$, the mean-value-$\chi^2$ in \eqref{eq:lmod_lexplicit}. These results provide a solid theoretical justification for Lancaster's adjustments. Furthermore, we describe the properties of the optimized statistic $Z$ in general.

Our third contribution is the proposal of an improved testing procedure designed to overcome the limitations of Lancaster's method. Broadly speaking, when testing $H_0$, we suggest a further adjustment of the approximating null distribution from $Y$ to an optimal continuous distribution, denoted as $\tilde{Y}$. This distribution minimizes the Wasserstein distance between the adjusted statistic $Z$ and a more general distribution family that includes $Y$. In the context of Fisher's combination test, we suggest modifying the approximating null distribution of statistics in \eqref{eq:lmod_lexplicit} -- \eqref{eq:sum12}. We demonstrate an improved fit when approximating their null distributions using gamma distributions that match the mean and variance of the corresponding statistics, rather than using chi-squared distributions as Lancaster proposed.
As depicted in Fig. \ref{fig:histlanc_midp}, the optimal gamma distributions (represented by solid curves) closely align with the true null distributions of $S_n$ and $\tilde{S}_n$ (the histograms). The corresponding testing procedures for $S_n$ and $\tilde{S}_n$ are asymptotically consistent as $n\to\infty$.
Compared to tests based on chi-squared distributions, our new procedure offers markedly more accurate type I error control and increased statistical power.

\section{Comparing discrete and continuous distributions with the Wasserstein metric}

\begin{figure}[h!]
	\centering
	\subfloat[$-2\log(P)$ vs. $\chi_2^2$]{
		\includegraphics[width=0.4\linewidth]{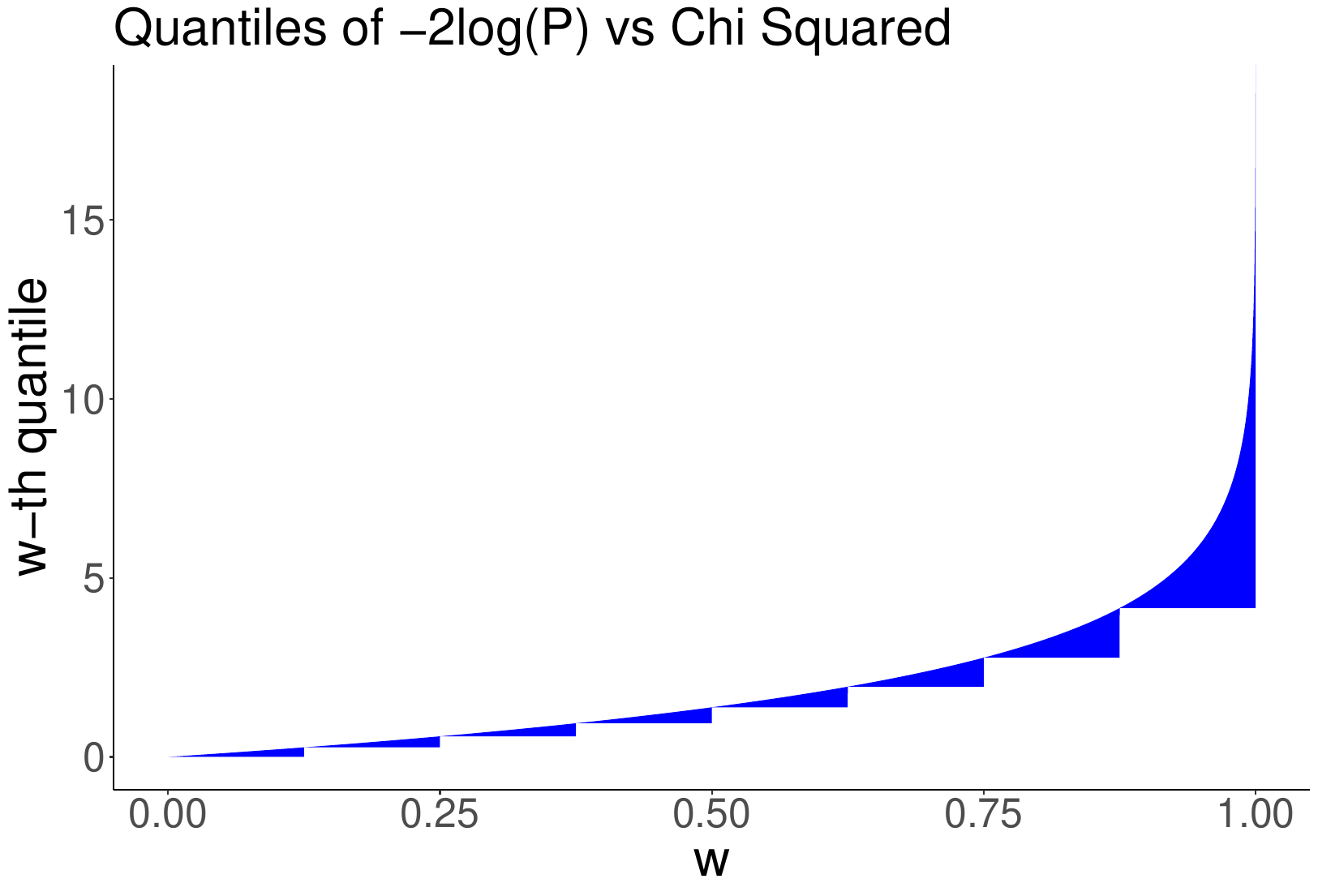}
	}
	\subfloat[$\ddot{Z}$ vs. $\chi_2^2$]{
		\includegraphics[width=0.4\linewidth]{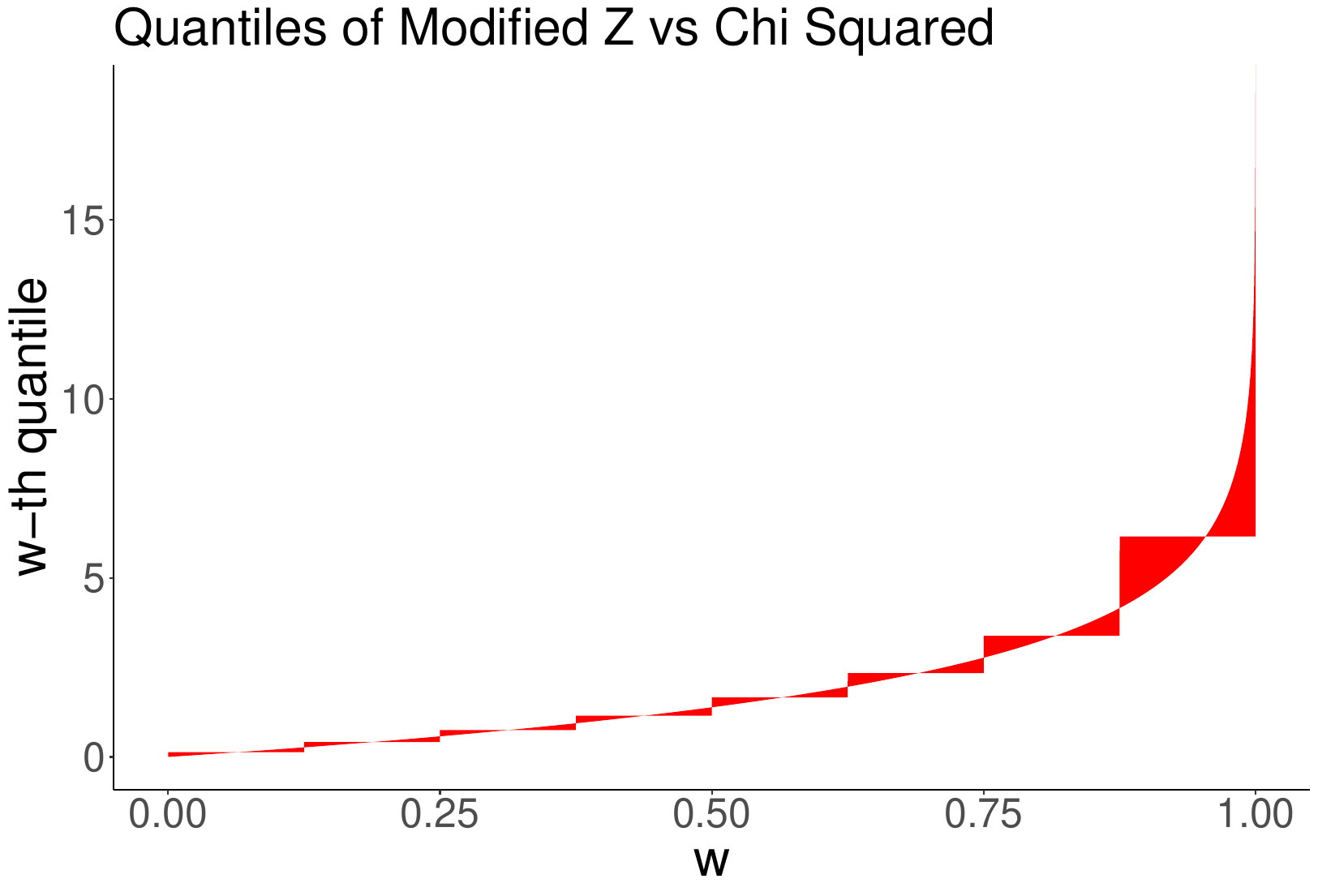}
	}
	\caption{$W_2^2(X,Y)$ in \eqref{eq:wass} is represented by the shaded area between the quantiles of $Y\sim\chi_2^2$ (smooth profile) and discrete $X$ (step profile). Left panel: $X=-2\log(P)$. Right: $X=\ddot{Z}$ is the mean-value-$\chi^2_2$ statistic in \eqref{eq:lmod_lexplicit} after adjustment. $P \in \{F_i = i/8, i=1, \ldots, 8\}$ with equal chance.
	}
	\label{fig:wdimprovement}
\end{figure}

In this section, we introduce the Wasserstein metric and its application in comparing discrete and continuous distributions. Specifically, given two probability measures $P$ and $Q$ in an Euclidean space $E$ and a constant $p>1$, the Wasserstein distance $W_p(P,Q)$ between $P$ and $Q$ measures the ``cost" of reallocating the distribution $P$ into the distribution $Q$  \citep{villani2003topics}: 
\begin{equation}
	W_p(P,Q)=\left(\inf_{\gamma \in C(P,Q)}\int_{E^2} |x-y|^pd\gamma (x,y)\right)^\frac{1}{p}, 
	\label{eq:wassgen}
\end{equation}
where $C(P, Q)$ is the set of all couplings of $P$ and $Q$, i.e., the set of all joint probability distributions in the product space $E^2$ with the marginals $P$ and $Q$. 
For convenience, we use the notations $W_p(P,Q)$ and $W_p(X,Y)$ interchangeably, where $X$ and $Y$ are random variables corresponding to $P$ and $Q$, respectively. 
When $X$ and $Y$ are real-valued one-dimensional random variables, an explicit formulation of the Wasserstein distance is  \cite[Chapter 5]{shao2006}
\begin{equation}
	W_p^p(X,Y)=\int_0^1|F^{-1}(w)-G^{-1}(w)|^pdw, 
	\label{eq:wass}
\end{equation}
where $F$ and $G$ are the respective CDF of $X$ and $Y$, and $F^{-1}$ and $G^{-1}$ denote the quantile functions, e.g., $F^{-1}(w)=\inf \{x: F(x) \geq w\}$. At $p=2$, $W_2^2(X,Y)$ represents the area between the quantile curves of $X$ and $Y$. Figure \ref{fig:wdimprovement} illustrates $W_2^2(X, Y)$ when $Y \sim \chi^2_2$ and $X$ is either the log-transformed discrete p-value $-2\log(P)$ (left panel) or the mean-value-$\chi^2_2$ statistic in \eqref{eq:lmod_lexplicit} (right panel). By adjusting $X$ to minimize its Wasserstein distance to $Y$, we can minimize the area of discrepancy between the two distributions.

The Wasserstein distance between two random variables is determined solely by their marginal distributions. However, the minimizer coupling in \eqref{eq:wassgen} is crucial for characterizing the joint distribution of a continuous distribution and its closest discrete distribution. For that, we construct a measure in Definition \ref{def:coupling} and prove that it is a suitable minimizer coupling through Theorem \ref{thm:wpdc}.

\begin{definition}
	\label{def:coupling}
	Let $X$ be a discrete random variable with increasing support and CDF $F$ as defined in \eqref{eq:disc_P_left}. Let $Y$ be a continuous random variable with invertible CDF $G$ and density $g$. For a set $A$, let $\delta_A(x)$ be the indicator function of the event $x\in A$. We define a probability measure $\pi$ as the unique extension of the function
	\begin{align*}
		\pi : \mathbb{R}^2 \rightarrow [0,1] \text{ and }
		\pi(A\times B) = \sum_{i \in \mathbb{N}} \delta_A(x_i)\int_{B\cap [G^{-1}(F_{i-1}), G^{-1}(F_i)]} g(y)dy
	\end{align*}
	to all Borel sets.  
\end{definition}

\begin{theorem}
	\label{thm:wpdc}
	The measure $\pi$ defined in Definition \ref{def:coupling}  is a coupling of the probability distributions of $X$ and $Y$, and it satisfies
	$$W_p^p (X,Y)= \int_{\mathbb{R}^2}|x-y|^p \pi(dx,dy).$$
	In other words, $\pi$ achieves the infimum among all couplings in the definition of \eqref{eq:wassgen}.
\end{theorem}

At $p=2$, Corollary \ref{cor:wdmoments} gives a convenient formula to calculate the distance between discrete and continuous distributions using their moments and covariance. 

\begin{corollary}
	\label{cor:wdmoments}
	Assuming the same notation and assumptions as in Theorem \ref{thm:wpdc}, we have
	$$W_2^2 (X,Y)= \Var(X)+\Var(Y)+(\E(X)-\E(Y))^2-2\Cov_\pi(X,Y),$$
	where $\Cov_\pi$ denotes the covariance with respect to the measure $\pi$ in Definition \eqref{def:coupling}:
	$$\Cov_\pi(X,Y)=\sum_{i \in \mathbb{N}} x_i\int_{G^{-1}(F_{i-1})}^{G^{-1}(F_i)} yg(y) dy-\E(X)\E(Y).$$
\end{corollary}

\section{Optimal adjustment of discrete statistics}

\subsection{A general framework}

We present a general framework for adjusting a discrete statistic $X$ based on a suitable continuous ``target" distribution represented by $Y$. In practice, one can choose $Y$ based on a continuous approximation of $X$. For instance, if $X=P$ is a discrete p-value, a reasonable choice is $Y \sim \text{Uniform}(0,1)$; if $X=-2\log(P)$, a choice is $Y \sim \chi^2_2$ \citep{lancaster1949combination}.

The goal is to find a discrete statistic $Z$ that takes values $z_i$ whenever $X=x_i$, for each $i$, with $z_i$ values being chosen so that the distribution of $Z$ is closest to $Y$ in Wasserstein metric. Let $supp(V)$ denote the support of a random variable $V$. We have: 
\begin{align}
	Z &= \mbox{argmin}_{V \in \mathcal{M}_X}\{W_p(V, Y)\}, \mbox{where } p>1 \mbox{ and}  \label{eq:zopti}\\
	\mathcal{M}_X &= \{V: supp(V)= \{v_i\}_{i \in \mathbb{N}} \mbox{ with }  \P(V=v_i)=\P(X=x_i) \text{ for all } i \in \mathbb{N}\}.  
	\label{eq:candidates}
\end{align}

This optimization procedure ``collapses" the density of $Y$ onto $Z$ while preserving the probability masses of $X$. To see this, 
using the same notations in Definition \ref{def:coupling}, the domain of $Y$ can be partitioned based on the distribution of $X$ as follows:
\begin{equation}
	A_i=[G^{-1}(F_{i-1}),G^{-1}(F_i)], \quad i \in \mathbb{N}.
	\label{eq:partition}
\end{equation}
We have $\P(Y \in A_i)=F_i-F_{i-1}= \P(X=x_i)=\P(Z=z_i)$. Moreover, Lemma \ref{lem:newm} below shows that $z_i \in A_i$ for all $i$. Therefore, the Wasserstein metric is an excellent tool for adjusting $X$ because it provides a one-to-one and monotone correspondence between observation $X=x_i$ and value $z_i$ in the hypothesis testing procedure. The statistic $Z$ is closest to $Y$ in distribution among all discrete random variables that naturally retain the probabilistic structure of $X$. 
\begin{lemma}
	\label{lem:newm}
	The minimizer $Z$ obtained by \eqref{eq:zopti} and \eqref{eq:candidates} belongs to $\mathcal{M}_X^* \subset \mathcal{M}_X$, where
	\begin{align*}
		\mathcal{M}_X^*=\{V: supp(V)=\{v_i\}_{i \in \mathbb{N}} \text{ with }  v_i \in A_i \mbox{ and } \P(V=v_i)=\P(X=x_i) \text{ for all }  i \in \mathbb{N}\}.
	\end{align*}
\end{lemma}

The values of $Z$ are given by Theorem \ref{thm:zopti}.

\begin{theorem}
	Consider a discrete random variable $X$ and a continuous random variable $Y$ with domain partition in (\ref{eq:partition}). For a fixed $p>1$, the $Z$ variable obtained by \eqref{eq:zopti} and \eqref{eq:candidates} has support $supp(Z) = \{z_i\}_{i \in \mathbb{N}}$ with
	$$z_i = \left(\E\left[Y^{p-1} | Y \in A_i\right]\right)^{\frac{1}{p-1}}.$$
	In particular, when $p=2$, we have $z_i=\E \left[Y | Y \in A_i \right]$.
	\label{thm:zopti}
\end{theorem}

According to Theorem \ref{thm:zopti},  $z_i^{p-1}$ is the average of $Y^{p-1}$ on the projected quantiles of $X$. Furthermore, we have $\E(Y^{p-1}|Z) = Z^{p-1}$, which in turn implies $\Cov(Y^{p-1},Z^{p-1})=\Var(Z^{p-1})$.  
When $p=2$, $z_i$ is simply the average of $Y$ in $A_i$, which means that the mid-p and mean-value-$\chi^2$ are both Wasserstein-optimal at $p=2$, as deduced and formally stated in Corollary \ref{cor:zopti2}.

\begin{corollary}
	If $X=P$ is the discrete p-value and $Y\sim \text{Uniform}(0,1)$, $Z=\tilde{P}$, the mid-p in \eqref{eq:midp_P}. If $X=-2\log(P)$ and $Y \sim \chi^2_2$, $Z=\ddot{Z}$, the mean-value-$\chi^2$ in \eqref{eq:lmod_lexplicit}.
	\label{cor:zopti2}
\end{corollary}

In the following, we will delve into the the properties of the optimally adjusted statistics. According to Corollary \ref{cor:zopti2}, the mid-p and mean-value-$\chi^2$ statistics are unbiased, meaning their means are equal to the means of their target distributions, i.e., $\E(\tilde{P})=0.5$ and $\E(\ddot{Z})=2$. However, it's important to note that their variances are smaller than their target distributions:
\begin{align*}
	\Var(\tilde{P}) &=\sum_{i \in \mathbb{N}} F_iF_{i-1}(F_i-F_{i-1})/4 < 1/12 =\Var(\text{Uniform}(0, 1)),\\
	\Var(\ddot{Z}) &= 4\sum_{i \in \mathbb{N}} ( F_i\log(F_i) -F_{i-1}\log(F_{i-1}))^2/(F_i-F_{i-1})< 4 = \Var(\chi^2_2).
\end{align*}

The comparison of the means and variances between the target distributions and the optimally adjusted statistics can be generalized 
for any $p>1$, as detailed in Corollary \ref{cor:unbiased}.


\begin{corollary} 
	\label{cor:unbiased}
	Under the same assumptions and notations as in Theorem \ref{thm:zopti}, we have the following: $\E(Z) <  \E(Y)$ for $p \in (1,2)$; $\E (Z) = \E(Y)$ for $p=2$; and $\E(Z)  > \E(Y)$ for $p >2$. Regarding the variance, $\Var (Z^{p-1}) <  \Var(Y^{p-1})$ for all $p>1$. 
	As a special case, when $p=2$ and using the domain partition in \eqref{eq:partition}, we have
	$$\Var(Z)=\Var(Y)-\sum_{i \in \mathbb{N}}\left(E(Y^2\big|A_i) -z_i^2\right)(F_i-F_{i-1}).$$
\end{corollary}

Corollary \ref{cor:unbiased} demonstrates that the Wasserstein-optimal adjustment $Z$ is unbiased with respect to $Y$ only when $p=2$. Therefore, $p=2$ is an appropriate choice for applying the Wasserstein metric to adjust a discrete statistic, in which case the variance of $Z$ is always smaller than the variance of $Y$. Interestingly, Corollary \ref{cor:wdmoments2} shows that the Wasserstein distance between a discrete statistic and its target continuous distribution equals the difference of their variances. 

\begin{corollary}
	\label{cor:wdmoments2}
	Under the same assumptions and notations as in Theorem \ref{thm:zopti} with $p=2$,
	$$W_2^2(Z,Y)=\Var(Y)-\Var(Z)=\E (\Var(Y|X)).$$
\end{corollary}

In particular, $W_2^2(\ddot{Z},\chi_2^2) = 4 - \Var(\ddot{Z})$ and $W_2^2(\tilde{P}, \text{Uniform}(0,1)) =1/12-\Var(\tilde{P})$. 
The equality $W_2^2(Z,Y)=\E (\Var(Y|X))$ indicates that the Wasserstein distance between $Y$ and $Z$ equals the average variance of $Y$ that is left unexplained by $X$. Moreover, as $Z$ is a discrete analogue of $Y$, the distance between them is attributed to the probability masses of $Z$ being concentrated at discrete points, while $Y$ possesses a spread-out density. As the grid of probability masses $\{F_i -F_{i-1}\}_{i \in \mathbb{N}}$ becomes finer, $Z$ approximates $Y$ more closely and so does its variance, leading to a decrease in distance. Corollary \ref{cor:closetc} provides an asymptotic result to describe this phenomenon.

\begin{corollary}
	Assuming the same conditions and notation as in Theorem \ref{thm:zopti}, when $p=2$, as $\Delta \equiv \sup_i \{F_i -F_{i-1}\}_{i \in \mathbb{N}} \rightarrow 0$, we have $\Var (Z) \rightarrow \Var (Y)$ and $W_2(Z, Y)\to 0$.
	\label{cor:closetc}
\end{corollary}

\subsection{Comparing the median-value and the mean-value $\chi^2$ statistics}

In the following, we discuss the properties of $\tilde{Z}$, the median-value-$\chi^2$ in \eqref{eq:medianchi}, in comparison with the mean-value-$\chi^2$ statistic $\ddot{Z}$. To simplify notations, we will denote $m = \E(\tilde{Z})$, $v = \Var(\tilde{Z})$, and $\nu = \Var(\ddot{Z})$ throughout the remainder of the paper.

\cite{lancaster1949combination} suggested using $\chi^2_2$ as the approximate null distribution for both statistics in hypothesis testing. As previously discussed, $\ddot{Z}$ is unbiased to $\chi^2_2$ since $\E(\ddot{Z}) = 2$, but $\tilde{Z}$ is biased with $m < 2$. Furthermore, $\tilde{Z}$ has a smaller variance: $v<\E(\tilde{Z} -2)^2 =4\sum_{i \in \mathbb{N}} {(\log ((F_i+F_{i-1})/2) +1)(F_i-F_{i-1})}^2/(F_i-F_{i-1}) \approx \nu < 4$ (the approximation is by the intermediate value theorem). These comparisons are expected because $Z$ is shown to be the distance minimizer to $\chi^2_2$ while $\tilde{Z}$ is not. Therefore, according to Corollary \ref{cor:zopti2}, we have $W_2^2(\tilde{Z},\chi_2^2) \geq W_2^2(\ddot{Z},\chi_2^2)$. In fact, direct calculations show that $W_2^2(\tilde{Z},\chi_2^2) =W_2^2(\ddot{Z},\chi_2^2)+\E[(\ddot{Z} -\tilde{Z})^2]$. This calculation results in Corollary \ref{cor:closetochi1}, which confirms that the inequality is strict.

\begin{corollary}\label{cor:closetochi1}
	For any discrete p-value $P$, it holds that $W_2^2(\tilde{Z},\chi_2^2) > W_2^2(\ddot{Z},\chi_2^2).$
\end{corollary}

The $\ddot{Z}$ and $\tilde{Z}$ statistics represent two distinct strategies for adjustment and transformation. $\ddot{Z}$ is the optimal adjustment for $-2\log P$, the log-transformation of the original discrete p-value $P$. Conversely, $\tilde{Z}=-2\log \tilde{P}$ is the log-transformation of the mid-p $\tilde{P}$, which is the optimal adjustment for $P$.
A more fundamental and broader question arises: should we opt to adjust the statistic through minimum distance before or after its transformation? Corollary \ref{cor:convex} offers a perspective to this question by comparing the relative magnitudes of these statistics under any convex monotone transformation.

\begin{corollary}
	\label{cor:convex}
	Assume that $T$ is a strictly monotonic and convex transformation, and let $p=2$. Denote $Z$ as the minimizer of \eqref{eq:zopti} and \eqref{eq:candidates} for a discrete-continuous pair $(X,Y)$, and $Z^T$ as the minimizer for the pair $(T(X),T(Y))$. It follows that $T(Z)< Z^T$ almost surely and that $\E(T(Z)) < \E(Z^T) =\E(T(Y))$.
\end{corollary}

Corresponding to our discussion, $X = P$ and $Y\sim$ Uniform(0,1), so the adjusted statistic $Z=\tilde{P}$ is the mid-p in \eqref{eq:midp_P}, and $T(Z) = \tilde{Z}$. To adjust after transformation, $T(X) = T(P) = -2\log P$, $T(Y)\sim \chi^2_2$, and $Z^T = \ddot{Z}$. Because  $\E(\tilde{Z}_j) < \E(\ddot{Z}_j) =2$, after combining these statistics by $S_n$ and $\tilde{S}_n$ in \eqref{eq:sum12}, we have $\E(\tilde{S}_n) < \E(S_n) = 2n$. Therefore, for Lancaster's testing procedure, where $\chi_{2n}^2$ is used as the null distribution, the unbiased $S_n$ (adjusting the statistic after transformation) is a better choice than the biased $\tilde{S}_n$ (adjusting the statistic before transformation). Similar relationships holds for any strictly monotonic and convex transformation function. 

\section{An improved testing procedure}
\subsection{Further adjustment of the null distribution}

According to Corollary \ref{cor:unbiased}, the $W_2$-optimally adjusted statistic $Z$ always has a smaller variance than the target distribution $Y$. If this target distribution is used as the null distribution in a testing procedure, it would lead to a conservative test. To rectify this issue, we propose a second-stage adjustment: the adjustment of $Y$ to a distribution $\tilde{Y}$ that is closer to $Z$ in Wasserstein distance. 

The optimal $\tilde{Y}$ should be achieved within a distribution family that satisfies two properties. Firstly, it should contain the $Y$ distribution to preserve the desired properties of $Y$ (e.g., being consistent with the continuous scenario). Secondly, when combining any number of discrete p-values in the form of \eqref{eq:sum12}, the distribution family should be additive. \cite{lancaster1949combination} leveraged the additive property of the chi-squared distribution by the same reasoning. A family of probability distributions $\mathcal{F}$ is said to be additive if $X_1+X_2\in \mathcal{F}$ whenever $X_1\in \mathcal{F}$ and $X_2\in \mathcal{F}$.

In the case of Fisher's combination type tests, we suggest the gamma family, which contains the chi-squared distributions and is additive. More specifically, among the $Gamma(\alpha,\beta)$ distributions, with the shape parameter $\alpha >0$ and the scale parameter $\beta>0$, we search for $\tilde{Y}$ such that $W_2^2(Z,\tilde{Y})$ is minimal. Theorem \ref{thm:optgamma} shows that the gamma distribution closest to any positive discrete random variable is the one that matches its mean and variance.

\begin{theorem}\label{thm:optgamma}
	Let $X$ be any non-negative discrete random variable with mean $\mu$ and variance $\sigma^2$, and $\mathcal{F}=\{Y \sim Gamma(\alpha,\beta): \alpha > 0, \beta>0\}$ be the family of gamma distributed random variables. Then,
	$$W^2_2(X, \mathcal{F})\equiv \inf_{Y \in \mathcal{F}}W^2_2(X,Y)=W^2_2(X,\tilde{Y}),$$
	where $\tilde{Y}\sim Gamma(\tilde{\alpha},\tilde{\beta})$ with $\tilde{\alpha}=\mu^2/\sigma^2$ and $\tilde{\beta}= \sigma^2/\mu$, i.e., $\E(\tilde{Y})=\E(X)$ and $\Var(\tilde{Y})=\Var(X)$.
\end{theorem}

For the mean-value-$\chi^2$ statistic $\ddot{Z}$ with $\E(\ddot{Z})=2$ and $\Var(\ddot{Z})=\nu$, we obtain the adjusted distribution $\tilde{Y} \sim Gamma(\tilde{\alpha},\tilde{\beta})$, where $\tilde{\alpha}=4/\nu$ and $\tilde{\beta}= \nu/2$. Unlike the $\chi^2_2$ distribution used in \cite{lancaster1949combination}, $\tilde{Y}$ shares the same mean and variance of $\ddot{Z}$, which makes the corresponding hypothesis testing procedure better at controlling the type I error rate.
Similarly, for the median-value-$\chi^2$ statistic $\tilde{Z}$ in \eqref{eq:medianchi}, with $\E(\tilde{Z})=m$ and $ \Var(\tilde{Z})=v$, we obtain $\tilde{Y} \sim Gamma(m^2/v, v/m)$, which provides a much improved null distribution over $\chi^2_2$ for $\tilde{Z}$. 


\subsection{Combining independent and identically distributed p-values}

Following the above discussion, since $S_n$ and $\tilde{S}_n$ statistics combine $\ddot{Z}_j$ and  $\tilde{Z}_j$, respectively, they are better approximated by suitable gamma distributions rather than $\chi^2_{2n}$. Specifically, when the discrete p-values $P_j$ are independent and identically distributed (i.i.d.), the additivity property of the gamma distribution yields the moment-matched approximations:
\begin{equation}
	\label{eq:approxSn}
	S_n \overset{D}{\approx} \text{Gamma} (4n/\nu, \nu/2) \quad \text{and} \quad \tilde{S}_n \overset{D}{\approx} \text{Gamma} (nm^2/v, v/m).
\end{equation}

These approximations guarantee asymptotically consistent tests as $n\to\infty$. Let $G(x; \alpha,\beta)$ denote the CDF of the gamma distribution, and let $q_{p; \alpha, \beta}=G^{-1}(p; \alpha, \beta)
=\inf\{x: G(x;\alpha,\beta)\geq p\}$ be the $p$-quantile for any $p\in (0, 1)$.

\begin{proposition} \label{thm:asymptoticpvalue}
	Let $\{X_i\}_{i \in \mathbb{N}}$ be a sequence of i.i.d. non-negative discrete random variables. Define $\alpha=[\E(X_1)]^2/\Var (X_1)$, $ \beta=\Var(X_1)/\E(X_1)$ and $T_n=\sum_{i=1}^nX_i$. For any $p \in (0,1)$, 
	$$
	\lim_{n \rightarrow \infty} \P(T_n < q_{p; n\alpha, \beta})=p.
	$$
	Furthermore, $G(T_n; n\alpha,\beta)$ converges in distribution to Uniform(0,1).
\end{proposition}

By Proposition \ref{thm:asymptoticpvalue}, for $S_n$ and $\tilde{S}_n$ in \eqref{eq:sum12} combining i.i.d. p-values and any $p \in (0,1)$, we have
$$
\lim_{n \rightarrow \infty} \P(S_n < q_{p; 4n/\nu, \nu/2})=p,\quad
\lim_{n \rightarrow \infty} \P(\tilde{S}_n < q_{p; nm^2/v, v/m})=p.
$$
This means $\P\left(S_n \geq q_{1-\alpha; 4n/\nu,\nu/2}\right) \to \alpha$, indicating that we can control the type I error rate well at any significance level $\alpha \in (0, 1)$ when $n$ is large. In contrast, using $\chi^2_{2n}$ as the null distribution will lead to systematic bias, since $\P\left(S_n \geq q_{1-\alpha; n, 2} \right) = \Phi\left(4\Phi^{-1}(1-\alpha)/\nu+ O(\sqrt{n})\right) \to 1$. 


Moreover, both $G(S_n; 4n/\nu, \nu/2)$ and $G(\tilde{S}_n; nm^2/v,v/m)$ converge in distribution to Uniform(0,1).
This property allows for approximating any bounded functions of the CDFs, in the sense that $\E(h(G(S_n; 4n/\nu, \nu/2))) \to \int_0^1 h(u)du$ for any bounded function $h$, based on the Portmanteau Lemma (cf. \cite{van2000asymptotic,Bill86}). In the Appendix, further discussion on these asymptotic convergences and their connection to the Wasserstein distance can be found after the proof of Proposition \ref{thm:asymptoticpvalue}.




\subsection{Combining independent, non-identically distributed p-values}

Consider a sequence of discrete p-values $\{P_j\}_{j \in \{1, \ldots n\}}$ obtained from independent but non-identically distributed statistics $X_j$ with cumulative distributions $F^{(j)}$, $j = 1, \ldots, n$. Let $\ddot{Z}_j$ and $\tilde{Z}_j$ be the corresponding statistics in \eqref{eq:lmod_lexplicit} and \eqref{eq:medianchi}, respectively, with parameters given by 
\begin{equation}
	\label{eq:paramnoniid}
	\nu_j=\Var(\ddot{Z}_j),\quad m_j=\E(\tilde{Z}_j),\quad v_j=\Var(\tilde{Z}_j).
\end{equation}
As we have shown, the optimal gamma approximations for $\ddot{Z}_j$ and $\tilde{Z}_j$ are $Gamma(4/\nu_j, \nu_j/2)$ and $Gamma(m_j^2/v_j, v_j/m_j)$, respectively.  The density function of the sum of independent but non-identically distributed gamma random variables is given by the digamma density, as derived by \cite{Mo1985}. However, exact computation of the density function is challenging as it involves an infinite sum, which requires further approximation steps to calculate quantiles.


We propose to approximate $S_n$ and $\tilde{S}_n$ using the optimal gamma distributions obtained through moment matching. Let $\bar{\nu}$, $\bar{m}$, and $\bar{v}$ be the averages of the moments in \eqref{eq:paramnoniid} over $j \in \{1, \ldots n\}$. We have $\E(S_n)=2n$, $\Var(S_n)=n\bar{\nu}$, $\E(\tilde{S}_n)=n\bar{m}$, and $\Var(\tilde{S}_n)=n\bar{v}$. According to Theorem \ref{thm:optgamma}, the best approximations among all possible gamma distributions are given by
\begin{equation}
	\label{eq:approxnoniid}
	S_n \overset{D}{\approx} Gamma(4n/\bar{\nu},\bar{\nu}/2),\quad \tilde{S}_n \overset{D}{\approx} Gamma(n\bar{m}^2/\bar{v},\bar{v}/\bar{m}).
\end{equation}

These approximations ensure that the type I error can be consistently controlled as $n\to \infty$ according to Proposition \ref{thm:asymptoticpvaluenoniid}. Proposition \ref{thm:asymptoticpvaluenoniid} extends the result of Proposition \ref{thm:asymptoticpvalue} from the i.i.d. scenario to the non-identically distributed scenario, under mild regularity conditions. 

\begin{proposition} \label{thm:asymptoticpvaluenoniid}
	For independent $\{\ddot{Z}_j\}_{j \in \mathbb{N}}$ and $\{\tilde{Z}_j\}_{j \in \mathbb{N}}$, if the parameter sequences in \eqref{eq:paramnoniid} satisfy the classic Lyapunov condition, i.e., for some $\delta, \tilde{\delta} >0$,
	$$
	\lim_{n \to \infty} \frac{\sum_{j=1}^n\E(|\ddot{Z}_j-2|^{2+\delta})}{(\sum_{j=1}^n \nu_j )^{1+\delta/2}}=0 
	\quad \text{and} \quad 
	\lim_{n \to \infty} \frac{\sum_{j=1}^n\E(|\tilde{Z}_j-m_j|^{2+\tilde{\delta}})}{(\sum_{j=1}^n v_j )^{1+\tilde{\delta}/2}}=0,
	$$
	then, for any $p \in (0,1)$,
	$$
	\lim_{n \rightarrow \infty} \P(S_n < q_{p; 4n/\bar{\nu},\bar{\nu}/2})=p
	\quad \text{and} \quad 
	\lim_{n \rightarrow \infty} \P(\tilde{S}_n < q_{p; n\bar{m}^2/\bar{v},\bar{v}/\bar{m}})=p.
	$$
\end{proposition}

A remark regarding the conditions in Proposition \ref{thm:asymptoticpvaluenoniid} is provided after its proof in Appendix section \ref{sec:proofs}. With this remark, we illustrate that these conditions can be easily satisfied in common scenarios, such as when the p-values are i.i.d. or when combining non-identical binomial p-values with proportions that are bounded away from 0 and 1.

\section{Numerical Studies}\label{sec:numerial}

\subsection{Combining independent and identically distributed p-values}

In this section, we conduct numerical studies to compare the performance of approximating the null distributions of $S_n$ and $\tilde{S}_n$ in (\ref{eq:sum12}) using gamma distributions as in \eqref{eq:approxSn} versus using $\chi^2_{2n}$ as suggested by \cite{lancaster1949combination}. The comparisons are made in three aspects: the fit to the empirical distributions of $S_n$ and $\tilde{S}_n$, the type I error, and statistical power.

Under the i.i.d. scenario, we consider p-values derived from i.i.d. discrete test statistics $X_j \sim \text{Binomial}(K,\theta)$, for $j=1,\ldots,n$. The aim is to test the global hypotheses: 
\begin{equation}
	\label{eq:H0A}
	H_0:\theta=\theta_0 \text{ vs. } H_A:\theta < \theta_0, 
\end{equation}
where we assume a left-sided alternative without loss of generality. Table ~\ref{table:paramfit} lists various choices of $\theta_0$ and $K$ values and the corresponding moments of the mean-value-$\chi^2$ statistic $\ddot{Z}_j$ in \eqref{eq:lmod_lexplicit} and the median-value-$\chi^2$ statistic $\tilde{Z}_j$ in \eqref{eq:medianchi}. Recall that $\E(\ddot{Z}_j)=2$ for all $j$.

The table indicates that $\tilde{Z}_j$ is biased relative to $\chi^2_2$, with $m = \E(\tilde{Z}_j)< 2$, aligning with Corollary \ref{cor:convex}. Although $\ddot{Z}_j$ is unbiased in relation to $\chi^2_2$, its variance is always smaller, that is, $\nu < 4$, consistent with Corollary \ref{cor:unbiased}. Moreover, the table confirms that $\nu= \Var(\ddot{Z}_j)$ is larger than $v = \Var(\tilde{Z}_j)$. Consequently, if $\chi^2_{2n}$ were employed for Fisher's combination test, as recommended by \cite{lancaster1949combination}, $S_n$ would be preferable to $\tilde{S}_n$ for integrating discrete p-values. Furthermore, as $K$ increases, particularly when $\theta_0$ is away from 0 or 1, the moments of $\ddot{Z}_j$ and $\tilde{Z}_j$ become closer to those of $\chi_2^2$, in line with Corollary \ref{cor:closetc}.

\begin{table*}
		\caption{Parameters of Binomial$(K,\theta_0)$ and the corresponding moments of $\ddot{Z}_j$ and $\tilde{Z}_j$}
		\label{table:paramfit}	
		\begin{tabular}{l c c c c c c c c c}
			\hline
			$K$        & \multicolumn{3}{c}{5}   & \multicolumn{3}{c}{10}   & \multicolumn{3}{c}{20} \\  \hline
			$\theta_0$ & 0.01 & 0.1   & 0.5  & 0.01 & 0.1  & 0.5  & 0.01 & 0.1  & 0.5  \\ \hline
			$\nu= \Var(\ddot{Z}_j)$     & 0.20 & 1.61   & 3.61 & 0.38 & 2.53 & 3.83 & 0.73 & 3.37 & 3.92 \\ 
			$m= \E(\tilde{Z}_j)$      & 1.41 & 1.63 & 1.92 & 1.44 & 1.77 & 1.96 & 1.50  & 1.89 & 1.98 \\ 
			$v=\Var(\tilde{Z}_j)$      & 0.10  & 0.96 & 3.19 & 0.19 & 1.74 & 3.63 & 0.38 & 2.74 & 3.81 \\ \hline
	\end{tabular}
\end{table*}

\begin{figure}
	\begin{centering}
		{\includegraphics[width=0.49\linewidth]{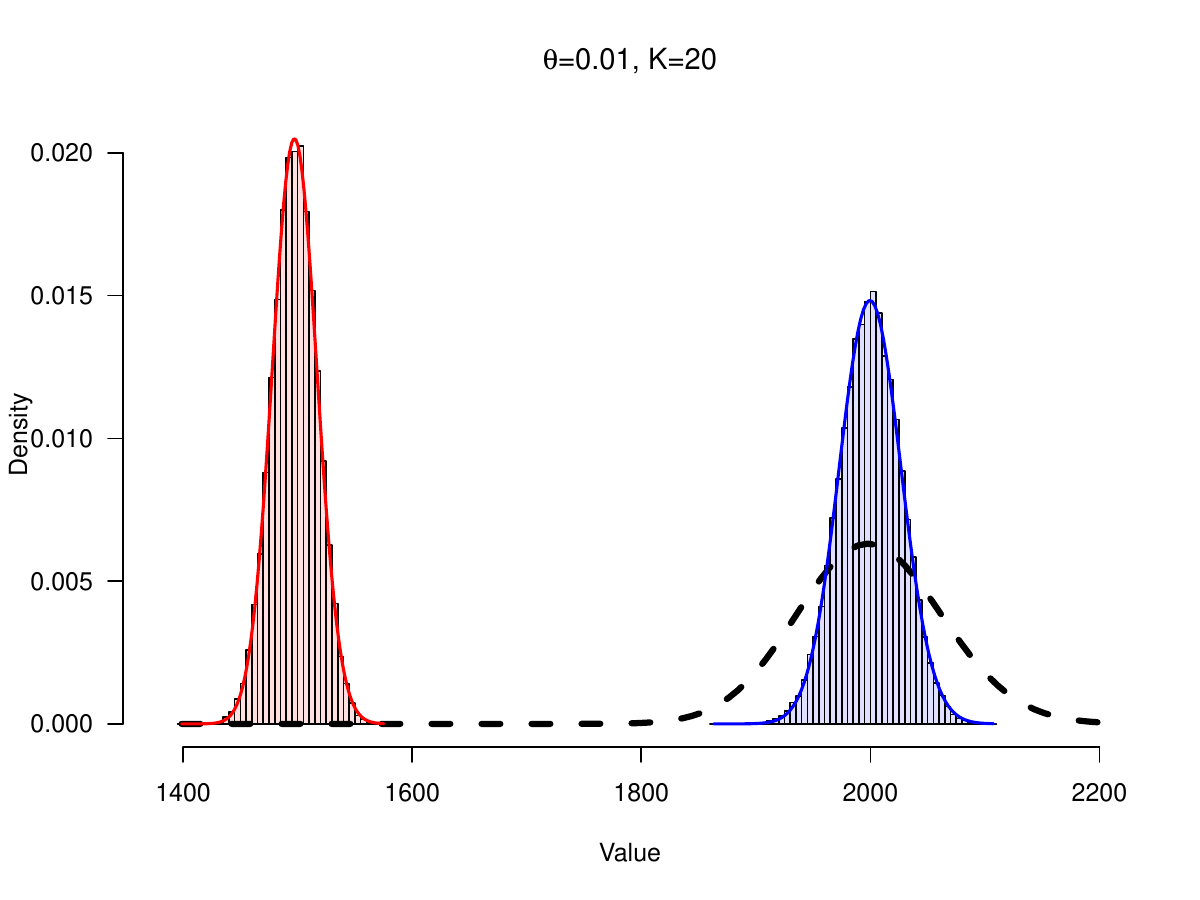}}
		{\includegraphics[width=0.49\linewidth]{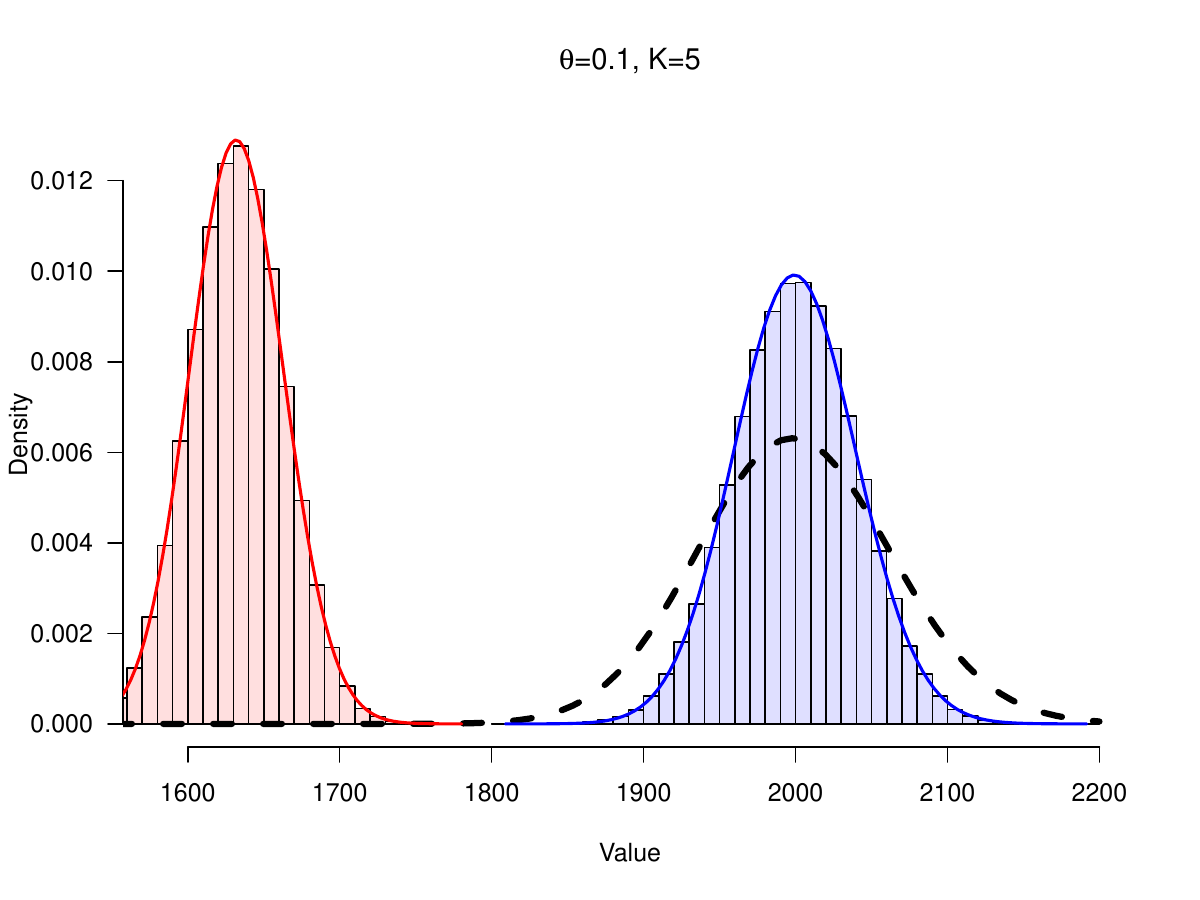}}
		\caption{
			The densities of $\chi_{2n}^2$ (illustrated by the dashed curve) and the optimal gamma distributions (solid curves) are contrasted with histograms of simulated statistics $S_n$ (right) and $\tilde{S}_n$ (left). These statistics combine $n=1000$ i.i.d. p-values from $X_j \sim Binomial(K, \theta_0)$. The gamma distributions clearly outperform $\chi_{2n}^2$ in fitting the null distributions of $S_n$ and $\tilde{S}_n$.	}
		\label{fig:histlanc_midp_more}	
	\end{centering}
\end{figure}

We assess the goodness-of-fit of the $\chi^2_{2n}$ distribution and the optimal gamma distributions to the empirical distributions of $S_n$ and $\tilde{S}_n$ using $10^7$ simulated statistics. Figure \ref{fig:histlanc_midp_more} displays the results for combining $n=1000$ left-sided binomial p-values at $(\theta_0,K)=(0.1,5)$ and $(0.01,20)$.
The figure illustrates that the variance of $S_n$ is significantly smaller than that of $\chi_{2n}^2$ when $K$ and $\theta$ are small. The distribution of $S_n$ is much better approximated by the gamma distribution. Similarly, $\tilde{S}_n$ and $\chi_{2n}^2$ differ in both mean and variance, whereas the proposed gamma distribution offers a better approximation. See figures \ref{fig:convg} and \ref{fig:quantilesconv} for additional parameter settings.

We evaluate the accuracy of type I error control for $S_n$ and $\tilde{S}_n$ using either $\chi^2_{2n}$ or the optimal gamma distributions by implementing the following rejection rules:
\begin{align}
	\label{eq:rejection_chisq_gamma}
	\begin{split}
		J_n &= I(S_n \geq \chi^2_{1-\alpha, 2n}), \quad
		J^*_n = I(\tilde{S}_n \geq \chi^2_{1-\alpha, 2n}),\\
		I_n &= I( S_n \geq q_{1-\alpha, 4n/\nu, \nu/2}), \quad
		I^*_n= I(\tilde{S}_n \geq q_{1-\alpha, nm^2/v, v/m}),
	\end{split}
\end{align}
where $\chi^2_{1-\alpha, 2n}$ denotes the $1-\alpha$ quantile of $\chi^2_{2n}$, and $I(A)$ represents the indicator function of event $A$. We compute the empirical type I error rates by averaging the indicators $\bar{I_n}$, $\bar{I^*_n}$, $\bar{J_n}$, and $\bar{J^*_n}$ over $20,000$ simulations. Accurate error control is demonstrated by how closely the empirical type I error rate matches the nominal value of $\alpha$. A larger (or smaller) empirical type I error rate than $\alpha$ suggests that the test is too liberal (or conservative).

\begin{figure}[h!]
	\begin{centering}
		{\includegraphics[width=1\linewidth]{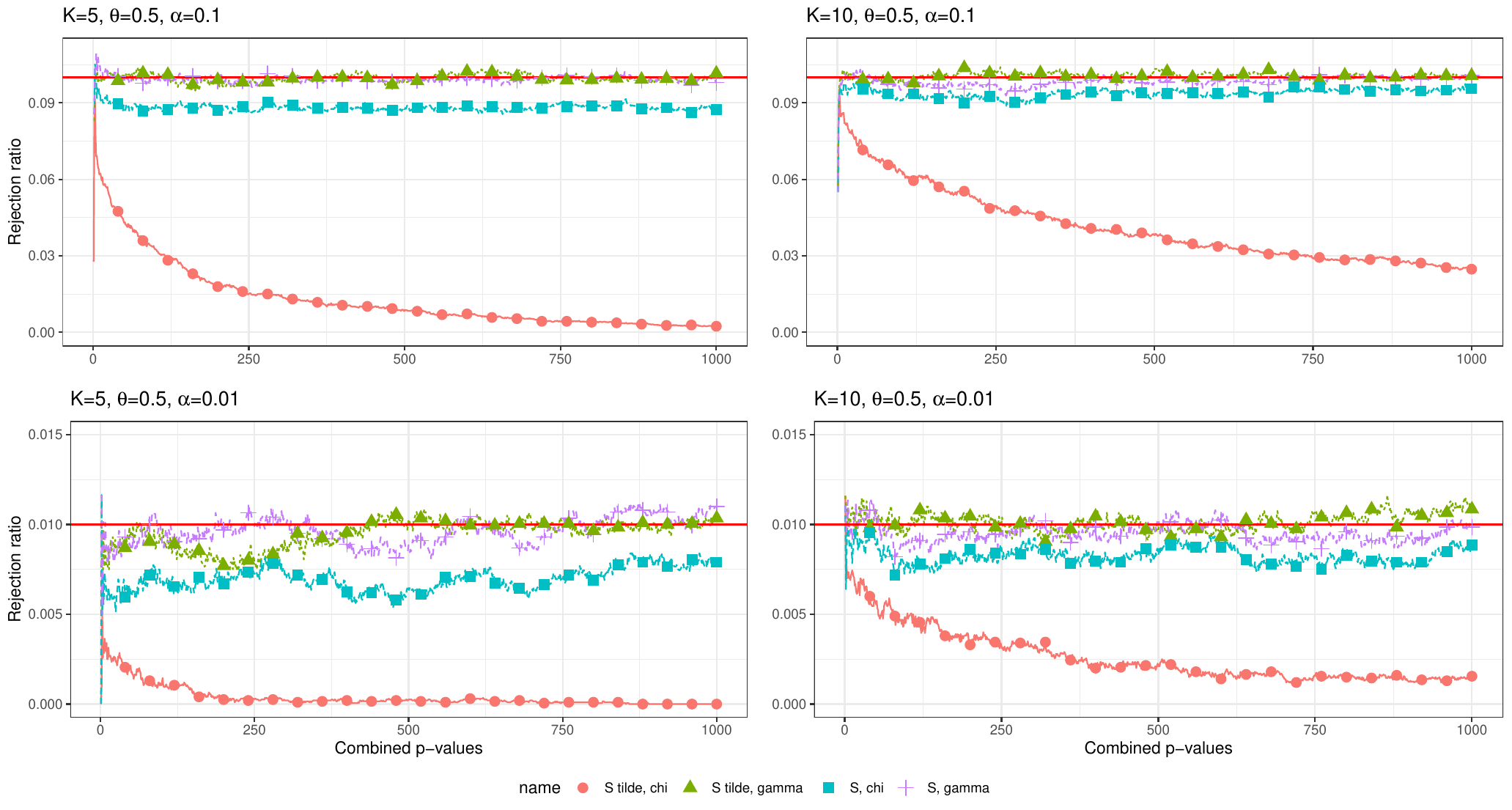}}
		\caption{Empirical type I error rates at $\theta_0=0.5$ and various $K$ and $\alpha$ values. 
			X-axis: Number of p-values combined. 
			The symbols for statistics and null distributions: dots: $\tilde{S}_n$ with $\chi^2_{2n}$; triangles: $\tilde{S}_n$ with optimal gamma; squares: $S_n$ with $\chi^2_{2n}$; and crosses: $S_n$ with optimal gamma.
		}
		\label{fig:pvcontrol05-part}
	\end{centering}
\end{figure}

Figure \ref{fig:pvcontrol05-part} compares the empirical type I error rates at $\theta_0=0.5$. It can be seen that using $\chi^2_{2n}$ to approximate the null distributions of $S_n$ and $\tilde{S}_n$ is too conservative: the type I error stabilizes at a value lower than the nominal $\alpha$ for $S_n$ and plummets to 0 for $\tilde{S}_n$. This observation is consistent with $S_n$ and $\tilde{S}_n$ having different moments from $\chi^2_{2n}$. On the other hand, using the optimal gamma distributions provides consistent type I error control, i.e., $\E(I_n)$ and $\E(I^*_n)$ approach $\alpha$ as $n$ increases, by Proposition \ref{thm:asymptoticpvalue}. Similar results are observed at other $K$ and $\theta_0$ values, as illustrated in Figure \ref{fig:pvcontrol05} in the Appendix. The conservativeness is more pronounced when $K$ is small or $\theta_0$ is close to 0.

\begin{figure}
	\begin{tabular}{p{0.42\linewidth} p{0.12\linewidth} p{0.42\linewidth}}
		\vspace{0pt} {\includegraphics[width=\linewidth]{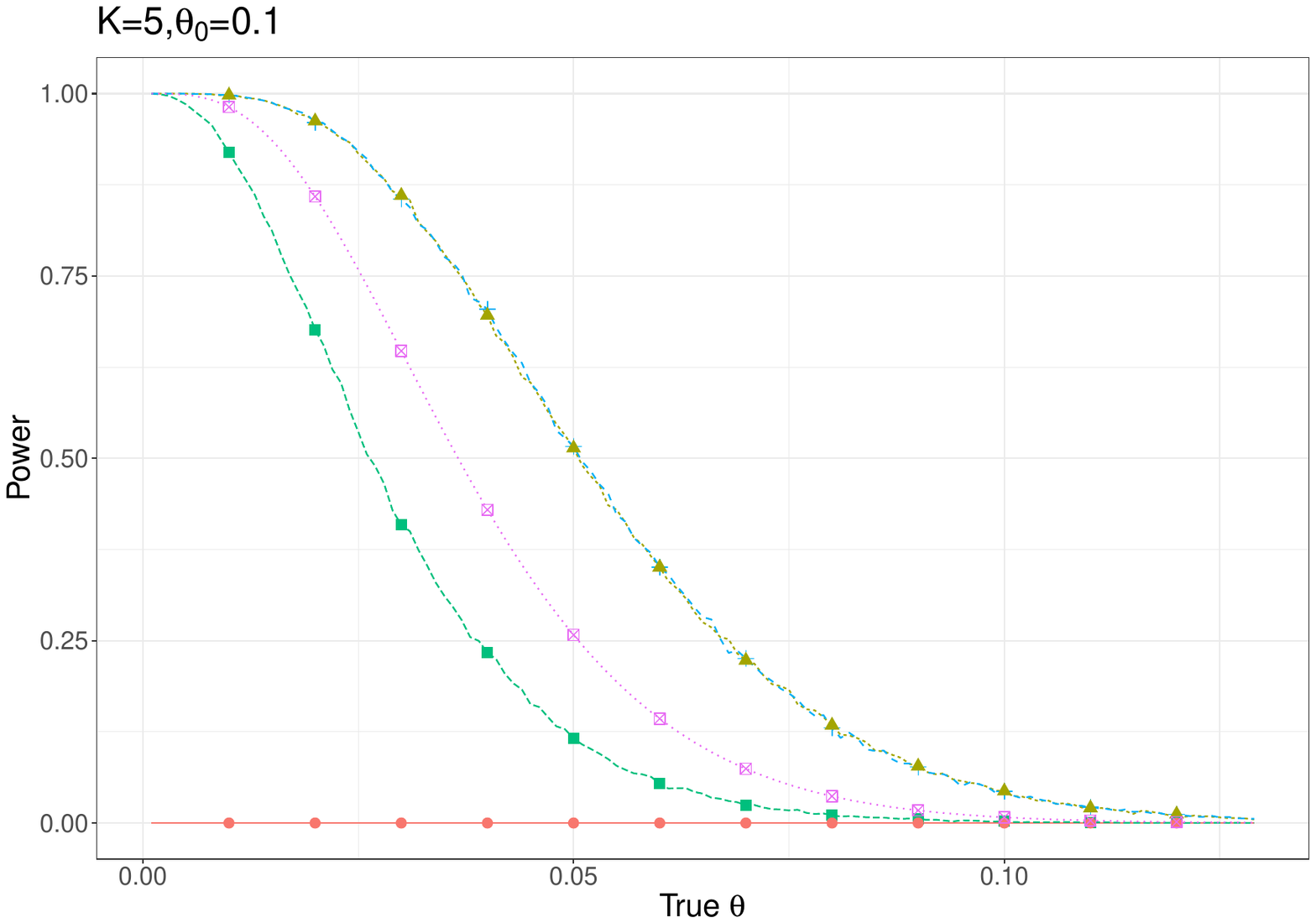}} &
		\vspace{20pt}
		{\includegraphics[width=\linewidth]{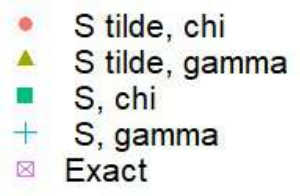}} &
		\vspace{0pt}
		{\includegraphics[width=\linewidth]{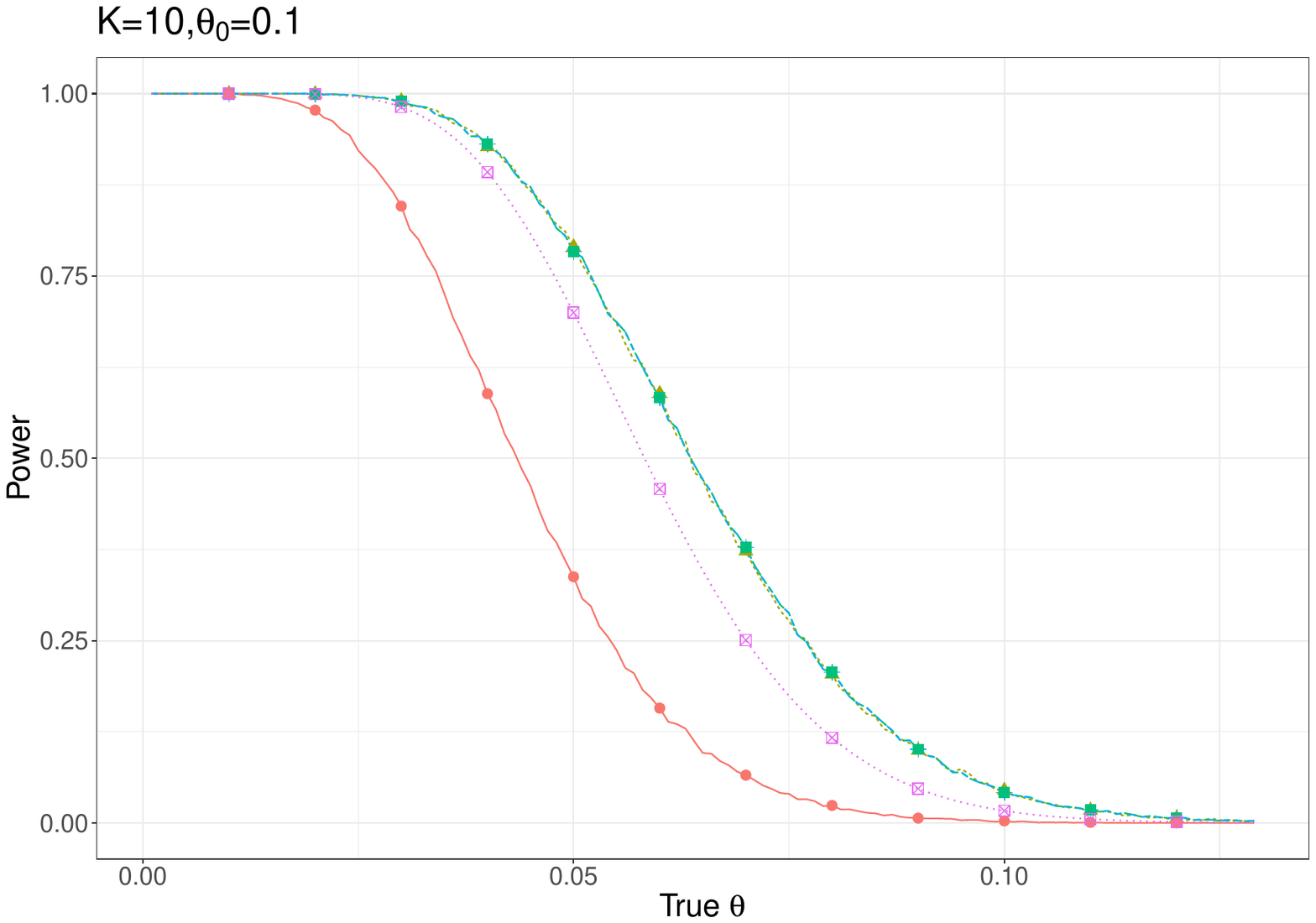}}   
	\end{tabular}
	\caption{Power curves over various actual $\theta$ values (x-axis) at size $\alpha=0.05$. 
		The null parameter is $\theta_0=0.1$ and $K=5$ (left panel) or $10$ (right). 
		Symbols: dots: $\tilde{S}_n$ with $\chi^2_{2n}$; triangles: $\tilde{S}_n$ with optimal gamma; solid squares: $S_n$ with $\chi^2_{2n}$; crosses: $S_n$ with optimal gamma; hollow squares: the exact test in \eqref{eq:pvalueExactBinom}. 
		The number of combined tests $n=40$; the number of simulations $N=10^4$. 
	}
	\label{fig:power}
\end{figure}

Under $H_A$ in \eqref{eq:H0A}, the empirical powers are the averages of $\bar{I_n}$, $\bar{I^*_n}$, $\bar{J_n}$, and $\bar{J^*_n}$ under the actual parameter $\theta$.
We also consider an exact binomial test, where the test statistic $T =  \sum_{j=1}^{n} X_j \sim \text{Binomial}(nK, \theta)$, since $X_j$'s are i.i.d. Binomial$(K, \theta)$.
The power function of this exact test is
\begin{equation}
	\label{eq:pvalueExactBinom}
	p(\theta)=\sum_{k=0}^{K_{\theta_0}} \binom{nK}{k}\theta^k(1-\theta)^{nK-k},
\end{equation}
where $K_{\theta_0}=\max \{k: \P(\text{Binomial}(nK,\theta_0)\leq k) \leq \alpha\}$. Figure \ref{fig:power} showcases a comparison of the power of the exact binomial test with those of $S_n$ and $\tilde{S}_n$ using the $\chi^2_{2n}$ or the optimal gamma distributions to control type I error at a level of $\alpha=0.05$ and with a null parameter of $\theta_0=0.1$. We simulated the data under the alternative hypothesis with $\theta\in (0.01, 0.12)$. Additional results for $\theta_0=0.5$ can be seen in Figure \ref{fig:power_theta0_0.5} of the Appendix.

The figure illustrates that the power of the exact test often surpasses the power of $S_n$ and $\tilde{S}_n$ when $\chi^2_{2n}$ is employed to control their type I error. This observation is consistent with the results documented by \cite{lancaster1949combination}. 
However, our proposed methods using the optimal gamma distributions have a higher power than the other three methods. The power curves for $S_n$ and $\tilde{S}_n$ using optimal gamma distributions are nearly identical, indicating the similar power.
Note that at $\theta=\theta_0$, the power is equivalent to the type I error rate. The power becomes less than the type I error rate when the actual value $\theta > \theta_0$, since the alternative hypothesis in \eqref{eq:H0A} is left-sided. 

\begin{figure}[h!]
	\begin{centering}
		{\includegraphics[width=0.45\linewidth]{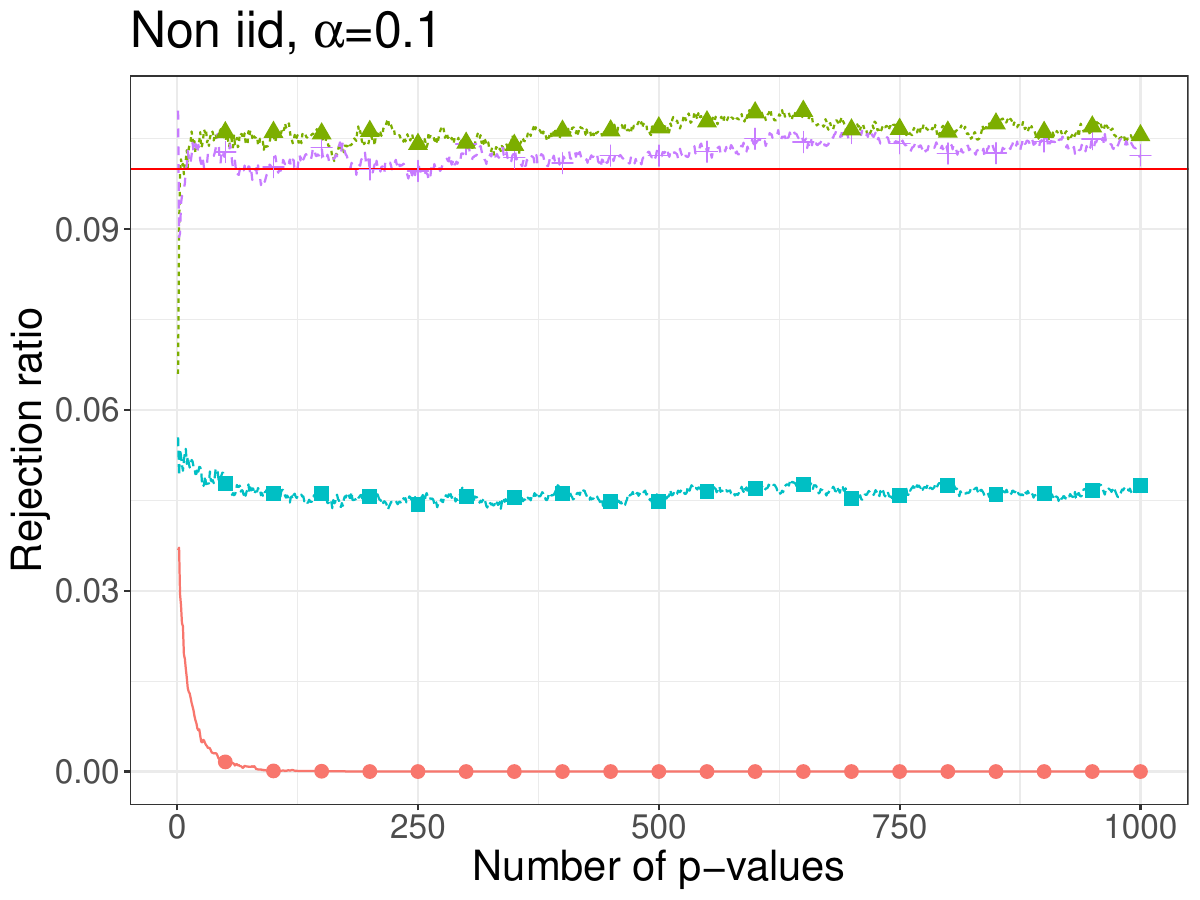}}
		{\includegraphics[width=0.45\linewidth]{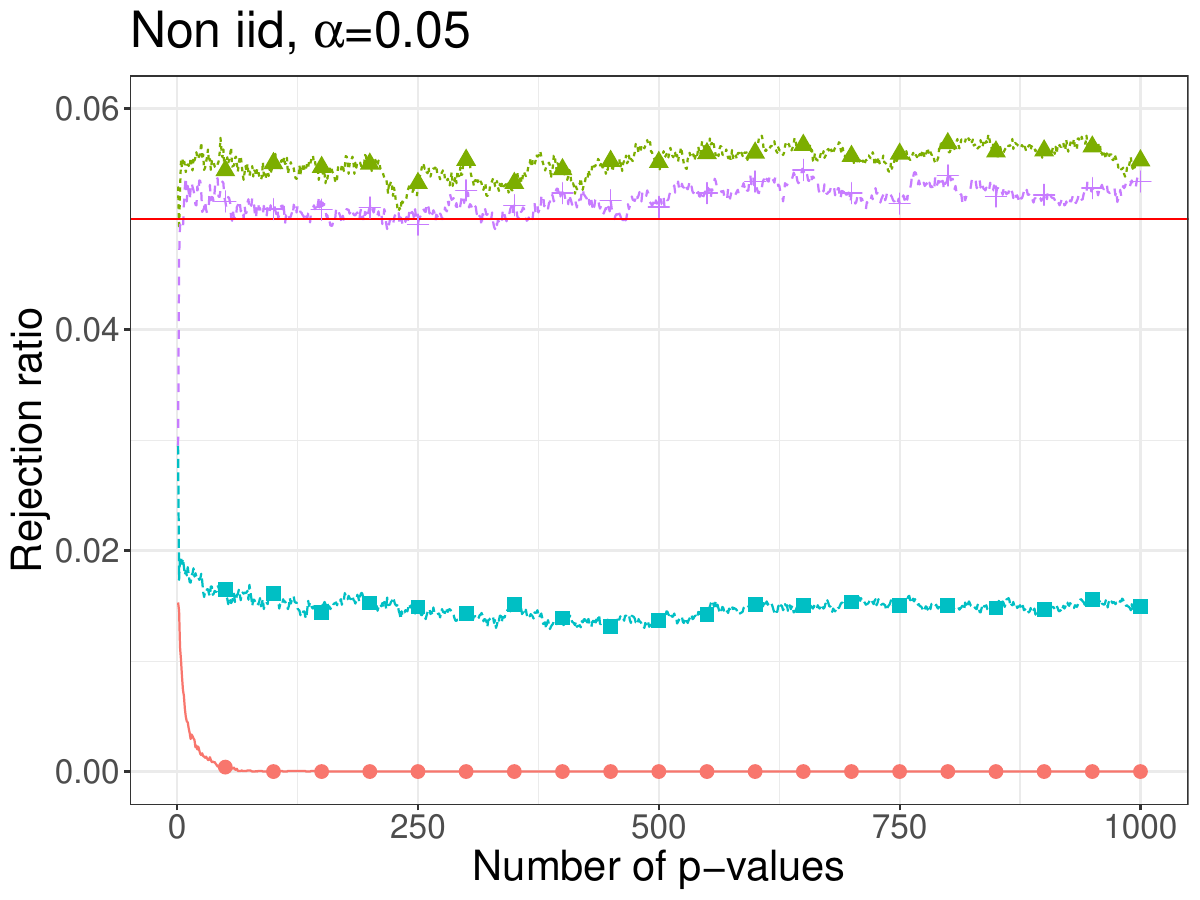}}\\
		{\includegraphics[width=0.4\linewidth]{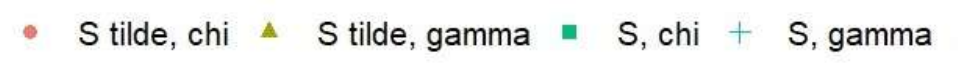}}
		\caption{Empirical type I error rates for combining independent but non-identically distributed p-values generated with equal frequency from one of the nine settings in Table \ref{table:paramfit}. Nominal level $\alpha=0.1$ (left panel) or $0.05$ (right).
			Two curves are close to $\alpha$ (indicated by the horizontal line): $S_n$ with optimal gamma (marked by crosses) and  $\tilde{S}_n$ with optimal gamma (triangles); two curves are significantly lower:   $S_n$ with $\chi^2_{2n}$ (squares) and $\tilde{S}_n$ with $\chi^2_{2n}$ (dots).
		}
		\label{fig:pvcontrolnoniid}
	\end{centering}
\end{figure}

\subsection{Combining independent, identically distributed hypergeometric p-values}

Our framework is general and readily applicable to arbitrary distributions. To demonstrate this, we consider p-values derived from i.i.d. discrete test statistics following the noncentral hypergeometric distribution \citep{wallenius1963biased}. This distribution corresponds to Fisher's exact test and is broadly applied in practical studies. Specifically, the discrete input statistics are $X_j \sim \text{NCHYG}(M,m,K,\omega)$, for $j=1,\ldots,n$, and the probability mass function of $X_j$ is proportional to 
\begin{equation}\label{eq:nchyg}
	P(X_j=x) \propto \binom{m}{x} \binom{M}{K-x} \omega^x 
	\quad \text{ for } \max\{0,K-M\} \leq x \leq \min\{K,m\}.
\end{equation}

This distribution can be viewed as the conditional distribution of $Y_1$ given $Y_1+Y_2=K$, where $Y_1 \sim \text{Binomial}(m,\pi_1)$, $Y_2 \sim \text{Binomial}(M,\pi_2)$ and $\omega=\frac{\pi_1(1-\pi_2)}{\pi_2(1-\pi_1)}$, initially considered in the context of a contingency table \citep{fisher1922contingency}. It is often used in case-control studies. For example, in the study of genetic associations \citep{lin2014association,lin2016beyond},  $m$ and $M$ are the numbers of cases and controls, respectively, $K$ is the number of mutations, and $x$ is the number of mutations in the case group. Under the null hypothesis of no association, the $K$ mutations are randomly distributed between the case and control groups. In this scenario, $\omega=1$ and the distribution reduces to the usual hypergeometric distribution. However, under the alternative hypothesis that an association exists, the mutations will have different probabilities of being present in the two groups, leading to $\omega \neq 1$.

In the study, we test the global hypotheses $H_0: \pi_1=\pi_2$ versus $H_A: \pi_1>\pi_2$ (i.e., cases have a higher mutation rate),  which is equivalent to
\begin{equation}
	\label{eq:H0cc}
	H_0:\omega=1 \text{ vs. } H_A:\omega >1.
\end{equation}
We choose the parameters $m=M=4000$ for a balanced study, and a small $K=5$ to emulate rare genetic variants and to achieve a far-from-continuous p-value distribution. Table \ref{table:paramfitcc} shows the distribution of $X_j$ and the p-value $P_j$ in \eqref{eq:disc_P_left} under the null hypothesis in \eqref{eq:H0cc}.

\begin{table*}[h!]
	\caption{Values of $X$ and corresponding values of $\ddot{Z}_j$ and $\tilde{Z}_j$}
	\label{table:paramfitcc}	
	\begin{tabular}{l|llllll}
		\hline
		Value of $X$           & 0      & 1      & 2      & 3      & 4      & 5      \\ \hline
		$P(X_j=x)$             & 0.0312 & 0.1562 & 0.3126 & 0.3126 & 0.1562 & 0.0312 \\
		$P_j=P(X_j\leq x)$         & 0.0312 & 0.1874 & 0.5    & 0.8126 & 0.9688 & 1      \\
		$\ddot{Z}_j$ when $X_j=x$  & 8.9339 & 4.6325 & 2.2096 & 0.8615 & 0.2341 & 0.0315 \\
		$\tilde{Z}_j$ when $X_j=x$ & 8.3203 & 4.427  & 2.136  & 0.8423 & 0.2315 & 0.0314\\
		\hline
	\end{tabular}
\end{table*}

\begin{figure}[h!]
	\begin{centering}
		{\includegraphics[width=0.9\linewidth]{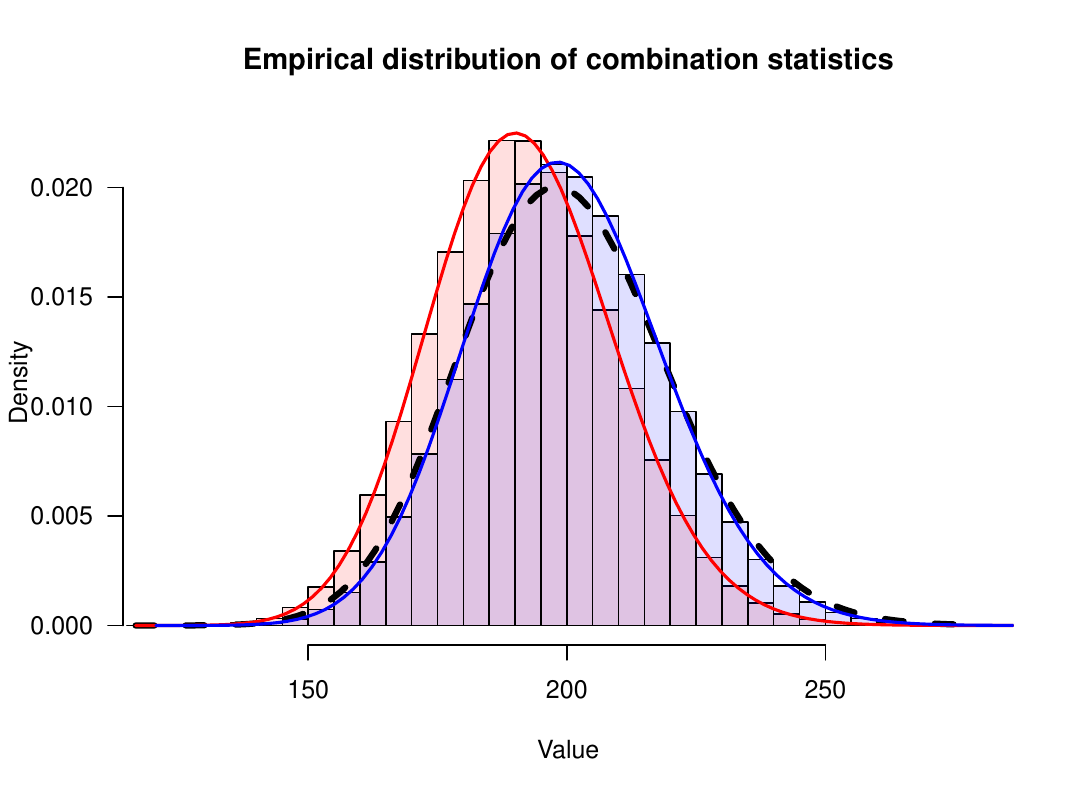}}
		\caption{
			Density of $\chi_{2n}^2$ (dashed curve) and optimal gamma distributions (solid curves), contrasted with histograms of simulated statistics $S_n$ (blue, left) and $\tilde{S}_n$ (red, right), combining $n=100$ i.i.d. hypergeometric p-values.	}
		\label{fig:histcc}	
	\end{centering}
\end{figure}

Under the null hypothesis, we have $E(\ddot{Z}_j)=2$, $\nu=\text{Var}(\ddot{Z}_j)=3.61$, $m=E(\tilde{Z}_j)=1.916$ and $v=\text{Var}(\tilde{Z}_j)=3.1765$. For $n=100$, we obtain $10^7$ simulated values of $S_n$ and $\tilde{S}_n$ to assess the goodness of fit of the optimal gamma and chi-squared distribution to their empirical distributions. Figure \ref{fig:histcc} displays the results. The figure illustrates that the variance of $S_n$ is similar to that of $\chi_{2n}^2$. However, the distribution of $S_n$ is better approximated by the gamma distribution at the tails. $\tilde{S}_n$ and $\chi_{2n}^2$ differ in both mean and variance, whereas the proposed gamma distribution offers a better approximation across all quantiles.

Table \ref{table:alphacc} shows the proportion of rejected tests using quantiles of the chi-squared distribution and the optimal gamma distributions at various nominal $\alpha$ levels. For both statistics, the empirical rejection rates are close to $\alpha$ under their gamma approximations. In contrast, using the quantiles of the chi-squared distribution leads to conservative decisions, especially when $\alpha$ is small.

\begin{table}[h!]
		\label{table:alphacc}	
		\begin{tabular}{lllll}
			\hline
			\multicolumn{1}{|l|}{Rejection for $S_n$}         & \multicolumn{1}{l|}{$\alpha=0.05$} & \multicolumn{1}{l|}{$\alpha=0.01$} & \multicolumn{1}{l|}{$\alpha=0.005$} & \multicolumn{1}{l|}{$\alpha=0.001$} \\ \hline
			\multicolumn{1}{|l|}{$\chi^2_{2n}$}               & 0.0406                             & 0.006                              & 0.0026                              & \multicolumn{1}{l|}{0.0003}                              \\ \cline{1-1}
			\multicolumn{1}{|l|}{Gamma}                       & 0,0487                             & 0.0086                             & 0.0044                              & \multicolumn{1}{l|}{0.0008}                              \\ \hline
			&                                    &                                    &                                     &                                     \\ \hline
			\multicolumn{1}{|l|}{Rejection for $\tilde{S}_n$} & \multicolumn{1}{l|}{$\alpha=0.05$} & \multicolumn{1}{l|}{$\alpha=0.01$} & \multicolumn{1}{l|}{$\alpha=0.005$} & \multicolumn{1}{l|}{$\alpha=0.001$} \\ \hline
			\multicolumn{1}{|l|}{$\chi^2_{2n}$}               & 0.0113                             & 0.0015                             & 0.0004                              & \multicolumn{1}{l|}{0}                                   \\ \cline{1-1}
			\multicolumn{1}{|l|}{Gamma}                       & 0.0484                             & 0.0086                             & 0.0044                              & \multicolumn{1}{l|}{0.0008}                              \\ \hline
		\end{tabular}
		\caption{Proportion of rejected combinations $S_n$ and $\tilde{S}_n$ at different $\alpha$ levels based on  chi-squared and the gamma distribution quantiles.}
	\end{table}
	
\begin{figure}[h!]
		\begin{tabular}{p{0.48\linewidth} p{0.48\linewidth}}
			\vspace{0pt} {\includegraphics[width=\linewidth]{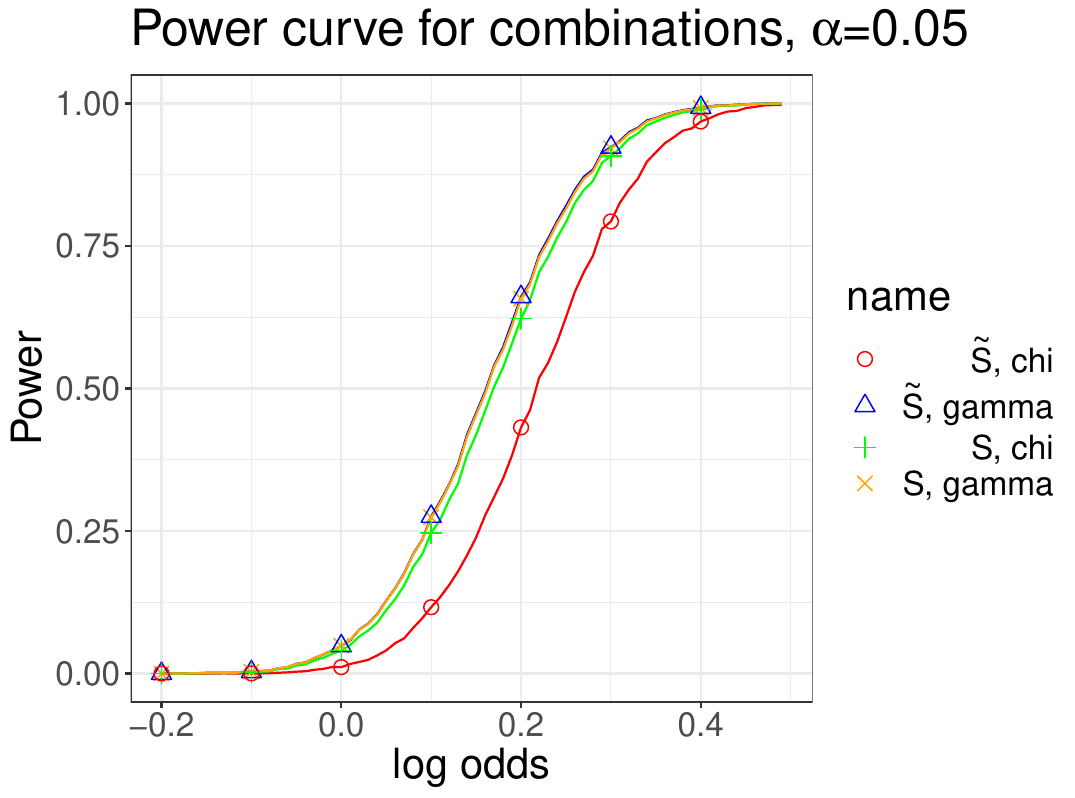}} &
			\vspace{0pt}
			{\includegraphics[width=\linewidth]{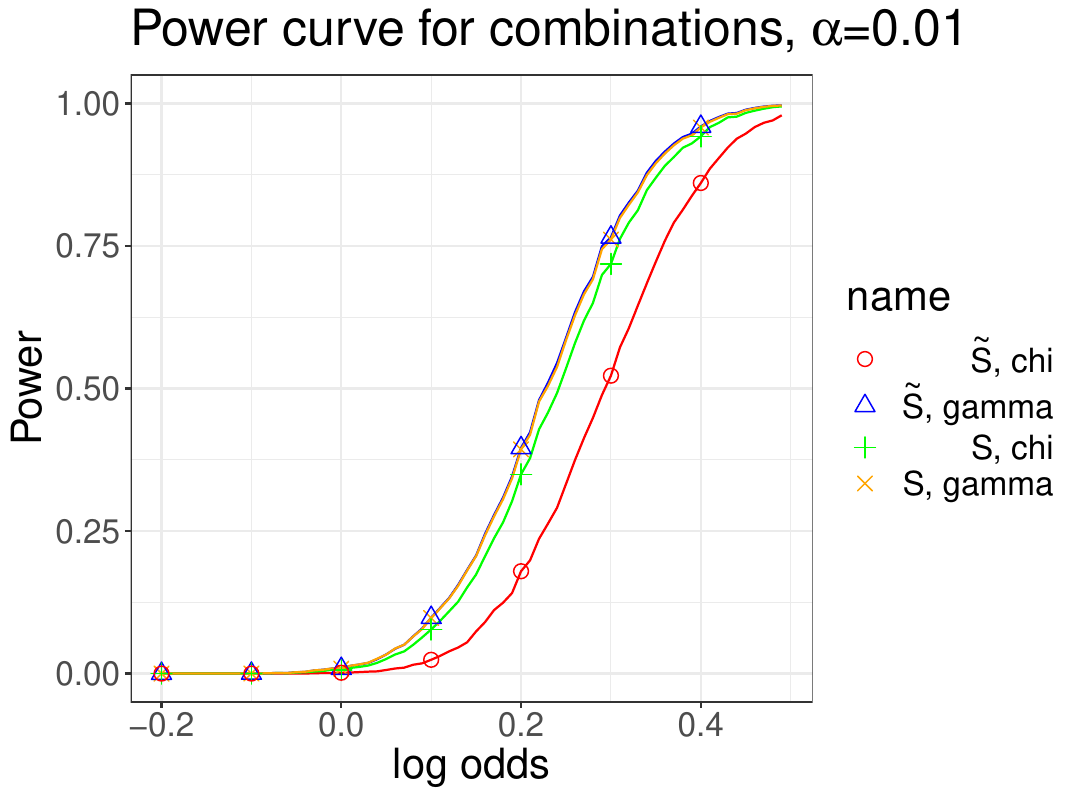}}   \\
			\vspace{0pt} {\includegraphics[width=\linewidth]{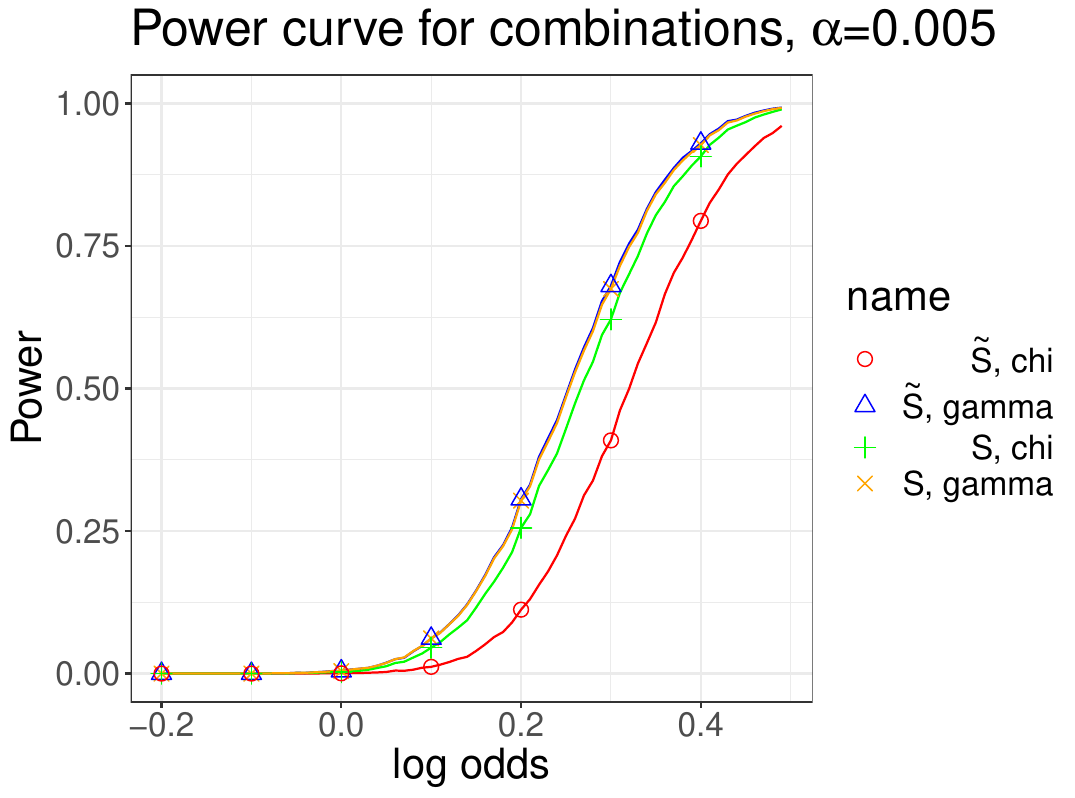}} &
			\vspace{0pt}
			{\includegraphics[width=\linewidth]{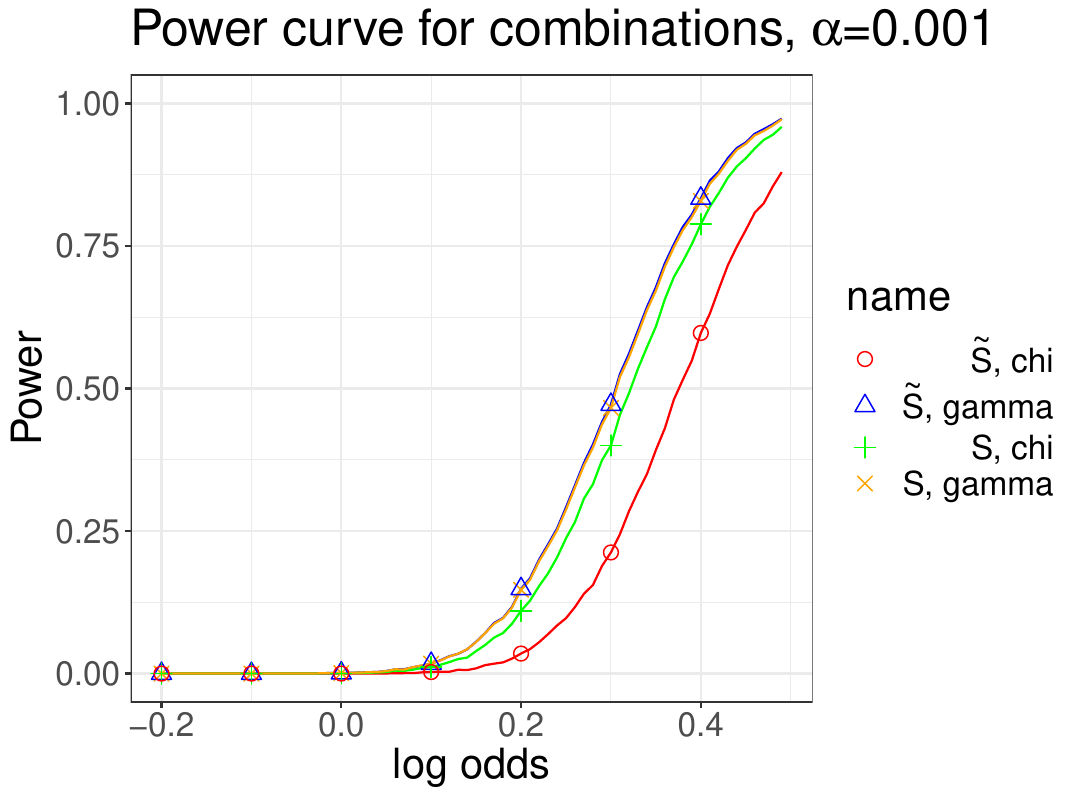}}   
		\end{tabular}
		\caption{Power curves over $\log(\omega)$ at various $\alpha$ levels.
			The number of combined tests is $n=100$; the number of simulations is $N=10^4$. 
		}
		\label{fig:powercc}
	\end{figure}

	To assess power, we simulate $N=10^4$ sets of $n=100$ statistics under \eqref{eq:nchyg} with $\log (\omega)\in (-1,1)$. Figure \ref{fig:powercc} illustrates that our proposed methods using the optimal gamma distributions have higher power than the chi-squared distribution.
	Power for $S_n$ and $\tilde{S}_n$ using the optimal gamma distribution is similar. The results are consistent across various $\alpha$ levels.

\subsection{Combining independent, non-identically distributed p-values}

In this section, we extend the numerical study to scenarios where $n$ independent but non-identically distributed p-values are combined. Specifically, we consider input discrete test statistics $X_j \sim \text{Binomial}(K_j,\theta_{0j})$, $j=1,\ldots,n$, under the null hypothesis. Each pair $(K_j, \theta_{0j})$ originates from one of the nine settings in Table \ref{table:paramfit} with equal frequency. To approximate the null distributions of $S_n$ and $\tilde{S}_n$, we employ both $\chi^2_{2n}$ and the optimal gamma distributions, as outlined in \eqref{eq:approxnoniid}.

The goodness of fit of the optimal gamma distributions for the empirical null distributions of $S_n$ and $\tilde{S}_n$ is seen in Appendix Figures \ref{fig:convg} and \ref{fig:quantilesconv}..
Furthermore, we evaluate the accuracy of type I error controls by the rejection given in \eqref{eq:rejection_chisq_gamma}, except using $\bar{\nu}$, $\bar{m}$, and $\bar{v}$ (the averages of the parameters across all distributions from Table \ref{table:paramfit}) to calculate the gamma quantiles.
Figure \ref{fig:pvcontrolnoniid} illustrates that the empirical error rates align closely with the nominal size $\alpha$ when we employ the optimal gamma null distributions. In contrast, when $\chi_{2n}^2$ are used as the null distribution, the error rates for both $S_n$ and $\tilde{S}_n$ can be seen overly conservative.

\begin{figure}
	\begin{centering}
		{\includegraphics[width=0.55\linewidth]{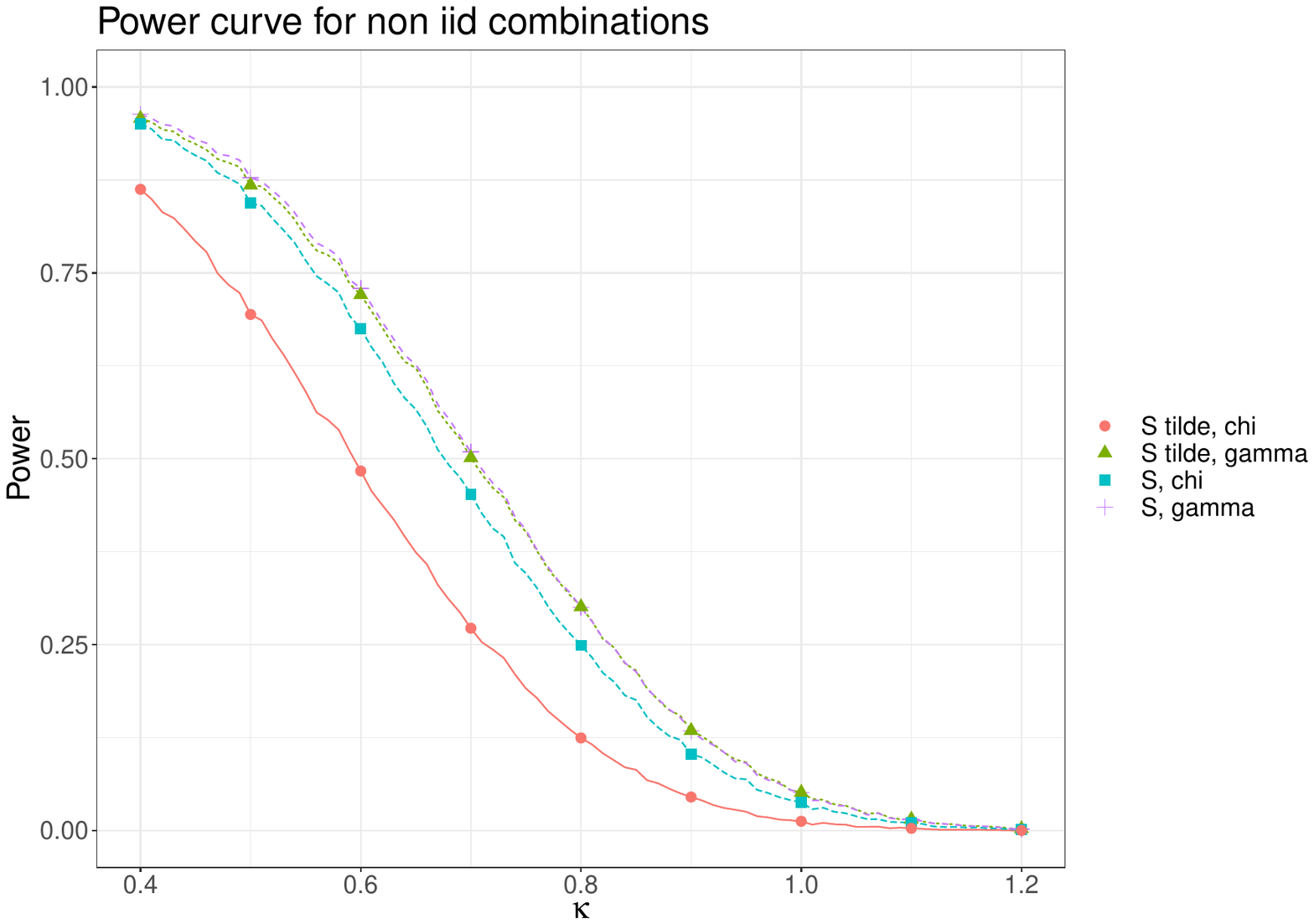}}
		\caption{Power curves over $\kappa\in (0.4, 1.2)$. Alternative parameters $\theta_j = \kappa \theta_{0j}$ based on null parameters $\theta_{0j}$ in Table \ref{table:paramfit}. Test level $\alpha=0.05$. 
			Symbols: dots: $\tilde{S}_n$ with $\chi^2_{2n}$; 
			Triangles: $\tilde{S}_n$ with optimal gamma;
			Squares: $S_n$ with $\chi^2_{2n}$; 
			Crosses: $S_n$ with optimal gamma. 
			$n=40$. $10^4$ simulations. 
		} 
		\label{fig:powerniid}
	\end{centering}
\end{figure}

Under the same null hypothesis and rejection rules, we assess the statistical power using an alternative hypothesis that $X_j \sim \text{Binomial}(K_j, \kappa \theta_{0j})$ over $\kappa \in (0.4, 1.2)$. The p-values are left-sided for consistency with previous studies.
We calculate the empirical power as the proportion of rejections across $10^4$ simulations. The power curves are illustrated in Figure \ref{fig:powerniid}. The figure suggests that, for both $S_n$ and $\tilde{S}_n$,  tests using the gamma distributions offer higher statistical power than those using the $\chi^2_{2n}$ distributions. We remark that the exact test based on statistic $T=\sum_{j=1}^nX_j$ examined earlier is not easily applicable when $X_j$'s are non-identical because the exact distribution of $T$ is no longer binomial. 

\section{Discussion}\label{sec:discu}

This paper presents a new framework for the combination of discrete statistics using a continuous gamma approximation. It addresses the issues pointed out by \cite{Berry}  concerning the unsuitability of the chi-squared distribution in Fisher's combination of discrete p-values. This novel approach provides a more accurate representation of the distributions and enhances the power of statistical tests while maintaining control over the type I error rates.

The framework expands on the fundamental concept of Lancaster's procedures, which use continuous distributions to approximate discrete statistics for testing purposes. We have demonstrated that suitable continuous distributions can effectively approximate properly adjusted discrete statistics. Nevertheless, there will always be a discrepancy between them due to the inherent discreteness. In Appendix \ref{sec:sfig}, we analyzed this discrepancy measured by the Wasserstein distance. We provided several lower bounds to account for the discrepancy due to the largest probability mass of the discrete statistic. The analysis help quantify the performance and limitation of the approximation process. It also illustrates how the optimal gamma approximation lessens the discrepancy compared to the chi-squared approximation. In Appendix \ref{sec:sfig}, we also discussed the limitations of directly using the original Fisher's statistic even if a gamma distribution could be used to approximate its distribution. These limitations further justify the use of the adjusted statistics for combining discrete p-values.

Drawing on this paper, an essential topic for future research involves addressing the challenge of combining dependent discrete p-values. This problem presents a steeper hurdle than the case with independent variables, given that the interdependence among p-values can significantly shape the behavior of the combined statistic. One potential solution to this problem is the development of methods capable of gauging the distance between random vectors that house dependent variables. These methods could facilitate the quantification of the discrepancy between the true joint distribution of the p-values and their approximation. 

\begin{acks}[Acknowledgments]
We thank the partial support of Contador and Wu by the U.S. National Science Foundation grants DMS-2113570 and DMS-1812082. 
\end{acks}

\begin{appendix}
	\section{Proofs}\label{sec:proofs}

	\subsection{Proof of Theorem \ref{thm:wpdc}}
	
	\begin{proof}
		To see that $\pi$ is a coupling, note that
		\begin{align*}
			\pi(A\times \mathbb{R})&=\sum_{i\in \mathbb{N}I} \delta_A(x_i)\int_{ [G^{-1}(F_{i-1}), G^{-1}(F_i)]} g(y)dy\\
			&=\sum_{i\in \mathbb{N}} \delta_A(x_i) \left(G(G^{-1}(F(x_i)))-G(G^{-1}(F(x_{i-1})))\right)\\
			&=\sum_{i\in \mathbb{N}} \delta_A(x_i) \left(F(x_i)-F(x_{i-1})\right)\\
			&=\sum_{i\in \mathbb{N}} \delta_A(x_i)\P(X=x_i)=\P(X\in A),
		\end{align*}
		and
		\begin{align*}
			\pi( \mathbb{R} \times B)&=\sum_{i\in \mathbb{N}} \int_{ B\cap [G^{-1}(F_{i-1}), G^{-1}(F_i)]} g(y)dy\\
			&=\sum_{i\in \mathbb{N}} \int_{ [G^{-1}(F(x_{i-1}), G^{-1}(F(x_i)]}\delta_B(y)g(y)dy\\
			&=\int_\mathbb{R} \delta_B(y)g(y)dy\\
			&=\int_Bg(y)dy =\P(Y \in B).
		\end{align*}
		So the marginals of $\pi$ are the distributions of $X$ and $Y$, respectively. Note also that for $w \in (F(x_{i-1}),F(x_i)]$, $F^{-1}(w)=\inf \{x:F(x)\geq w\}=x_i$.  We have
		\begin{align*}
			\int_{\mathbb{R}^2} |x-y|^p \pi(dx,dy)&=\sum_{i\in \mathbb{N}} \int_{ [G^{-1}(F_{i-1}), G^{-1}(F_i)]} |x_i-y|^pg(y)dy\\ 
			&= \sum_{i\in \mathbb{N}} \int_{F(x_{i-1})}^{F(x_i)} |x_i-G^{-1}(w)|^pdw \\
			&= \sum_{i\in \mathbb{N}} \int_{F(x_{i-1})}^{F(x_i)} |F^{-1}(w)-G^{-1}(w)|^pdw\\
			&=\int_{0}^{1} |F^{-1}(w)-G^{-1}(w)|^pdw,
		\end{align*}
		which is the one dimensional formulation of Wasserstein distance in \eqref{eq:wass}. It is therefore proven that $\int_{\mathbb{R}^2} |x-y|^pd\pi (x,y)=\inf_{\gamma \in C(P,Q)}\int_{\mathbb{R}^2} |x-y|^pd\gamma (x,y)$, thus $\pi$ achieves the minimum expected distance among all copulas, as desired.
	\end{proof}

	\subsection{Proof of Corollary \ref{cor:wdmoments}}
	
	\begin{proof}
		With measure $\pi$ defined in Definition \ref{def:coupling}, the distance between $X$ and $Y$ can be written as
		\begin{align*}
			W^2_2(X,Y)=&E_\pi((X-Y)^2)\\
			=& E_\pi(X^2+Y^2-2XY)\\
			=&\iint (x^2+y^2-2xy)\pi(x,y)dxdy.
		\end{align*}
		As $\pi$ is a coupling, it holds that
		$$\int x^2 \left(\int\pi(x,y)dy\right)dx=\sum_{i \in \mathbb{N}}x_i^2\P(X=x_i)=\Var(X)+\E(X)^2,$$
		$$\int y^2 \left(\int\pi(x,y)dx\right)dy=\int y^2dG(y)=\Var(Y)+\E(Y)^2,$$
		$$E_\pi(XY)=\iint xy\pi(x,y)dxdy=\sum_{i \in \mathbb{N}} x_i\int_{G^{-1}(F(x_{i-1}))}^{G^{-1}(F(x_{i}))} yg(y) dy.$$
		
		So
		$$W^2_2(X,Y)=\Var(X)+\Var(Y)+(\E(X)-\E(Y))^2-2\Cov_\pi(X,Y).$$
		
		The result follows from $\Cov_\pi (X,Y)=\E(X(Y-\mu_Y))$ after observing that
		\begin{align*}
			\sum_{i \in \mathbb{N}} \mu_X\int_{G^{-1}(F_{i-1})}^{G^{-1}(F_i)} (y-\mu_Y)g(y) dy&=\mu_X\sum_{i \in \mathbb{N}}\int_{G^{-1}(F(x_{i-1}))}^{G^{-1}(F(x_i))} (y-\mu_Y)g(y) dy\\
			&=\mu_X\int_{-\infty}^\infty yg(y)dy-\mu_X\mu_Y\int_{-\infty}^\infty g(y)dy\\
			&=\mu_X\mu_Y-\mu_X\mu_Y=0.
		\end{align*}
		
		Where $\mu_X$ and $\mu_Y$ denote the first moments of $X$ and $Y$, respectively.
	\end{proof}

	\subsection{Proof of Lemma  \ref{lem:newm}}
	
	\begin{proof}
		We prove the result by contradiction. Clearly, $\mathcal{M}_X^* \subset \mathcal{M}_X$. Let $V \in \mathcal{M}_X$ and suppose that $V \notin \mathcal{M}_X^*$. 
		This means that, for some $i$ it must hold that $V=v_i$ when $X=x_i$ but $v_i \notin A_i=[G^{-1}(F_{i-1}),G^{-1}(F_i)]$. It's clear that $G^{-1}$ is increasing as $G$ is a CDF and is increasing. We consider two cases:
		
		If $v_i<G^{-1}(F_{i-1})$, then $|v_i-a|>|G^{-1}(F_i)-a|$ for all $a \in A_i$ and

		$$\int_{F_{i-1}}^{F_{i}}|v_i-G^{-1}(w)|^pdw > \int_{F_{i-1}}^{F_{i}}|x_i-G^{-1}(w)|^pdw.$$
		
		Similarly, if $v_i>G^{-1}(F_i)$, then

		$$\int_{F_{i-1}}^{F_{i}}|v_i-G^{-1}(w)|^pdw > \int_{F_{i-1}}^{F_{i}}|G^{-1}(F_i)-G^{-1}(w)|^pdw.$$
		
		Let $I=\{i: v_i \notin A_i \}$. Define 
		
		$$V_0=\begin{cases}
			v_i & i\notin I,\\
			G^{-1}(F_{i-1}) & z_i \leq G^{-1}(F_{i-1}) ,\\
			G^{-1}(F_i) & z_i \geq G^{-1}(F_i).
		\end{cases}$$

		It follows that
		\begin{itemize}
			\item $V_0$ has the same probability masses than $Z$.
			\item $V_0 \in \mathcal{M}_X^*.$ 
			\item $W^p_p(V_0,Y)<W^p_p(V,Y).$
		\end{itemize}
		So, for any variable $V \in \mathcal{M}_X$, if $V \notin \mathcal{M}_X^*$ one can build another variable $V_0 \in \mathcal{M}_X^*$ that is closer to $Y$ than $V$ is. Therefore it holds that
		$$\inf_{V \in \mathcal{M}_X} W_p^p(V,Y)=\inf_{V \in \mathcal{M}_X^*}W_p^p(V,Y).$$
	\end{proof}

	\subsection{Proof of Theorem \ref{thm:zopti}}
	
	\begin{proof}
		
		Denote by $F$, $F_X$, and $G$ the cumulative distribution functions of $Z$, $X$, and $Y$, respectively. By definition, $F(z_i)=F_X(x_i)$ for any index $i$ and therefore $F^{-1}(w)=z_i$ when $w \in (F_{i-1},F_{i})$  according to the result of Lemma \ref{lem:newm}. Thus, by \eqref{eq:wass} we have
		
		$$W_p^p(Z, Y)=\sum_{i \in \mathbb{N}} \int_{F_{i-1}}^{F_{i}}|z_i-G^{-1}(w)|^pdw.$$
		
		So, in order to minimize $W_p^p(Z,Y)$, each of the (positive) integrals in the above sum should be minimized. Observe that
		
		$$\int_{F_{i-1}}^{F_{i}}|z_i-G^{-1}(w)|^pdw=\int_{A_i}|z_i-y|^pg(y)dy,$$
		where $g(y)=G'(y)$ is the density of the target distribution $Y$. For each $i$, the integral operator 
		
		$$T(w)=\int_{A_i}|w-y|^pg(y)dy=\E(|w-Y|^p|A_i)\P(Y\in A_i)$$
		is continuously differentiable on $w$, with derivative given by 
		\begin{align*}
			T'(w)=&\int_{G^{-1}(F_{i-1})}^{w}p(w-y)^{p-1}g(y)dy - \int_w^{G^{-1}(F_{i})}p(y-w)^{p-1}g(y)dy\\
			=&p\left[ \E((w-Y)^{p-1}\bm{1}_{w>Y}|A_i)-\E((Y-w)^{p-1}\bm{1}_{w<Y}|A_i) \right]\P(Y\in A_i).
		\end{align*}
		
		By convexity, $T(w)$ will be minimized when
		
		$$\E((w-Y)^{p-1}\bm{1}_{w>Y}|A_i)=\E((Y-w)^{p-1}\bm{1}_{w<Y}|A_i),$$
		with the equality achieved when $w^{p-1}=\E\left(Y^{p-1}|A_i\right)$. The result follows.
	\end{proof}

	\subsection{Proof of Corollary \ref{cor:zopti2}}
	\begin{proof}
		Following Theorem \ref{thm:zopti}, when $p=2$, a straightforward deduction gives 
		
		$$z_i=\E(Y|A_i)=\frac{\int_{G^{-1}(F_{i-1})}^{G^{-1}(F_{i})}yg(y)dy}{\int_{G^{-1}(F_{i-1})}^{G^{-1}(F_{i})}g(y)dy}=\frac{\int_{G^{-1}(F_{i-1})}^{G^{-1}(F_{i})}yg(y)dy}{\P(Z=z_i)}.$$
		
		If $X=P$ is the discrete p-value and $Y\sim \text{Uniform}(0,1)$ with $g(y)=1$, we have 
		$$z_i=\frac{[G(G^{-1}(F_{i}))]^2-[G(G^{-1}(F_{i-1}))]^2}{2[G(G^{-1}(F_{i}))-G(G^{-1}(F_{i-1}))]}=\frac{F_i+F_{i-1}}{2},$$
		which is the $i$-th value in the support of the mid-p statistic in (\ref{eq:midp_P}). 
		
		Similarly, if $X=-2\log(P)$ and $Y \sim \chi^2_2$ with $g(y)=e^{-y/2}/2$ for $y>0$, we have
		$$z_i=\frac{2(F_i-F_{i-1}) - 2(F_i\log F_i-F_{i-1} \log F_{i-1}) }{(F_i-F_{i-1})},$$
		which is Lancaster's mean-value-$\chi^2$ in (\ref{eq:lmod_lexplicit}). 
	\end{proof}

	\subsection{Proof of Corollary \ref{cor:unbiased} }
	\begin{proof}
		With the same notation as in Theorem \ref{thm:zopti}, for the Wasserstein metric of order $p=2$,
		\begin{align*}
			\E(Z)&=\sum_{i \in \mathbb{N}}z_i\P(Z=z_i)\\
			&=\sum_{i \in \mathbb{N}} E\left(Y|Y \in A_i \right)\P(Z=z_i)\\
			&=\sum_{i \in \mathbb{N}} \int_{G^{-1}(F_{i-1})}^{G^{-1}(F_{i})}yg(y)dy\\ 
			&= \int_\mathbb{R}yg(y)dy  =\E(Y).
		\end{align*}
		And therefore $\Var(Y)-\Var(Z)=\E(Y^2)-\E(Z^2).$
		Write
		$$\E(Y^2)=\int y^2 g(y)dy =\sum_{i \in \mathbb{N}} \int_{A_i} y^2 g(y)dy, \quad \E(Z^2)=\sum_{i \in \mathbb{N}} z_i^2 p_i =\sum_{i \in \mathbb{N}} z_i^2 (F_i-F_{i-1}).$$
		
		The variance difference becomes
		\begin{align*}
			\E(Y^2)-\E(Z^2) &=\sum_{i \in \mathbb{N}} \left[\int_{A_i} y^2 g(y)dy - z_i^2 (F_i-F_{i-1}) \right] \\
			&= \sum_{i \in \mathbb{N}} (F_i-F_{i-1})  \left[\frac{\int_{A_i} y^2 g(y)dy}{\int_{A_i} g(y)dy} - z_i^2 \right]\\
			&= \sum_{i \in \mathbb{N}} (F_i-F_{i-1})\left(E(Y^2\big|A_i) -z_i^2\right).
		\end{align*}
		By Jensen's conditional inequality, we have that
		$$\E(Y^2|A_i)\geq \E(Y|A_i)^2=z_i^2,$$
		meaning that the sum in $\E(Y^2)-\E(Z^2)$ consists of all non-negative terms, and as $Y$ is non constant in $A_i$ for some $i$, $\E(Y^2|A_i) > \E(Y|A_i)^2=z_i^2$ and strict inequality follows for the variance when $p=2$.

		Considering the Wasserstein metric with $p > 2$, we have that the function $\phi (x) = x^{p-1}$ is strictly convex and increasing. So, using Jensen's inequality in $A_i$ we have
		$$\phi(z_i)=\E(\phi(Y)|A_i)\geq \phi(\E(Y|A_i)) \implies z_i \geq \phi^{-1} (\phi(\E(Y|A_i)))=\E(Y|A_i).  $$

		Since $Y$ is continuous, whenever $\P(A_i)>0$ it holds that $\P(Y\neq Z|Z)=\P(Y\neq z_i|A_i)=1$, so the above inequalities are strict. These strict inequalities imply that 
		\begin{align*}
			\E(Z)&=\sum_{i \in \mathbb{N}}z_i \P(Z=z_i)\\
			&>\sum_{i \in \mathbb{N}}\E(Y|A_i) \P(A_i)\\
			&=\sum_{i \in \mathbb{N}}\int_{A_i} yg(y)dy=\E(Y).
		\end{align*}  
		
		The inequality statement for $p \in (1,2)$ follows similarly: $\phi (x) = x^{p-1}$ is strictly concave and thus $-\phi$ is convex for $p<2$. Therefore, $\E(Z) <  \E(Y)$ for $p \in (1,2)$.
		
		In terms of the variance, for any $p>1$ and $\phi (x) = x^{p-1}$, the law of total variance guarantees that, with $\mathcal{A}$ being the sigma-algebra given by the partition $(A_i)_{i \in \mathbb{N}}$,
		$$\Var(Y^{p-1})=\Var(\phi(Y)) = \Var (\E(\phi(Y)|\mathcal{A}))+\E(\Var (\phi(Y)|\mathcal{A})) \geq  \Var (\E(\phi(Y)|\mathcal{A})) =\Var(Z^{p-1}). $$

		To see that this inequality is strict, observe that in each $A_i$, $$\Var (\phi(Y)|A_i)=\int_{A_i} \left(\phi(y)-\frac{\int_{A_i} \phi(z) g(z)dz}{\P(A_i)}\right)^2 g(y)dy>0$$ where the integrand is positive since $\phi$ is non-constant over $A_i$, thus $\E(\Var (\phi(Y)|\mathcal{A}))>0$ and the desired inequality follows.
	\end{proof}

	\subsection{Proof of Corollary \ref{cor:wdmoments2}}
	\begin{proof}
		Following Corollary \ref{cor:unbiased}, we have $\E(Y)=\E(Z)$. Using the formula for $W_2^2 (Z,Y)$ given by Corollary \ref{cor:wdmoments} and the formula for $z_i$ given by Theorem \ref{thm:zopti}, we have
		\begin{align*}
			W_2^2 (Z,Y)&= \Var(Z)+\Var(Y)-2\left(\sum_{i \in \mathbb{N}} z_i\int_{G^{-1}(F_{i-1})}^{G^{-1}(F_i)} yg(y) dy-\E(Y)\E(Z)\right)\\
			&= \Var(Z)+\Var(Y)-2\left(\sum_{i \in \mathbb{N}} z_i[z_i\P(Z=z_i)]-\E(Z)^2\right)\\
			&= \Var(Z)+\Var(Y)-2\left(\E(Z^2)-\E(Z)^2\right)\\
			&= \Var(Y)-\Var(Z).
		\end{align*}

		According to the law of total variance \citep{Bill86}, $\Var(Y)= \E (\Var(Y|X))+\Var(\E(Y|Z))=\E (\Var(Y|X))+\Var(Z)$ and thus $W_2^2(Z,Y)=\E (\Var(Y|X))$.
	\end{proof}

	\subsection{Proof of Corollary \ref{cor:closetc}}
	\begin{proof}
		We observe that, as in the proof of Corollary \ref{cor:unbiased} and using the form of the points in the support of the minimizer $z_i$ given in Corollary  \ref{cor:zopti2}, for each $i$,
		\begin{align*}
			\int_{A_i} y^2 g(y)dy - z_i^2 (F_i-F_{i-1}) &=\int_{A_i} (y^2-z_i^2) g(y)dy \\
			&\leq 2G^{-1}(F_i) (F_i-F_{i-1})^2.
		\end{align*}
		
		Adding the above terms across $i \in \mathbb{N}$ yields
		\begin{align*}
			\E(Y^2)-\E(Z^2) &=\sum_{i \in \mathbb{N}} \left[\int_{A_i} y^2 g(y)dy - z_i^2 (F_i-F_{i-1}) \right] \\
			&\leq \sum_{i \in \mathbb{N}}  2G^{-1}(F_i) (F_i-F_{i-1})^2\\
			&\leq  2\Delta \sum_{i \in \mathbb{N}}  G^{-1}(F_i) (F_i-F_{i-1}).
		\end{align*}
		
		When $\Delta \rightarrow 0$, the Riemman sum on the right above converges to 
		$$\sum_{i \in \mathbb{N}}  G^{-1}(F_i) (F_i-F_{i-1})\rightarrow \int_0^1G^{-1}(w)dw=\E(Y)<\infty,$$
		
		and, as $\Var(Y)\geq \Var(Z)$ and $	\E(Y)=\E(Z)$, using the result of Corollary 4 it follows that as $\Delta \rightarrow 0$,
		$$0\leq W_2^2 (Z,Y)=\Var(Y)-\Var(Z)= \E(Y^2)-\E(Z^2)\leq 2\Delta \sum_{i \in \mathbb{N}}  G^{-1}(F_i) (F_i-F_{i-1}) \rightarrow 0,$$
		which shows that $W_2^2 (Z,Y)\to 0$.
	\end{proof}

	\subsection{Proof of Corollary \ref{cor:closetochi1}}
	\begin{proof}
		Denote by $z_i$ and $\tilde{z}_i$ the values of $\ddot{Z}$ and $\tilde{Z}$ when $P=F_i$. We have, 
		\begin{align*}
			W_2^2(\tilde{Z},\chi_2^2) 
			& =8+v+m^2 +4\sum_{i \in \mathbb{N}}\log\left(\frac{F_i+F_{i-1}}{2}\right)\int_{-2\log F_i}^{-2\log F_{i-1}}\frac{y}{2}e^{-y/2}dy \nonumber\\
			&=8+v+m^2 +8\sum_{i \in \mathbb{N}}\log\left(\frac{F_i+F_{i-1}}{2}\right) (F_i-F_{i-1} -F_i\log F_i +F_{i-1}\log F_{i-1}) \nonumber\\
			&=4+v+(m-2)^2 -8\sum_{i \in \mathbb{N}}\log\left(\frac{F_i+F_{i-1}}{2}\right) (F_i\log F_i -F_{i-1}\log F_{i-1} )\\
			&=4+v+(m-2)^2 -2\sum_{i \in \mathbb{N}}\left(-2\log\frac{F_i+F_{i-1}}{2}\right) \left(-2\frac{F_i\log F_i -F_{i-1}\log F_{i-1}}{F_i -F_{i-1}} \right)(F_i -F_{i-1})\\
			&=4+v+(m-2)^2 -2\sum_{i \in \mathbb{N}}\tilde{z}_i (z_i-\E(\ddot{Z}))\\
			&= 4+v+(m-2)^2 -2\E[\tilde{Z}(\ddot{Z}-\E(\ddot{Z}))].
		\end{align*}
		Furthermore, since $\E[\tilde{Z}(\ddot{Z}-\E(\ddot{Z}))]=\Cov(\tilde{Z},\ddot{Z})$, we have
		$$v+\nu -2\E[\tilde{Z}(\ddot{Z}-\E(\ddot{Z}))]=\Var(\ddot{Z})+\Var(\tilde{Z})-2\Cov(\tilde{Z},\ddot{Z}) = \Var(\ddot{Z} -\tilde{Z}) \geq 0.$$
		As $(m-2)^2>0 $ and $W_2^2(Z,\chi_2^2)= 4-\nu$ according to Corollary \eqref{cor:wdmoments2}, we get the desired difference
		$$W_2^2(\tilde{Z},\chi_2^2)-W_2^2(\ddot{Z},\chi_2^2)=\Var(\ddot{Z} -\tilde{Z}) +(m-2)^2 >0.$$
	\end{proof}

	\subsection{Proof of Theorem \ref{thm:optgamma}}

	\begin{proof}
		The logic of this proof is as follows: We prove that, for a continuous approximation to discrete $X$ that is arbitrarily close to it, minimizing Wasserstein distance to a gamma distribution is achieved by matching the first two moments according to the Euler-Lagrange formulation for integral functional minimization. We conclude using convergence in distribution granted by convergence of Wasserstein metric.
		
		We start by constructing, for arbitrary $\epsilon>0$, a continuous random variable $X_\epsilon$ that is within $\epsilon$ Wasserstein distance of discrete $X$, i.e.,  $W_2(X_\epsilon, X) \leq \epsilon$. To do so, let $\kappa:\mathbb{R}\rightarrow \mathbb{R}_+$ be an infinitely differentiable, non-negative function with $\kappa\equiv 0$ outside of $[0,1]$ and $\int_0^{1}\kappa(x)dx=1$.  As $X\geq0$, using $x_0=0$ and $\{p_i\}_{i \in \mathbb{N}}$ its probability mass function given by \eqref{eq:disc_P_left} in the main document, we define $h_\epsilon(x):\mathbb{R}\rightarrow \mathbb{R}_+$  and $H_\epsilon(x):\mathbb{R}\rightarrow \mathbb{R}_+$, respectively, as 
		$$
		h_\epsilon(x)=\frac{1}{\epsilon}\sum_{i \in \mathbb{N}} p_i\kappa\left(  \frac{x-x_i}{\epsilon} \right); \quad 
		H_\epsilon(x)=\int_0^x h_\epsilon(z)dz.
		$$
		$H_\epsilon$ is increasing since $h_\epsilon$ is positive. Also, $H_\epsilon(0)=0$ and
		$$
		\lim_{x \to \infty} H_\epsilon(x)
		=\sum_{i \in \mathbb{N}} p_i\int_{x_i}^{x_i+\epsilon}\frac{1}{\epsilon}\kappa\left(\frac{x-x_i}{\epsilon}\right)dx
		=\sum_{i \in \mathbb{N}} p_i\left(\int_0^1\kappa(w)dw\right)
		=1.
		$$
		Thus, $h_\epsilon$ is a probability density function and $H_\epsilon$ is its corresponding cumulative distribution function, and we can define a positive continuous random variable $X_\epsilon$ with this density. If $x\in (x_{i}+\epsilon, x_{i+1})$ for any $i$ (that is, if $x$ is not immediately within $\epsilon$ units to the right of a discontinuous point of $F_X$), $h_\epsilon (x)= 0$, and therefore, by Theorem \ref{thm:wpdc},
		$$W_2^2(X_\epsilon, X)=\sum_{i \in \mathbb{N}}\int_{x_{i}}^{x_{i}+\epsilon}(x_{i}-y)^2h_\epsilon(y)dy\leq \epsilon^2\sum_{i \in \mathbb{N}}p_i=\epsilon^2.$$
		Therefore, we have that $W_2(X_\epsilon, X) \leq \epsilon$ as desired. 
		
		Now we show that the continuous $X_\epsilon$ and the gamma random variable closest to it share the same first two moments. For fixed $\epsilon>0$ and the corresponding cumulative distribution function $H_\epsilon$ of $X_\epsilon$ defined in the first part of the proof, we define an operator $L_\epsilon$ from $x$ in $\mathbb{R}$, a twice differentiable continuous function $f$ and its derivative $f'$ to the real numbers by 
		$$L_\epsilon(x,f,f'):=(H_\epsilon^{-1} \circ f(x)-x)^2f'(x).$$  
		
		Observe that, for any $\epsilon>0$, $L_\epsilon$ is differentiable. Using $f=G$ and $f'=g$, i.e., the cumulative distribution function and the density of $\tilde{Y}\sim Gamma(\alpha,\beta)$, we have
		\begin{align*}
			W_2^2(X_\epsilon,\tilde{Y})&=\int_0^1|H_\epsilon^{-1}(w)-G^{-1}(w)|^2 dw\\
			&=\int_{0}^{\infty}(H_{\epsilon}^{-1}[G(y)]-y)^2g(y)dy\\
			&=\int_{0}^{\infty} L_\epsilon(y,G,g)dy.
		\end{align*}
		
		In order to obtain parameters $\alpha$ and $\beta$ for the distribution of $\tilde{Y}$ closest to $X_\epsilon$, we use the Euler-Lagrange one-dimensional formulation on the integrand of $W_2^2(X_\epsilon,\tilde{Y})=\int_{0}^{\infty} L_\epsilon(y,G,g)dy$. According to the Euler-Lagrange formulation, the functional form $g$ minimizing this integral satisfies the equality $\frac{\partial L_\epsilon}{\partial f}=\frac{d}{dx}\frac{\partial L_\epsilon}{\partial f' }$. This equation yields
		
		$$(H_\epsilon^{-1})'[G(x)](H_\epsilon^{-1} \circ G(x)-x)(1-g(x))= (H_\epsilon^{-1} \circ G(x)-x).$$
		Using standard measure theory arguments, this equality is equivalent to the fact that for every Lebesgue integrable function $\theta$, the following equality holds:
		\begin{equation}
			\label{eq:ELcondition}
			\int(H_\epsilon^{-1})'[G(y)](H_\epsilon^{-1}[G(y)]-y)(1-g(y)) \theta(y)dy=\int(H_\epsilon^{-1}[G(y)]-y)\theta(y)dy,
		\end{equation}
		By definition, for any given $i$, when $y\in (G^{-1}(F_X(x_{i}))+\epsilon, G^{-1}(F_X)(x_{i+1}))$,  $H_{\epsilon}^{-1}[G(y)]=x_i$ is constant and its derivative $(H_\epsilon^{-1})'[G(y)]=0$ on such intervals. Thus, upon taking $\theta_i$ supported on $(G^{-1}(F_X(x_{i}))+\epsilon, G^{-1}(F_X)(x_{i+1}))$ (i.e., $\theta_i$ is zero everywhere else), and observing that the left hand side integrand is null, we get 
		\begin{equation}
			\label{eq:variationalg}
			x_i \int_{G^{-1}(F_X(x_i))+\epsilon}^{G^{-1}(F_X(x_{i+1}))}  \theta_i(y)dy=\int_{G^{-1}(F_X(x_i))+\epsilon}^{G^{-1}(F_X(x_{i+1}))}  y\theta_i(y)dy.
		\end{equation}
		
		Define $J_{i,\epsilon}=(G^{-1}(F_X(x_{i-1}))+\epsilon, G^{-1}(F_X(x_{i})))$ and let $A_\epsilon$ be the event where $X_\epsilon$ is further than $\epsilon$ to the left of a discontinuity of $F_X$, this is,
		$$A_\epsilon:=\cup_{i \in \mathbb{N}}\{X_\epsilon \in J_{i,\epsilon} \}.$$
		Observe that $A_\epsilon$ is nonempty for small enough $\epsilon$ and that $X$ and $X_\epsilon$ have the same probability distribution 
		in $A_\epsilon$. Using $\theta_i(y)=g(y)I(J_{i,\epsilon})(y)$ in  \eqref{eq:variationalg} one obtains 
		
		$$\E(X I(A_\epsilon))= \E(X_\epsilon I(A_\epsilon))= \sum_{i \in \mathbb{N}} x_i \int_{G^{-1}(F_X(x_i))+\epsilon}^{G^{-1}(F_X(x_{i+1}))}  \theta_i(y)dy =\int_{A_\epsilon} yg(y)dy . $$
		
		Similarly, choosing now $\theta_i(y)=yg(y)I(J_{i,\epsilon})(y)$ in \eqref{eq:variationalg}, one obtains
		$$\E(X^2_\epsilon I(A_\epsilon))=\E([X_\epsilon I(A_\epsilon)]^2)=\int_{A_\epsilon} y^2g(y)dy, $$
		which implies 
		$$\Var(X I(A_\epsilon))=\int_{A_\epsilon} y^2g(y)dy-\left[\int_{A_\epsilon} yg(y)dy \right]^2.$$
		
		Thus, for fixed $\epsilon$,  $W^2_2(X_\epsilon,\tilde{Y})$ achieves its single stationary point by choosing $\tilde{Y}=Y_\epsilon$ with shape and rate parameters matching the first two moments of $X_\epsilon$. This point is a global minimum because $W^2_2(X_\epsilon,\tilde{Y})$ grows to infinity as either $\alpha$ or $\beta$ grows to infinity, and it increases to $W^2_2(X_\epsilon,\delta_0)=\E(X_\epsilon^2)>W^2_2(X_\epsilon,\tilde{Y})$ when either  $\alpha$ or $\beta$ decreases to zero. 
		
		Since $A_\epsilon$ decreases to a countable set as $\epsilon \rightarrow 0$, by the monotone convergence theorem, the moments of $X I(A_\epsilon)$ grow to the moments of $X$ and the parameters of the density $g$ of $\tilde{Y}$ are given by 
		$$\lim_{\epsilon \to 0}\E(X I(A_\epsilon))= \E(X)=\int_0^\infty yg(y)dy=\alpha \beta,$$
		and 
		$$\lim_{\epsilon \to 0}\Var(X I(A_\epsilon)) = \Var (X)=\int_0^\infty [y-\alpha\beta]^2g(y)dy=\alpha\beta^2.$$
		
		Thus, the unique stationary point of $\int_0^\infty L(y,g,G)dy=\lim_{\epsilon \to 0}\int_0^\infty L_\epsilon(y,g,G)dy$ is achieved for $g$ and $G$ with shape and rate parameters $\alpha=\frac{[\E(X)]^2}{\Var(X)}$ and $\beta =\frac{\Var (X)}{\E(X)}$. To see that this point is indeed a minimizer of $W_2(X,Y)$, observe that, by the triangle inequality, one has for any $\epsilon>0$ and  $Y\sim Gamma(\alpha,\beta)$,
		
		\begin{align*}
			W_2(X,Y)&\leq   W_2(X,X_\epsilon)+W_2(X_\epsilon,Y_\epsilon)+W_2(Y,Y_\epsilon)\\
			&\leq \epsilon + \inf_{\hat{Y}\sim Gamma}W_2(X_\epsilon,\hat{Y})+W_2(Y,Y_\epsilon).
		\end{align*}
	
		As $\epsilon \to 0$, $Y_\epsilon \to Y$ in distribution which implies $W_2(Y,Y_\epsilon)\to 0$ (since the Wasserstein distance metrizes weak convergence). Moreover, as $X_\epsilon \to X$ in distribution (i.e., $H_\epsilon \to F_X$ pointwisely), an application of Fatou's Lemma \citep[Theorem 16.3]{Bill86} gives $\lim_{\epsilon \to 0}\inf_{\hat{Y}\sim Gamma}W_2(X_\epsilon,\hat{Y})\leq \inf_{\hat{Y}\sim Gamma}W_2(X_,\hat{Y})$. Therefore, upon taking limits as $\epsilon \to 0$, the above inequality yields:
		
		$$	W_2(X,Y)\leq  \inf_{\hat{Y}\sim Gamma}W_2(X,\hat{Y}).$$
		
		This is, taking $Y$ with the same first two moments as $X$ achieves the minimum distance among all gamma distributions.
	\end{proof}

	\subsection{Proof of Corollary \ref{cor:convex}}
	To prove Corollary \ref{cor:convex}, we will use Auxiliary Lemma \ref{thm:equalpart} given below. It is proven afterwards. 
	\begin{lemma}
		\label{thm:equalpart}
		Let $T$ be a strictly monotonic function, $X$ be a discrete random variable, and $Y$ be a continuous random variable. Denote by $F_i^T$ the probability $\P(T(X)\leq T(F_i))$and by $A_i^T$ the partition they induce on the domain of $T(Y)$ as in \eqref{eq:partition}. Then, 
		$$\{F^T_i\}_{i \in \mathbb{N}}=\begin{cases}
			\{F_i\}_{i \in \mathbb{N}} & \mbox{ if  } T \mbox{ is increasing},\\
			\{1-F_{i-1}\}_{i \in \mathbb{N}} & \mbox{ if  } T \mbox{ is decreasing},\\
		\end{cases}$$
		and
		$$\{A_i^T\}_{i \in \mathbb{N}}=\{T(A_i)\}_{i \in \mathbb{N}}.$$
	\end{lemma}
	
	\textbf{Proof of Corollary \ref{cor:convex}}
	\begin{proof}
		If $T(x_i)$ is the $i$-th value in the support of $T(X)$ and $z_i^T$ is the $i$-th value in the support of $Z^T$, using the same notation as in the proof of Lemma \ref{thm:equalpart}, one has that for $p=2$,
		\begin{align*}
			z_i^T&= \E(T(Y)|T(Y) \in A_i^T)\\
			&=\E(T(Y)|T(Y) \in T(A_i))\\
			&=\E(T(Y)|Y \in A_i)\\
			&\geq T\left[\E(Y|Y \in A_i)\right]\\
			&= T(x_i),\\
		\end{align*}
		where we used the invertibility of $T$, the result of Lemma \ref{thm:equalpart}, and Jensen's conditional inequality. Since this is valid for any value in the support of $T(Z)$ and the corresponding value in the support of $Z^T$, the result follows.
	\end{proof}

	\bigskip
	\noindent {\bf Proof of Auxilliary Lemma \ref{thm:equalpart}}
	
	\begin{proof}
		Denote by $G$ the cumulative distribution function of $Y$. Strict monotonicity implies that $T^{-1}$ exists. We distinguish between the two cases for the monotonicity of $T$:
		\begin{itemize}
			\item If $T$ is increasing:\\
			One has for the $i$-th value in the domain of $X$.
			$$F_i=P(X\leq x_i)=P(T(X)\leq T(x_i))= F^T_i.$$
			
			The cumulative distribution function of $T(Y)$ is given by $G \circ T^{-1}$ and therefore, 
			\begin{align*}
				A_i^T &= [(G \circ T^{-1})^{-1}(F_{i-1}^T),(G \circ T^{-1})^{-1}(F_i^T)]\\
				&=[T \circ G^{-1}(F_{i-1}), T \circ G^{-1}(F_{i})]\\
				&=T(A_i).\\
			\end{align*}
			
			\item If $T$ is decreasing:\\
			We rewrite $\{x_i\}_{i \in \mathbb{N}}$ in decreasing order and with that we have that the domain of $T(X)$, $\{T(x_i)\}_{i \in \mathbb{N}}$, is in increasing order. The $i$-th quantile of $T(X)$ is then given by
			\begin{equation}
				\label{eq:cdfT}
				F^T_i=P(T(X)\leq T(x_i))=P(X\geq x_i)=1-P(X< x_{i})=1-F_{i-1}
			\end{equation}

			The cumulative distribution function of $T(Y)$, $\tilde{G}$, is now given by $$\tilde{G}(x)=1-G \circ T^{-1}(x),$$ and therefore, its inverse is given by
			$$\tilde{G}^{-1}(y)=T \circ G^{-1}(1-y).$$ 
			
			With that, we have that the $i$-th element of partition of the domain of $T(Y)$ is given by
			\begin{align*}
				A_i^T &= [\tilde{G}^{-1}(F_{i-1}^T), \tilde{G}^{-1}(F_{i}^T)]\\
				&=[T \circ G^{-1}(1-F^T_i), T \circ G^{-1}(1-F^T_{i-1})]\\
				&=T(A_{i}),\\
			\end{align*}
			where we used the previous result and the fact that the quantiles of $X$ are also in decreasing order.
		\end{itemize}
		
	\end{proof}

	\subsection{Proof of Proposition \ref{thm:asymptoticpvalue}}
	\begin{proof}
		Let $\{X_i\}_{i \in \mathbb{N}}$ be an i.i.d. sequence of positive random variables and define $\alpha=[\E(X_1)]^2/\Var (X_1)$, $ \beta=\Var(X_1)/\E(X_1)$ and $T_n=\sum_{i=1}^nX_i$. Let $\epsilon >0$ be an arbitrary small number. By the central limit theorem, one can find $n_0$ such that for $n\geq n_0$, 
		$$\left|\P(T_n > q_{p; n\alpha, \beta})-\int_{\frac{(q_{p; n\alpha, \beta})-n\alpha\beta}{\beta\sqrt{n\alpha} }}^\infty \frac{1}{\sqrt{2\pi}}e^{-t^2/2} dt\right|\leq \frac{\epsilon}{2}.$$
		
		The definition of $q_{p; n\alpha, \beta}$ implies that 
		\begin{align*}
			1-p&=\int_{q_{p; n\alpha, \beta}}^\infty \frac{1}{\Gamma(n\alpha)\beta^{n\alpha}}x^{n\alpha-1}e^{-x/\beta}dx\\
			&=\int_{q_{p; n\alpha, \beta}/\beta}^\infty \frac{1}{\Gamma(n\alpha)}y^{n\alpha-1}e^{-y}dy\\
			&=\int_{\frac{(q_{p; n\alpha, \beta})-n\alpha\beta}{\beta\sqrt{n\alpha} }}^\infty 
			\frac{\sqrt{n\alpha}}{\Gamma(n\alpha)}(\sqrt{n\alpha}t+n\alpha)^{n\alpha-1}e^{-\sqrt{n\alpha}t-n\alpha}dt.
		\end{align*}
		
		Using Stirling's approximation and a second order approximation on $h(t)=e^{-\sqrt{n}t}$, one gets that there exists $n_1$ such that for $n\geq n_1$, 
		$$e^{-\sqrt{n\alpha}t}\left|\frac{\sqrt{n\alpha}(\sqrt{n\alpha}t+n\alpha)^{n\alpha-1}}{\Gamma(n\alpha)}e^{-n\alpha}-\frac{1}{\sqrt{2\pi}}\left( 1+\frac{t}{\sqrt{n\alpha}}\right) \right|\leq e^{-t^2/2}\frac{\epsilon}{4}(1+O(n^{-1/2})).$$
		
		Therefore,
		\begin{align*}
			\left|1-p -\int_{\frac{(q_{p; n\alpha, \beta})-n\alpha\beta}{\beta\sqrt{n\alpha} }}^\infty \frac{1}{\sqrt{2\pi}}e^{-t^2/2} dt\right| &\leq \frac{\epsilon}{4\sqrt{2\pi}}\left|\int_{\frac{(q_{p; n\alpha, \beta})-n\alpha\beta}{\beta\sqrt{n\alpha} }}^\infty e^{-t^2/2} dt (1+O(n^{-1/2}))\right| \\
			&\leq\frac{\epsilon}{2\sqrt{2\pi}}\int_\mathbb{R}e^{-t^2/2} dt (1+O(n^{-1/2})) < \epsilon/2,
		\end{align*}
		where the last inequality holds for sufficiently large $n$, which shows that the difference between $\P(T_n > q_{p; n\alpha, \beta})$ and $1-p$ can be made arbitrarily small for any i.i.d sum, proving the first statement of the theorem. 
		
		To check the convergence in distribution, one has that as $G(\cdot; n\alpha, \beta)$   is strictly increasing in its argument,
		$$\P(T_n \leq q_{p; n\alpha, \beta})=\P(G(T_n; n\alpha,\beta) \leq G(q_{p; n\alpha, \beta};  n\alpha,\beta)) =\P(G(T_n; n\alpha,\beta) \leq p), $$
		since $ G(q_{p; n\alpha, \beta}, n\alpha,\beta))=p$, for every $p \in (0,1)$. The previous result implies that 
		$$\lim_{n\to\infty} \P(G(T_n; n\alpha,\beta) \leq p) \to p.$$ 
		That is, for any given $\epsilon >0$, with sufficiently large $n$, $|\P(G(T_n; n\alpha,\beta) \leq p)-p|\leq \epsilon$. Therefore, the cumulative distribution function $G(T_n; n\alpha,\beta)$ converges to the cumulative distribution function of Uniform(0, 1).
	\end{proof}
	
	\bigskip
	{\bf Remarks}: 
	\begin{itemize}
		\item[$\bullet$]As the distribution of $Y\sim Gamma(n\alpha, \beta)$ is the same as the distribution of $\sum_{i=1}^nY_i$ with $Y_i \sim Gamma(\alpha, \beta)$ a sequence of independent and identically distributed random variables, the Central Limit Theorem \citep[Theorem 27.1]{Bill86} guarantees that as $n$ increases, $n^{-1/2}(S_n-Gamma(4n/\nu,\nu/2))$ converges in distribution to $\mathcal{N}(0,2\nu)$ and similarly $n^{-1/2}(\tilde{S}_n-Gamma(nm^2/v,v/m))$ converges in distribution to $\mathcal{N}(0,2v)$. This convergence can provide standard bounds on the approximation error like the Chebyshev bound $\P(|S_n-Gamma(4n/\nu,\nu/2)|>\sqrt{n}\varepsilon)\leq 2\nu\varepsilon^{-2}$ and $\P(|S_n-Gamma(nm^2/v,v/m)|>\sqrt{n}\varepsilon)\leq 2v\varepsilon^{-2}$. 
		
		\item[$\bullet$] Convergence in distribution further justifies the proposed approximation in the sense that the Wasserstein distance between the normalized combination statistic and its approximiation becomes arbitrarily small for large enough $n$. This is, as the Wasserstein distance metrizes weak convergence of random variables \citep{Clement2008WD}, one has that for $\{X_j\}_{j \in \mathbb{N}}$ a sequence of independent and identically distributed non-negative discrete random variables with mean $\mu$, variance $\sigma^2$ and $\E(X_j^4)<\infty$, $T_n=\sum_{j=1}^nX_j$ and $Y_n \sim Gamma(n\mu^2/\sigma^2, \sigma^2/\mu)$ one has that
		\begin{equation}
			\label{eq:distto0}
			\lim_{n \rightarrow \infty} W_2(n^{-1/2}T_n,n^{-1/2}Y_n)=0,
		\end{equation}
		a convergence result that applies to the combination statistics in \eqref{eq:sum12} and their gamma approximations in \eqref{eq:approxSn}. 
		
		\item[$\bullet$]The result in \eqref{eq:distto0} provides a rate $n^{-1/2}$ at which $n^{-1/2}S_n$ and  $n^{-1/2}\tilde{S}_n$ asymptotically resemble $n^{-1/2}Gamma(4n/\nu, \nu/2)$ and $n^{-1/2}Gamma(m^2n/v, v/m)$. Replacing $n^{-1/2}$ with $n^{-q}$, any rate of convergence of order $q<0.5$ is unachievable as a corollary of Slutsky's Theorem \citep[Chapter 7.3]{grimmett2001probability}, which guarantees that for $q<1/2$ the statistic $n^{-q}(T_n-Y_n)$ diverges, yielding $0.5$ to be the optimal convergence rate. The problem of convergence rate of a sum statistic has been studied in other contexts like empirical distributions \citep{GKconv} in optimal transport for sequences converging to a fixed probability distribution \citep{aude2016stochastic}, where convergence rates of order $q<0.5$ have similarly been proven unachievable.

	\end{itemize}
	\subsection{Proof of Proposition \ref{thm:asymptoticpvaluenoniid}}
	\begin{proof}
		The conditions on the statement of the theorem are the Lyapunov Conditions for non-i.i.d. sums $S_n$ and $\tilde{S}_n$.  We prove the result for $S_n$, for $\tilde{S}_n$ it can be proven analogously.
		
		The Lyapunov condition guarantees (see \cite{Bill86}) that, as $n \to \infty$, for any $p \in (0,1)$
		$$\left|\P(S_n > q_{p; 4n/\bar{\nu},\bar{\nu}/2})-\int_{\frac{(q_{p; 4n/\bar{\nu},\bar{\nu}/2})-2n}{\sqrt{n/\bar{\nu}} }}^\infty \frac{1}{\sqrt{2\pi}}e^{-t^2/2} dt\right|\to 0.$$
		
		For fixed $\epsilon >0$, one can find $n_0$ such that the above quantity is smaller than $\epsilon/2$ and, using Stirling's approximation in similar fashion than in the proof of Theorem \ref{thm:asymptoticpvalue}
		
		\begin{align*}
			& \left|1-p -\int_{\frac{(q_{p; 4n/\bar{\nu},\bar{\nu}/2})-2n}{\sqrt{n/\bar{\nu}} }}^\infty \frac{1}{\sqrt{2\pi}}e^{-t^2/2} dt\right|\\  \leq & \frac{\epsilon}{4\sqrt{2\pi}}\left|\int_{\frac{(q_{p; 4n/\bar{\nu},\bar{\nu}/2})-2n}{\sqrt{n/\bar{\nu}} }}^\infty e^{-t^2/2} dt (1+O(n^{-1/2}))\right| \\
			< &\epsilon/2 \mbox{ for large enough }n.
		\end{align*}
		Putting these two inequalities together, one achieves

		$$|\P(S_n \leq q_{p; 4n/\bar{\nu},\bar{\nu}/2})-p|\leq \epsilon.$$
		
		As $\epsilon >0$ is arbitrary, the result follows.  
	\end{proof}
	
	\bigskip
	\textbf{Remark:} 
	
	Observe that the regularity conditions in the theorem are essentially the Lyapunov conditions described in \cite{Bill86} (page 362). They guarantee that the normalized sums converge to a standard normal in distribution provided that the sequences of adjusted p-values is not affected by too large or too small increments at a time, which is attained by controlling the $2+\delta$ central moment. 
	
	As one example, any sequence of i.i.d. random variables satisfy that for arbitrary $\delta>0$,
	$$\lim_{n \to \infty} \frac{\sum_{j=1}^n\E(|X_j-\E(X_j)|^{2+\delta})}{(\sum_{j=1}^n \Var(X_j) )^{1+\delta/2}}=\lim_{n \to \infty}\frac{n\E(|X_1-\E(X_1)|^{2+\delta})}{n^{1+\delta/2}\Var(X_j)^{1+\delta/2}}=0.$$
	
	As another example, consider $\ddot{Z}_j$ and $\tilde{Z}_j$ were obtained from Binomial$(1,p_j)$, i.e., $F^{(j)}_0=0, F^{(j)}_1=1-p_j, F^{(j)}_2=1$. For $S_n = \sum_{j=1}^n \ddot{Z}_j$, we have 
	$$
	\Var(\ddot{Z}_j) = \nu_j=4\frac{(1-p_j)\log^2 (1-p_j)}{p_j}
	$$ and 
	$$
	\E(|\ddot{Z}_j-2|^{2+\delta})=(2|\log (1-p_j)|)^{2+\delta}\left[1-p_j+\left(\frac{1-p_j}{p_j}\right)^{2+\delta}\right],
	$$ 
	and therefore  
	$$
	\frac{\sum_{j=1}^n\E(|\ddot{Z}_j-2|^{2+\delta})}{(\sum_{j=1}^n \nu_j )^{1+\delta/2}}= \frac{\sum_{j=1}^n |\log (1-p_j)|^{2+\delta}\left[1-p_j+\left(\frac{1-p_j}{p_j}\right)^{2+\delta}\right]}{\left(\sum_{j=1}^n \frac{(1-p_j)\log^2 (1-p_j)}{p_j} \right)^{1+\delta/2}}
	$$ 
	will converge to 0 if and only if $\Var (S_n)=\sum_{j=1}^n \frac{(1-p_j)\log^2 (1-p_j)}{p_j} \to \infty$. 
	An analogous result can be obtained for $\tilde{S}_n = \sum_{j=1}^n \tilde{Z}_j$.

	\section{Further numerical Studies}\label{sec:appnumerical}
	
	\subsection{Combining independent and identically distributed p-values}
	
\begin{figure}[h!]
		\begin{centering}
			{\includegraphics[width=0.32\linewidth]{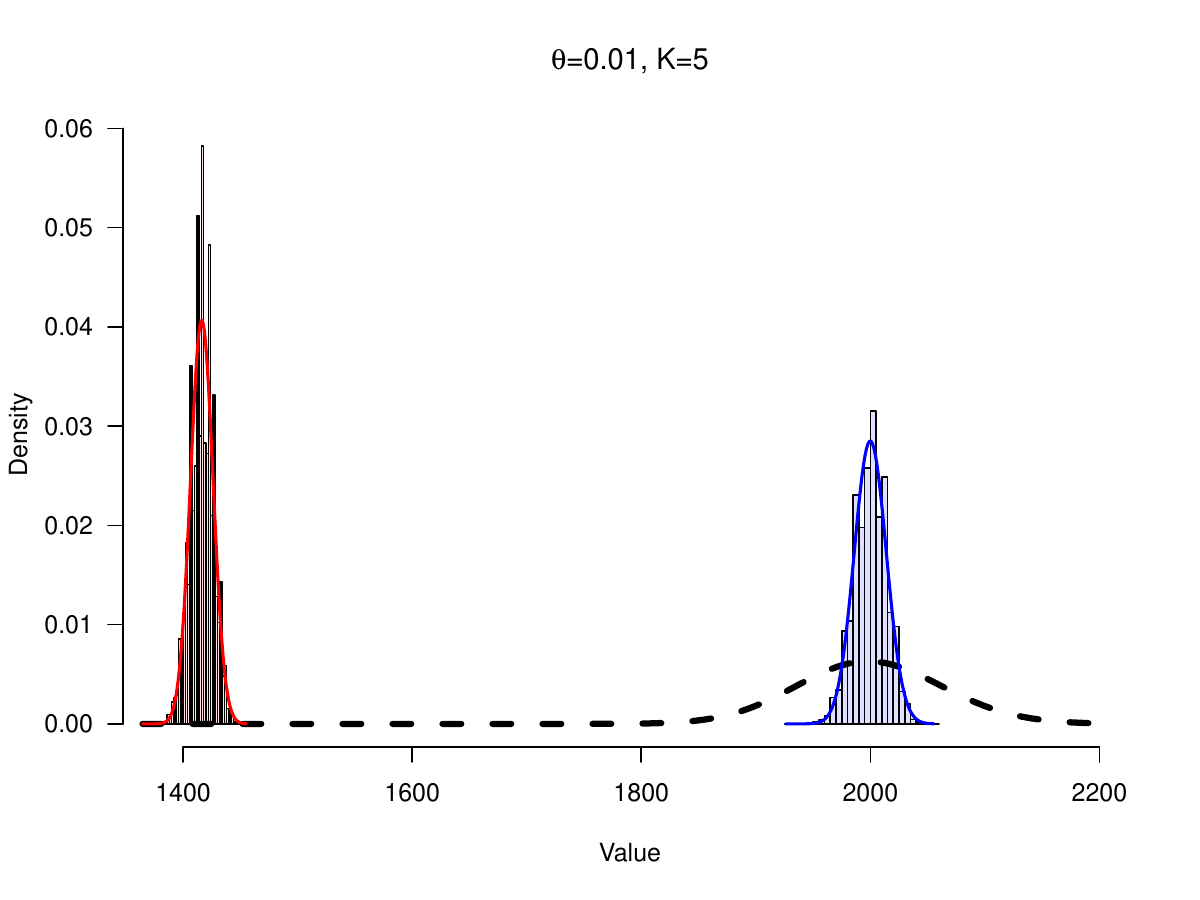}}
			{\includegraphics[width=0.32\linewidth]{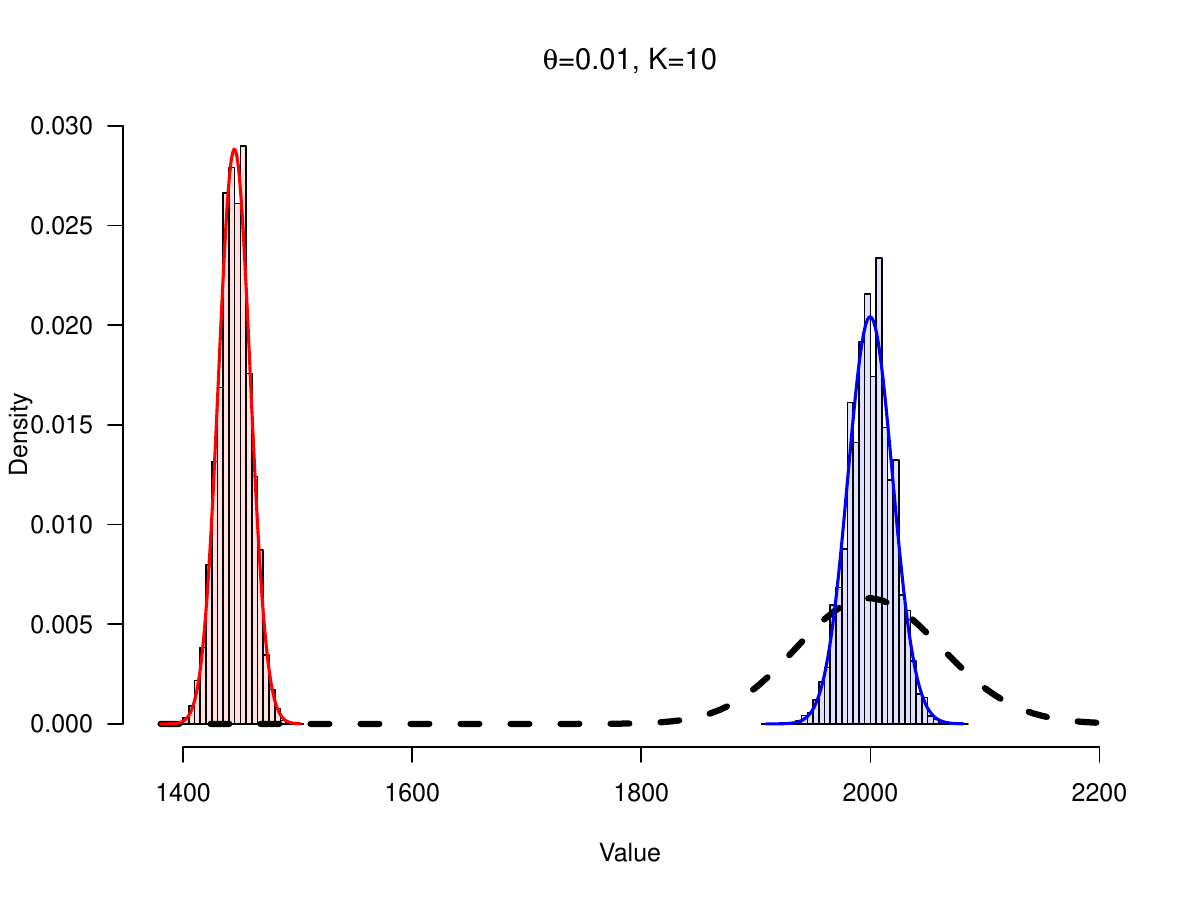}}
			{\includegraphics[width=0.32\linewidth]{histgofth001k20}}\\
			{\includegraphics[width=0.32\linewidth]{histgofth01k5}}
			{\includegraphics[width=0.32\linewidth]{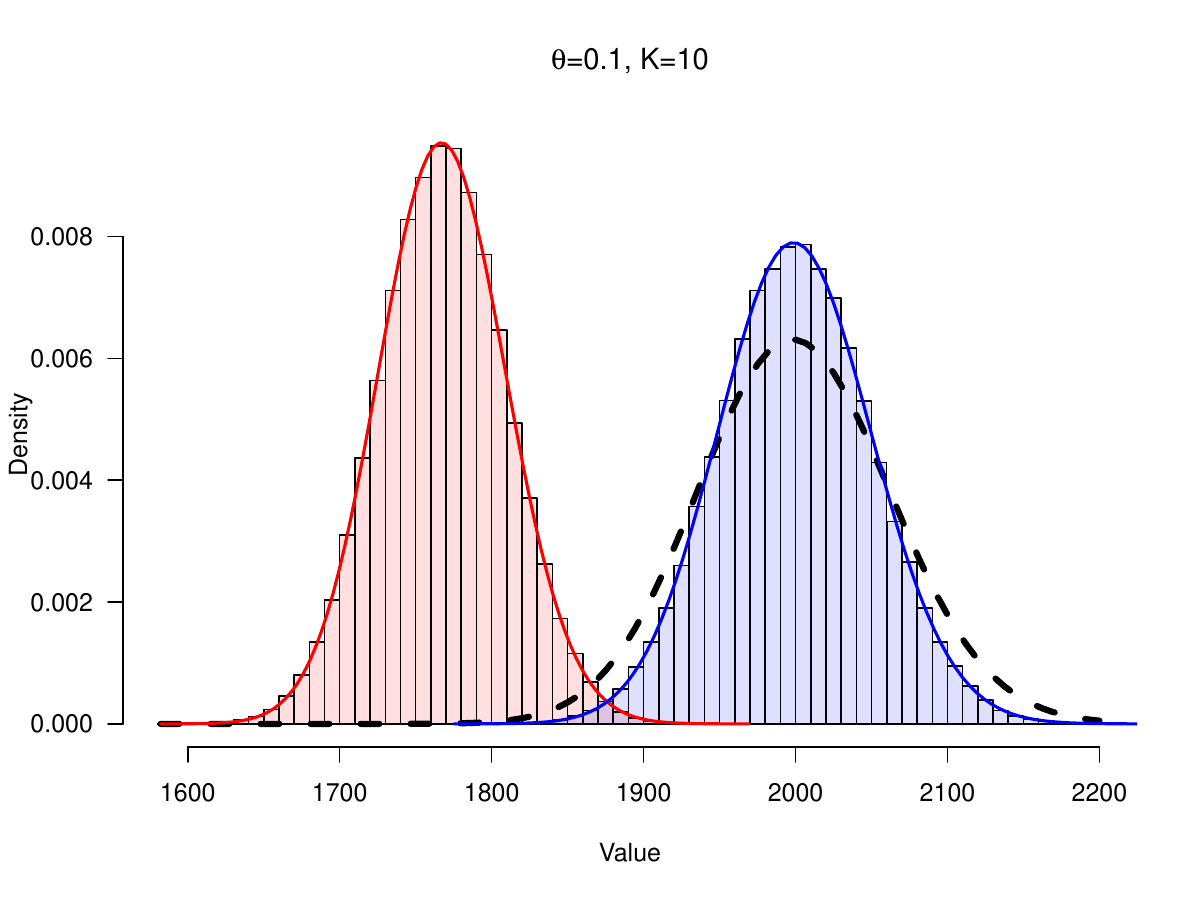}}
			{\includegraphics[width=0.32\linewidth]{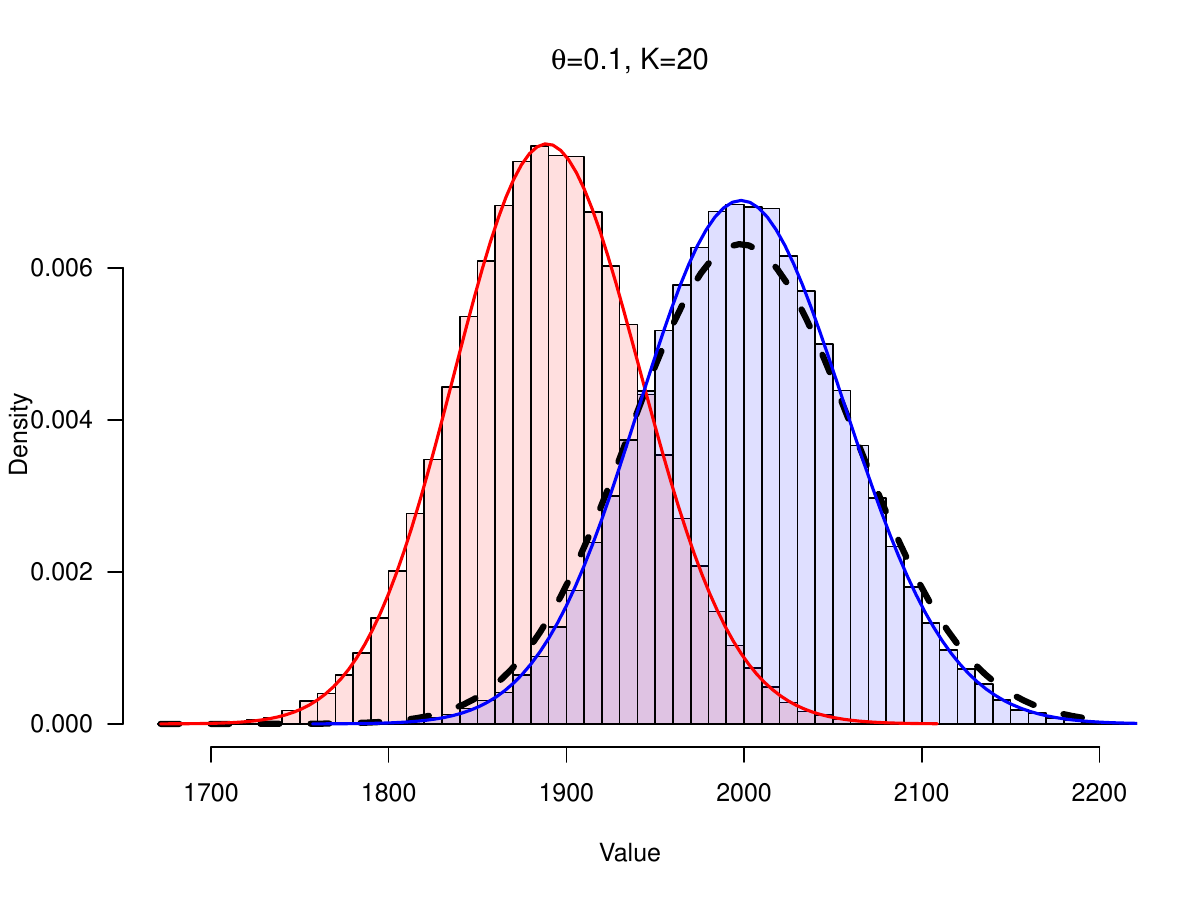}}\\
			{\includegraphics[width=0.32\linewidth]{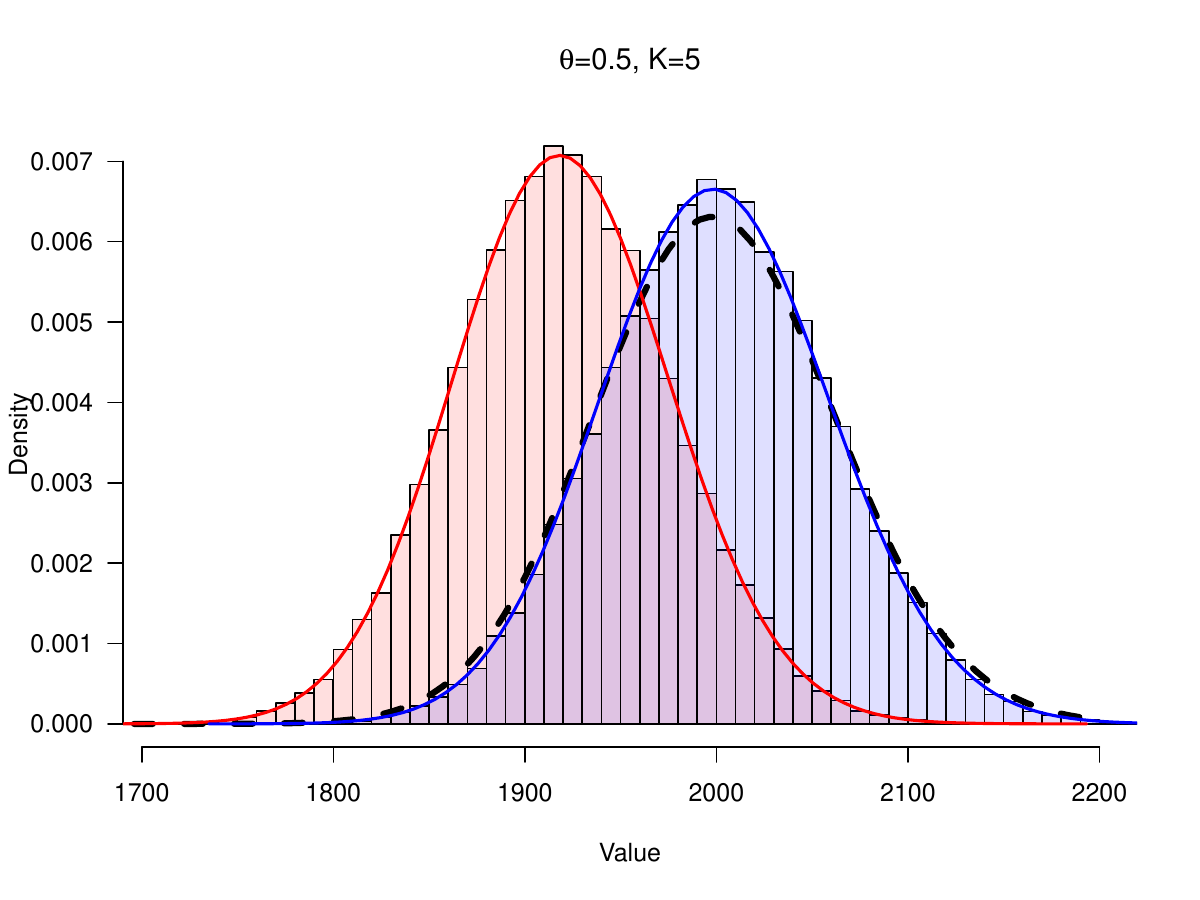}}
			{\includegraphics[width=0.32\linewidth]{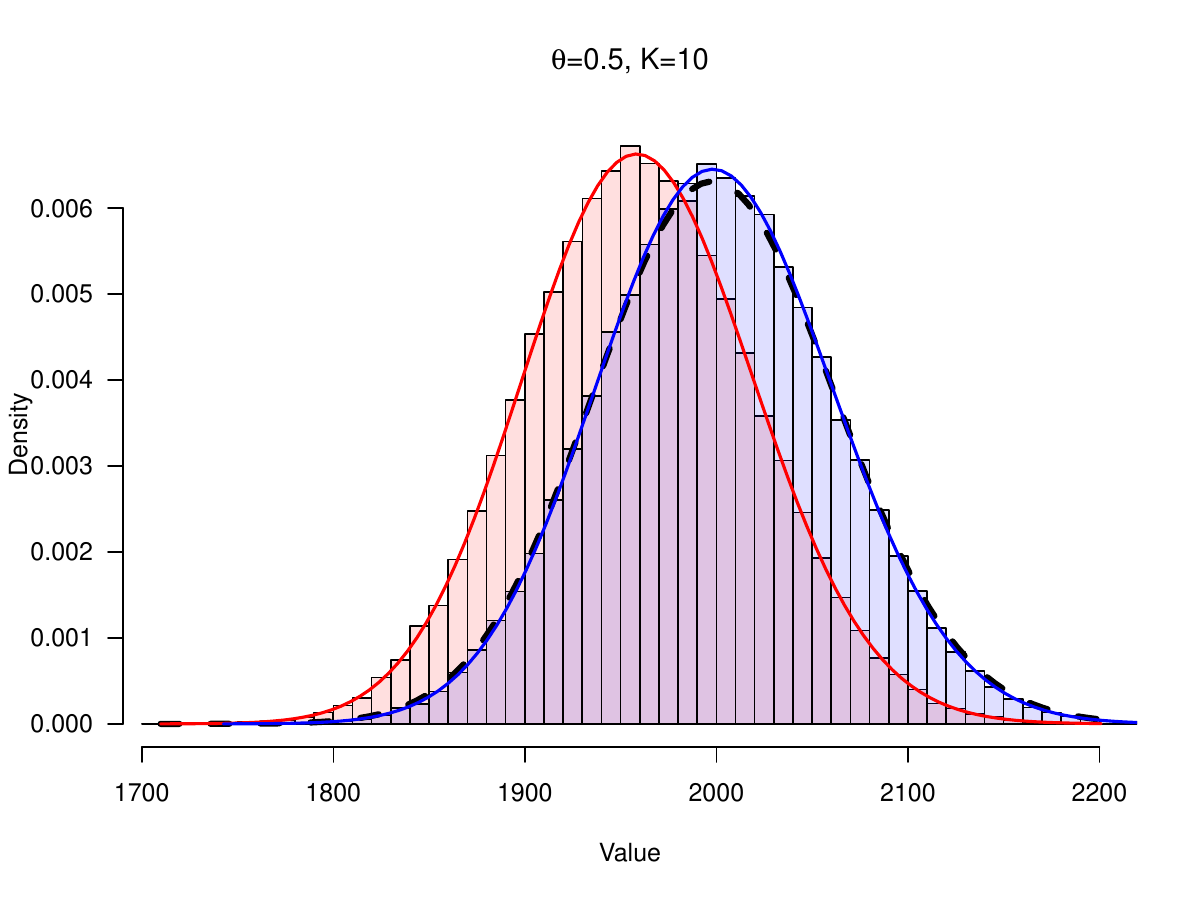}}
			{\includegraphics[width=0.32\linewidth]{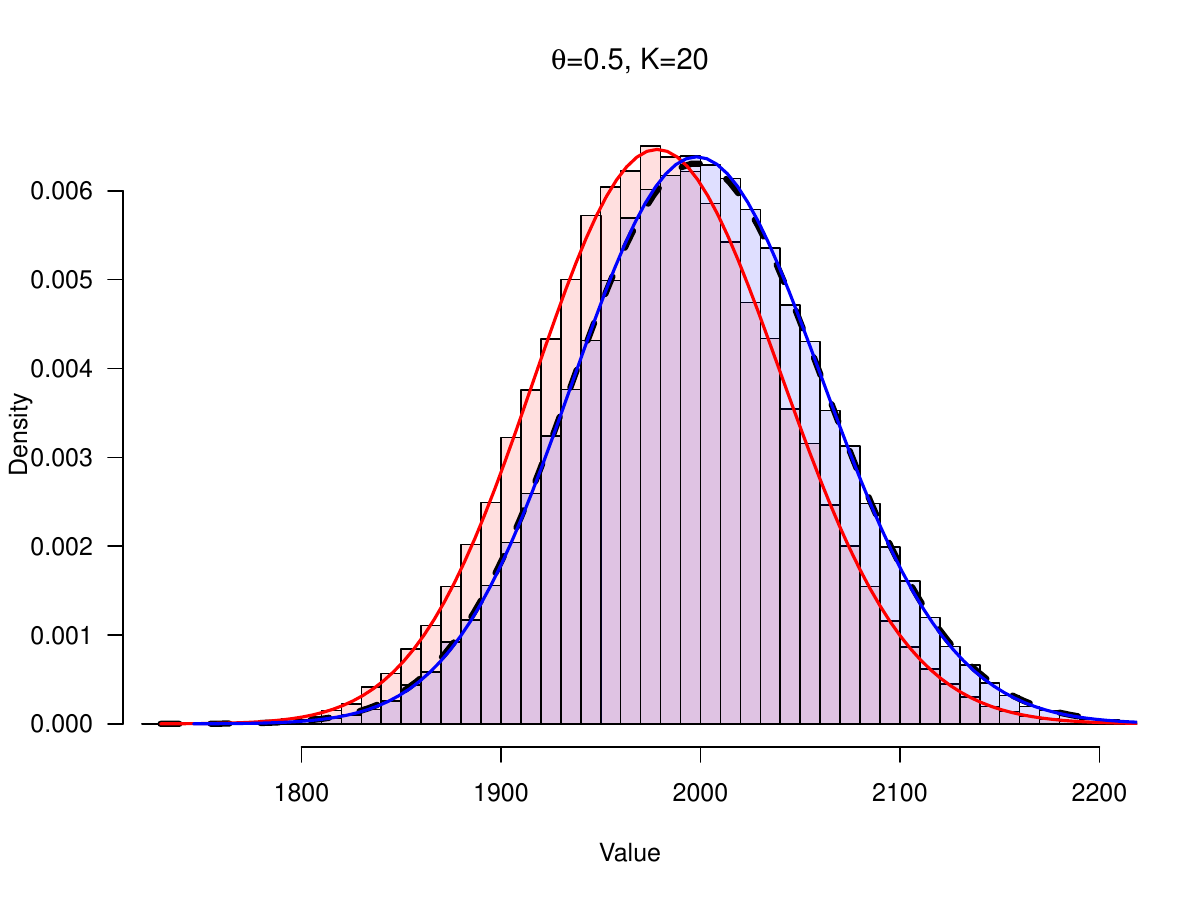}}
			\caption{Empirical distribution of $S_n$ (blue histogram, $10^7$ simulations) and $\tilde{S}_n$ (red histogram, $10^7$ simulations) compared with the densities of Gamma($4n/\nu,\nu/2$) (blue solid curve), Gamma($m^2n/v,v/m$) (red solid curve) and $\chi_{2n}^2$ (black dashed curve). Each panel corresponds to one combination of $\theta_0$ and $K$ values in Table \ref{table:paramfit}. $n=1000$. 
			}
			\label{fig:histgof}
			
		\end{centering}
	\end{figure}

\begin{figure}[h!]
		\begin{centering}
			{\includegraphics[width=0.32\linewidth, height=3.5cm]{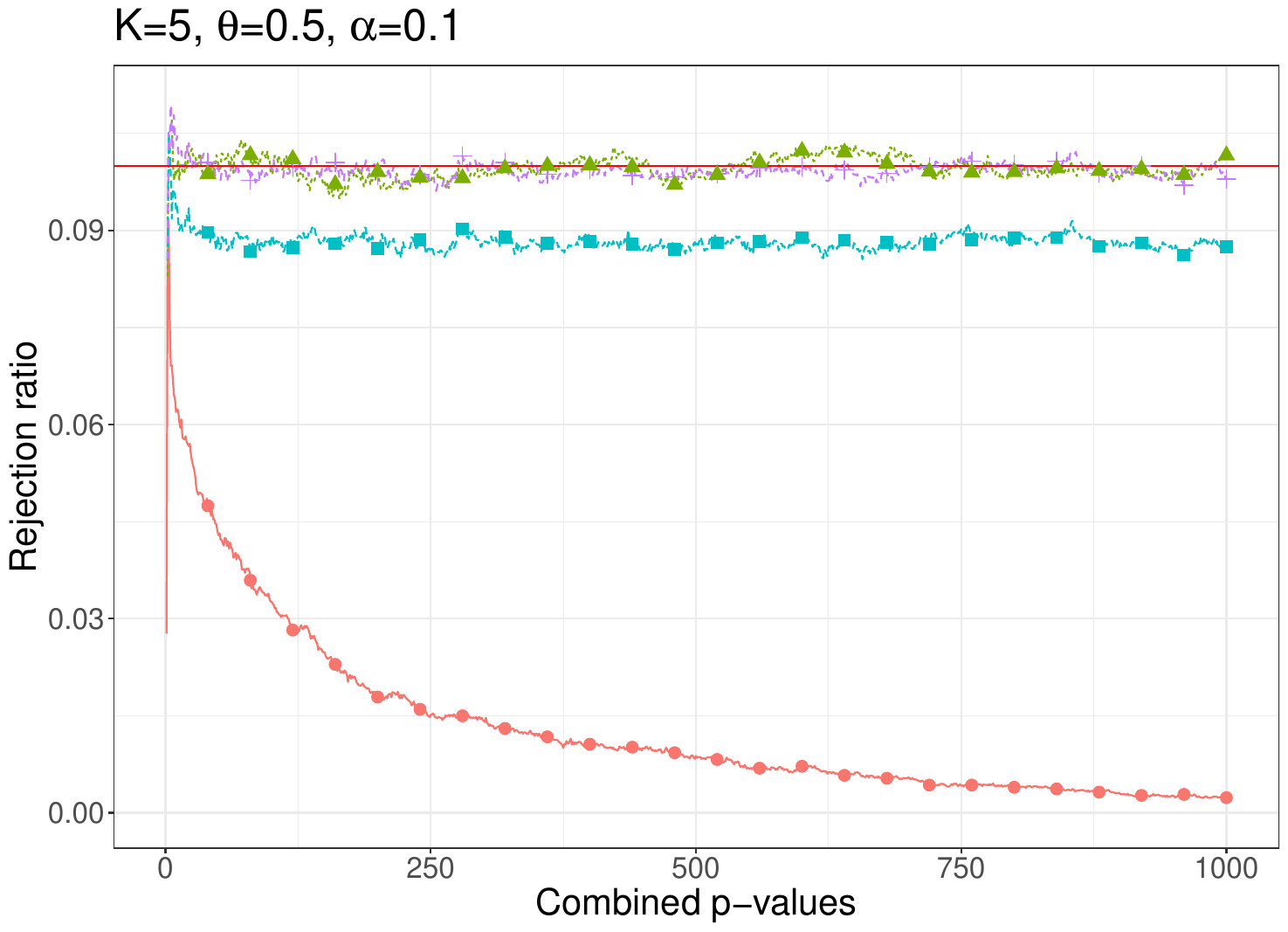}}
			{\includegraphics[width=0.32\linewidth, height=3.5cm]{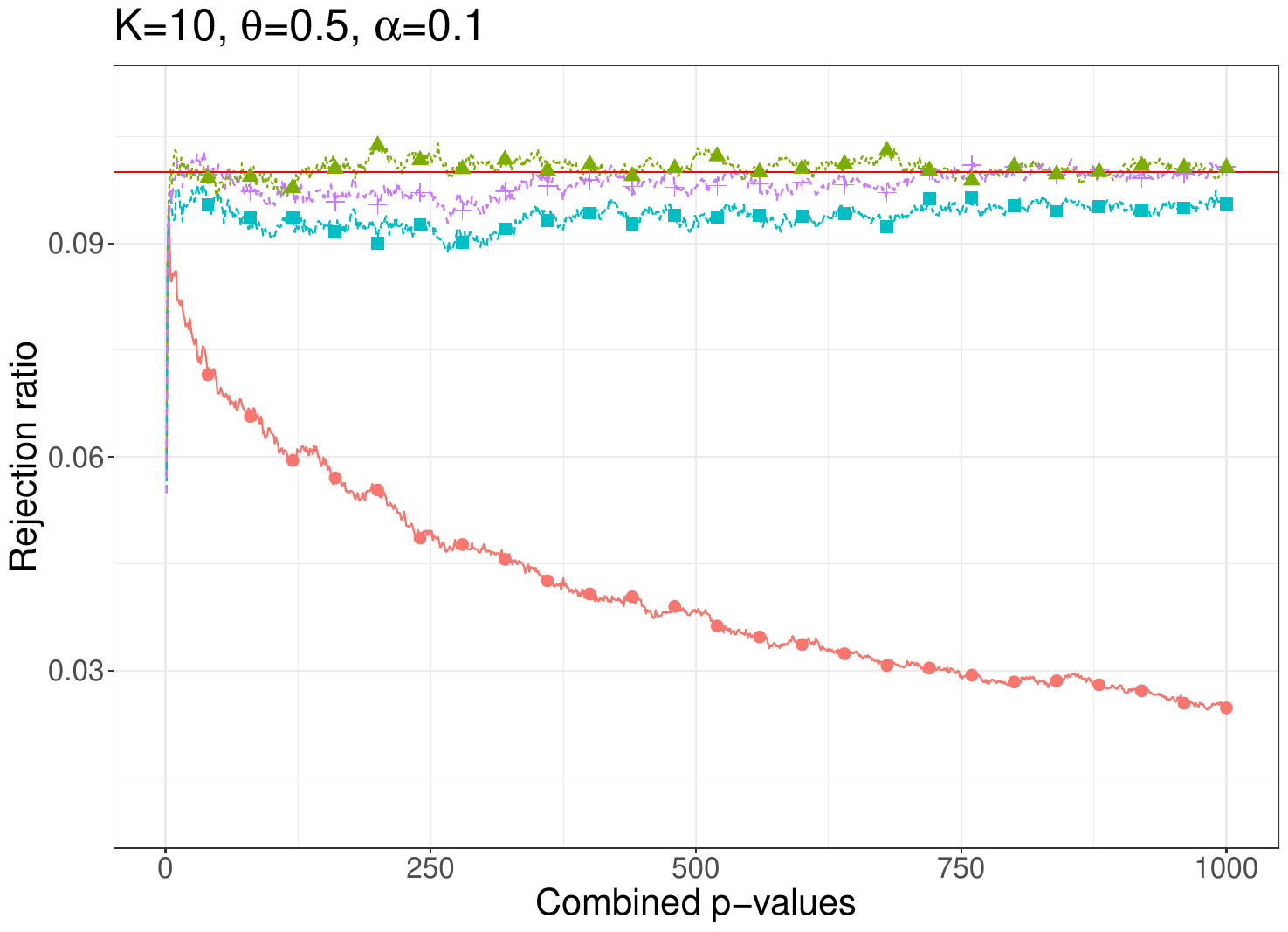}}
			{\includegraphics[width=0.32\linewidth, height=3.5cm]{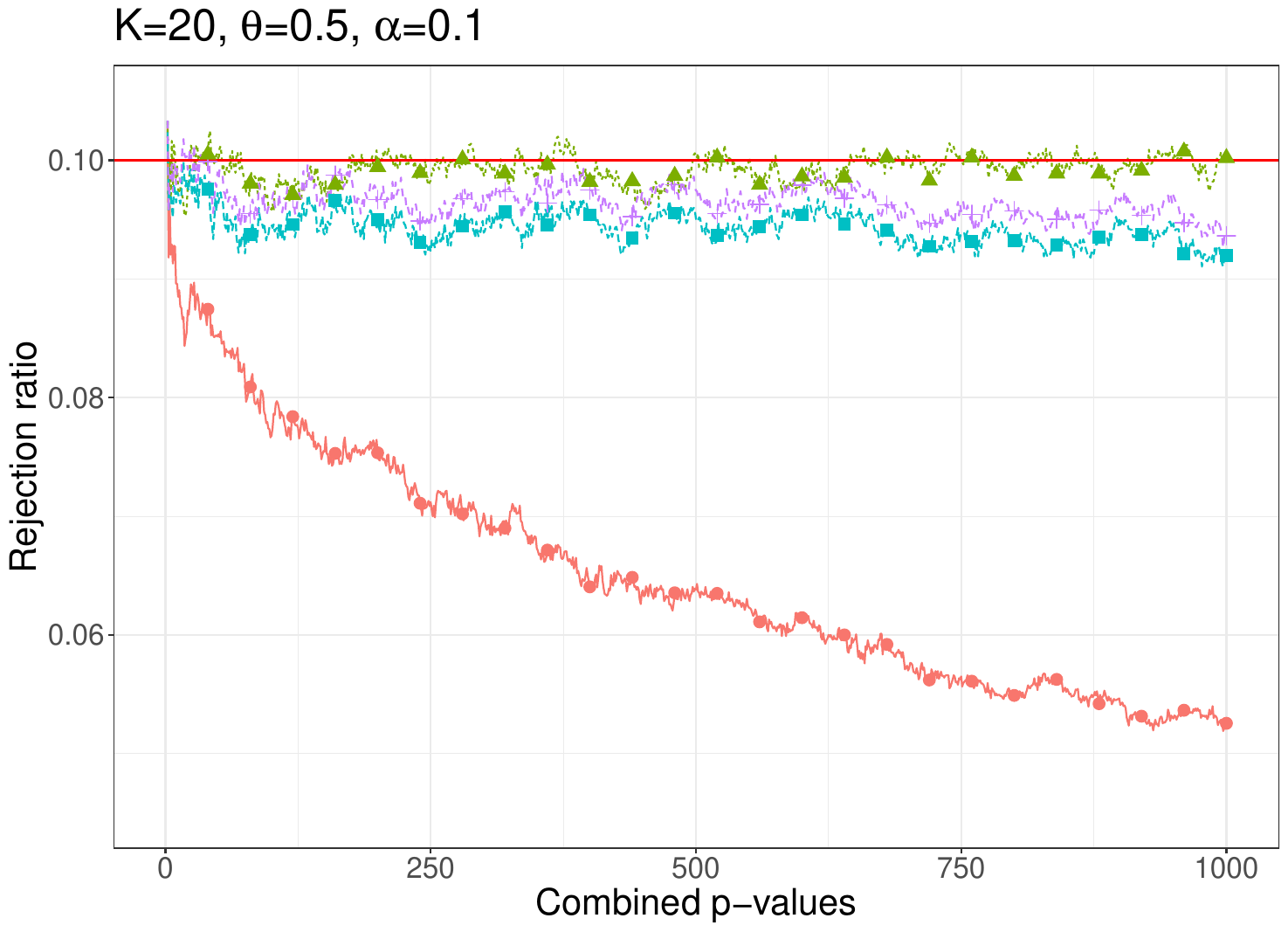}}\\
			{\includegraphics[width=0.32\linewidth, height=3.5cm]{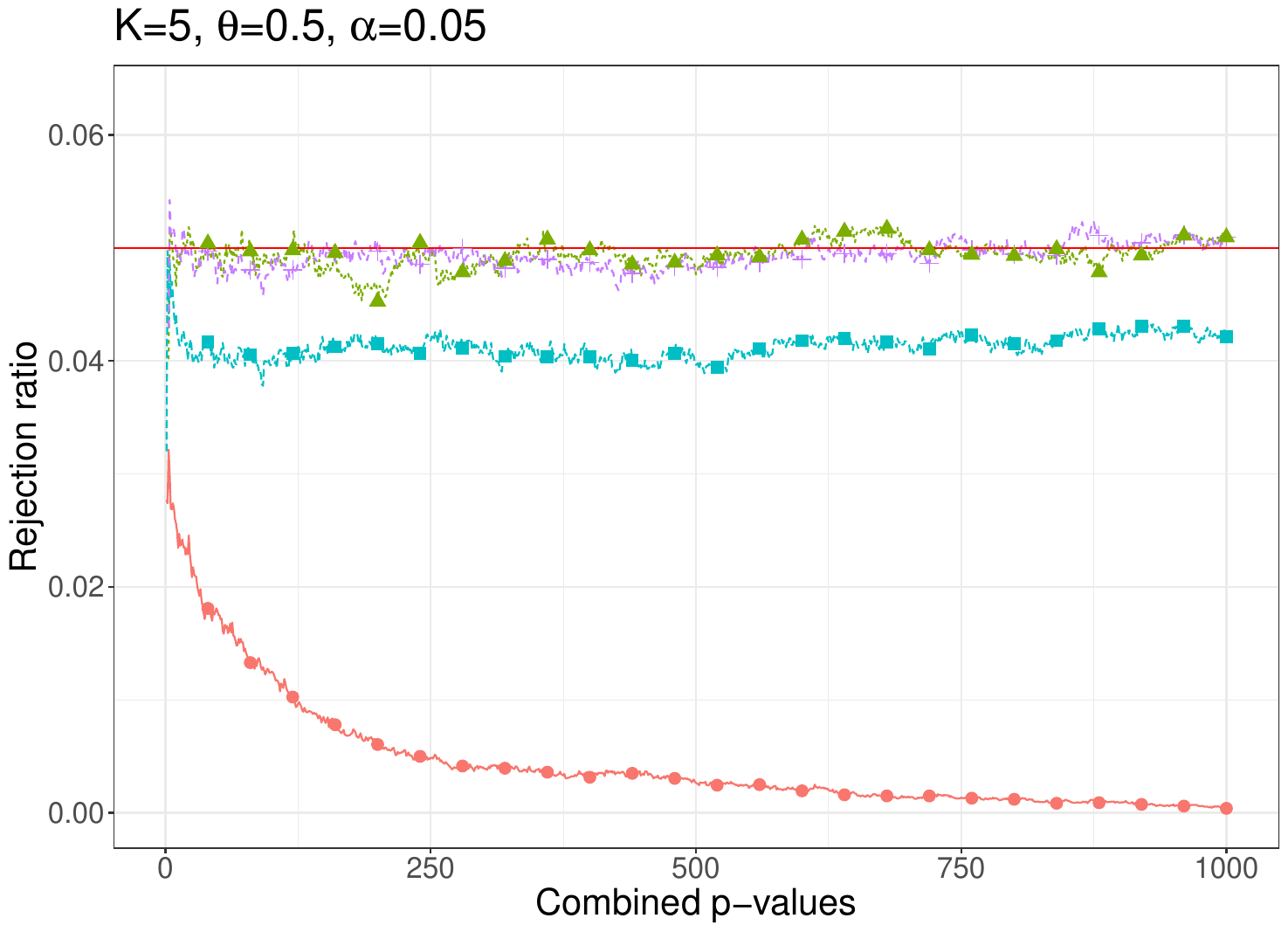}}
			{\includegraphics[width=0.32\linewidth, height=3.5cm]{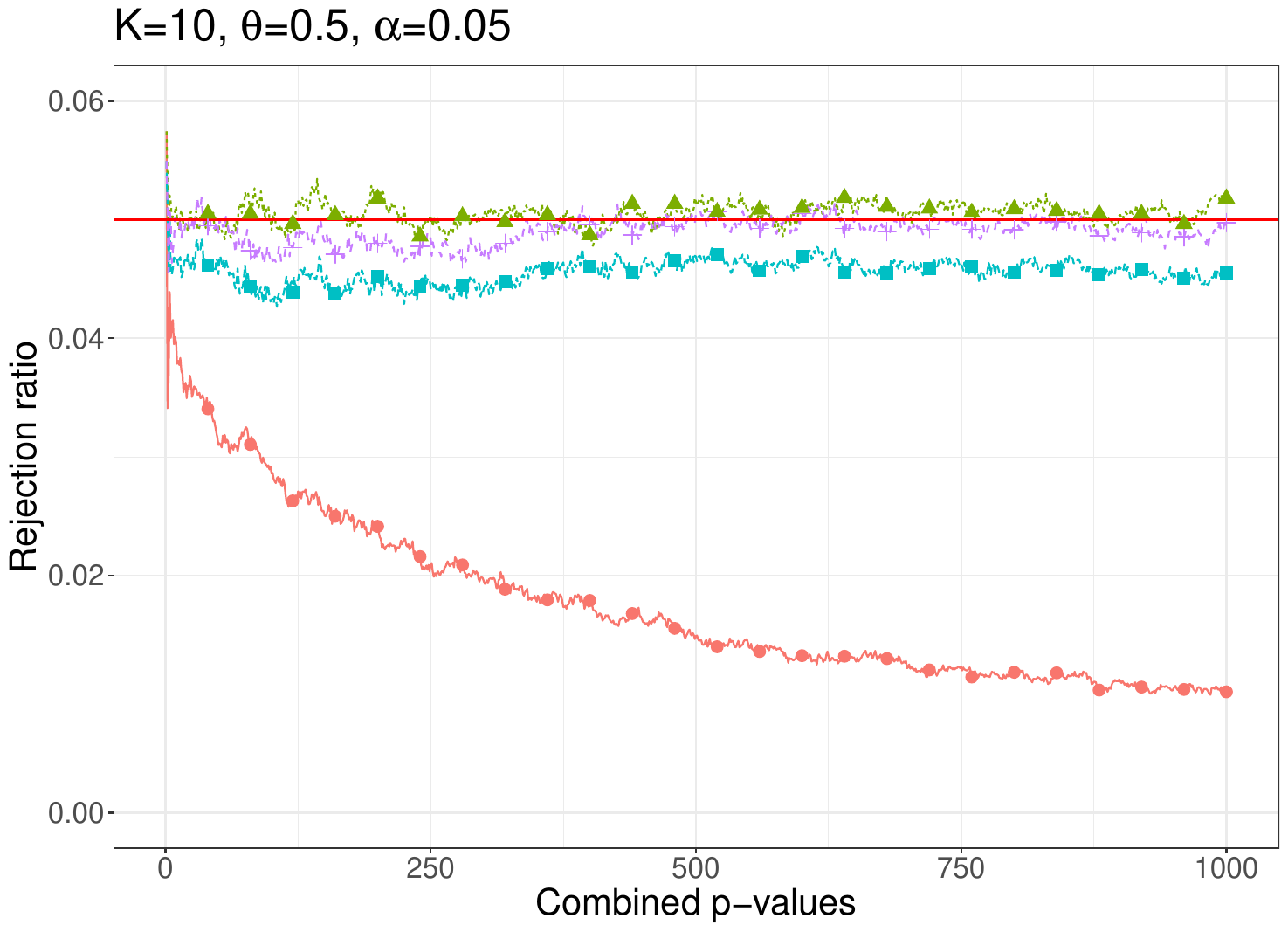}}
			{\includegraphics[width=0.32\linewidth, height=3.5cm]{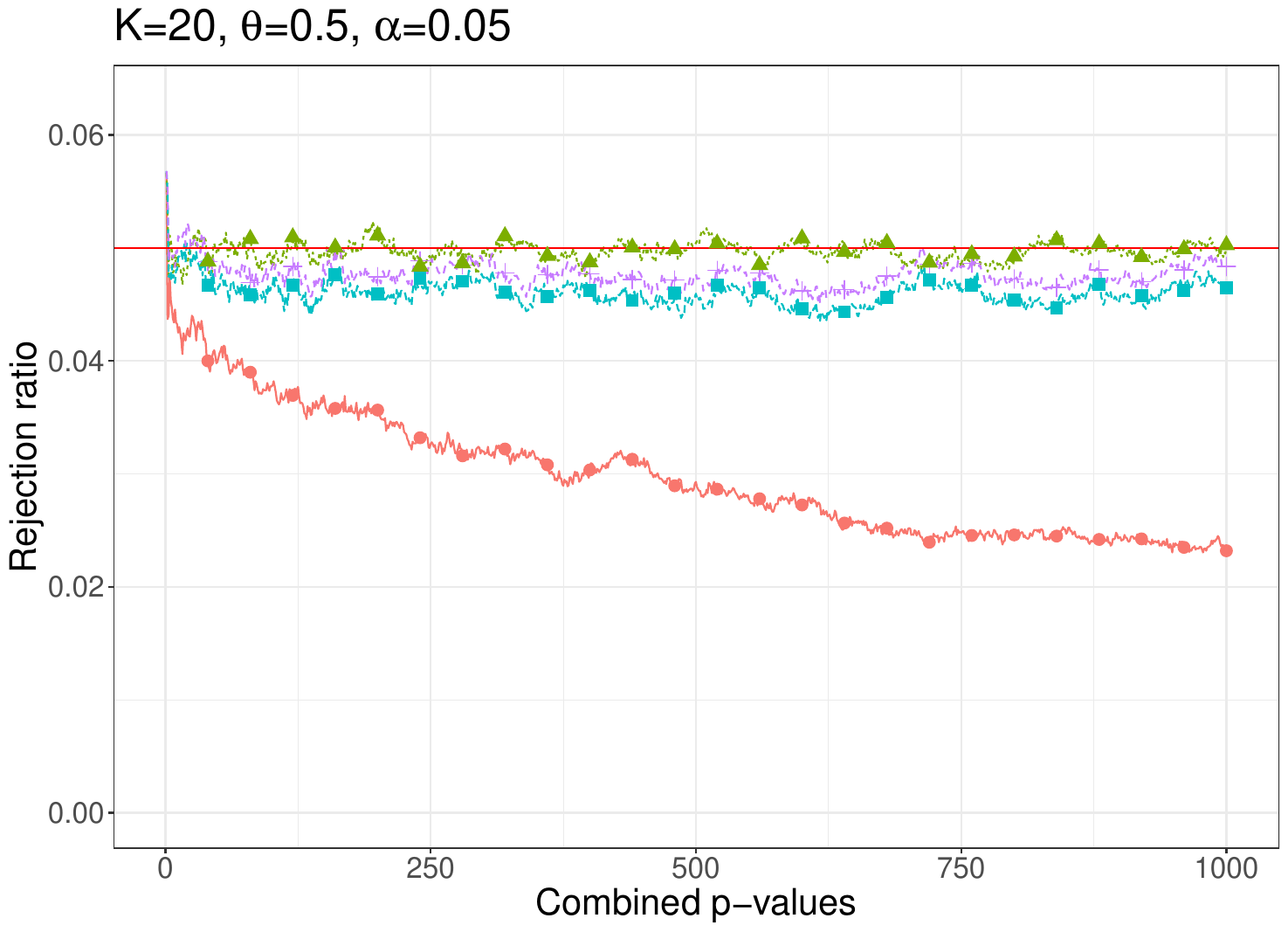}}\\
			{\includegraphics[width=0.32\linewidth, height=3.5cm]{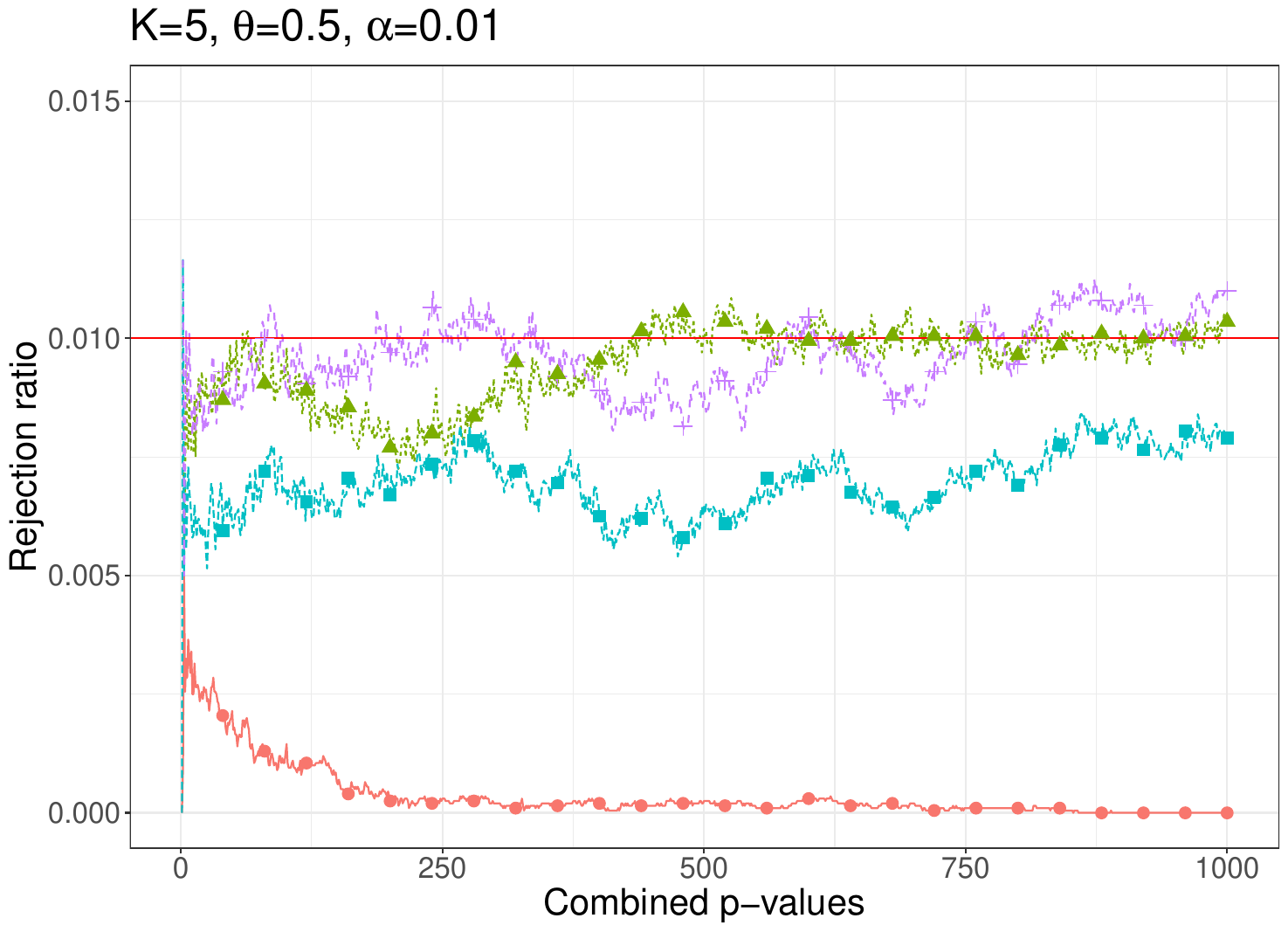}}
			{\includegraphics[width=0.32\linewidth, height=3.5cm]{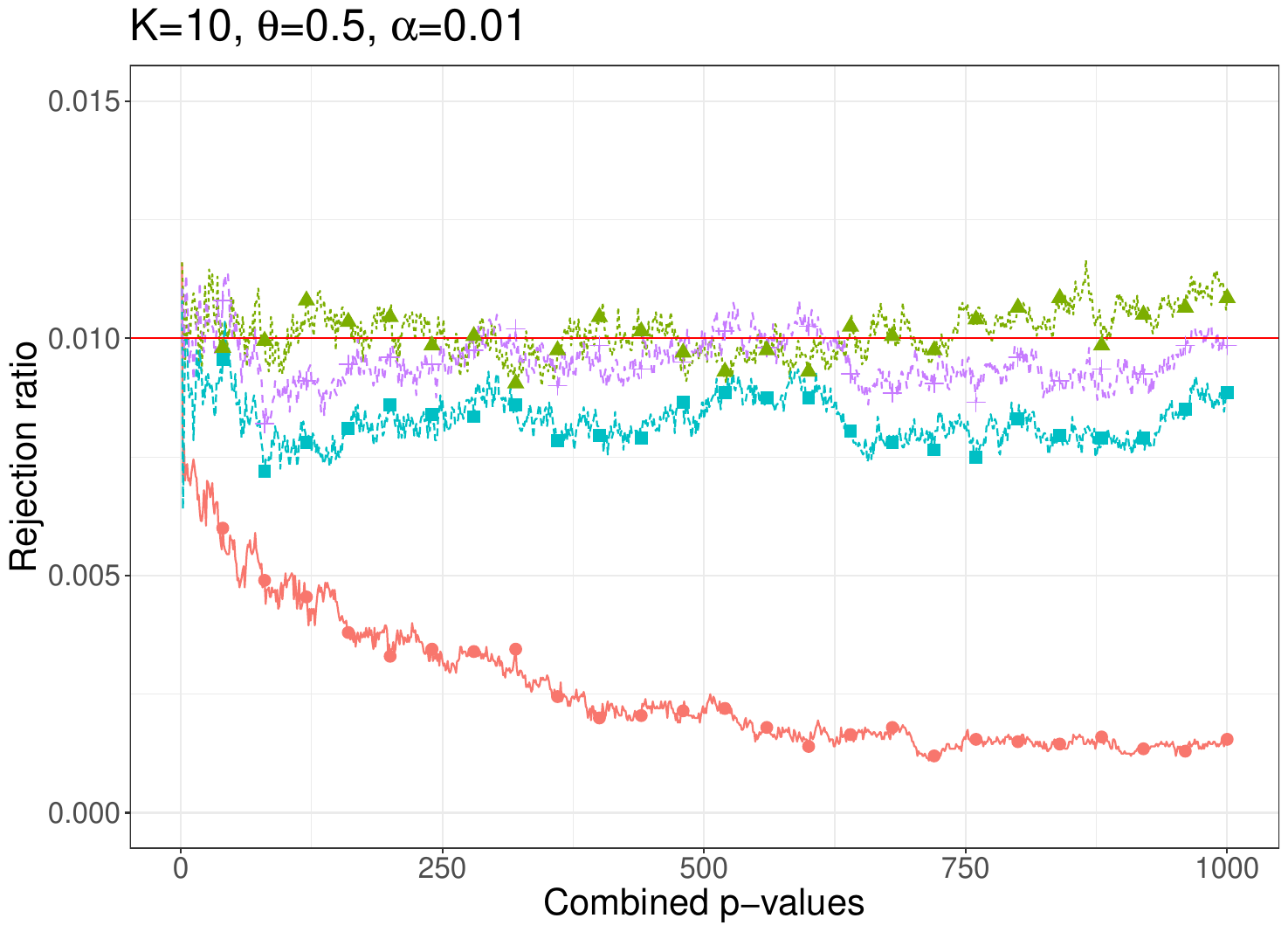}}
			{\includegraphics[width=0.32\linewidth, height=3.5cm]{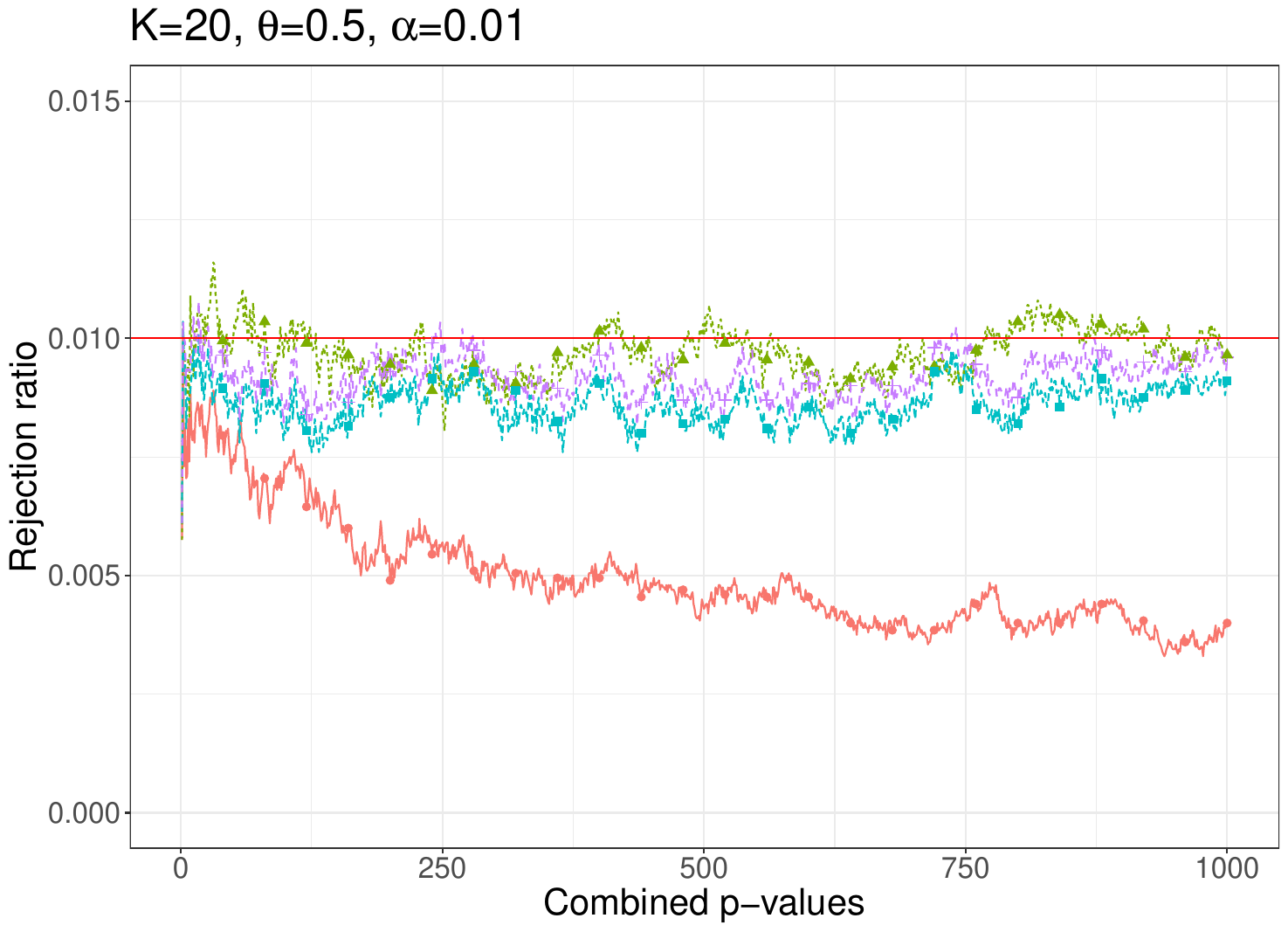}}\\
			{\includegraphics[width=0.32\linewidth, height=3.5cm]{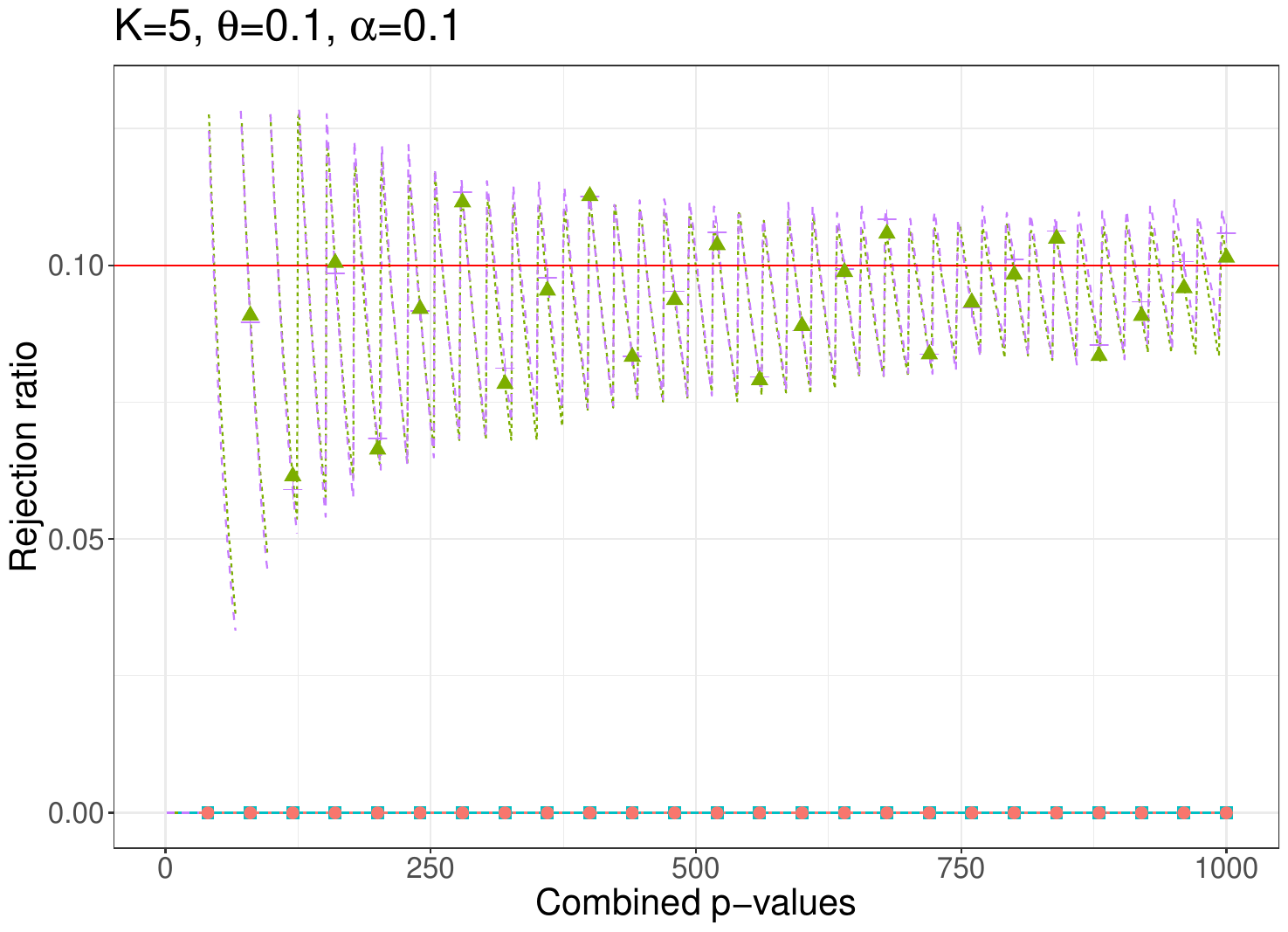}}
			{\includegraphics[width=0.32\linewidth, height=3.5cm]{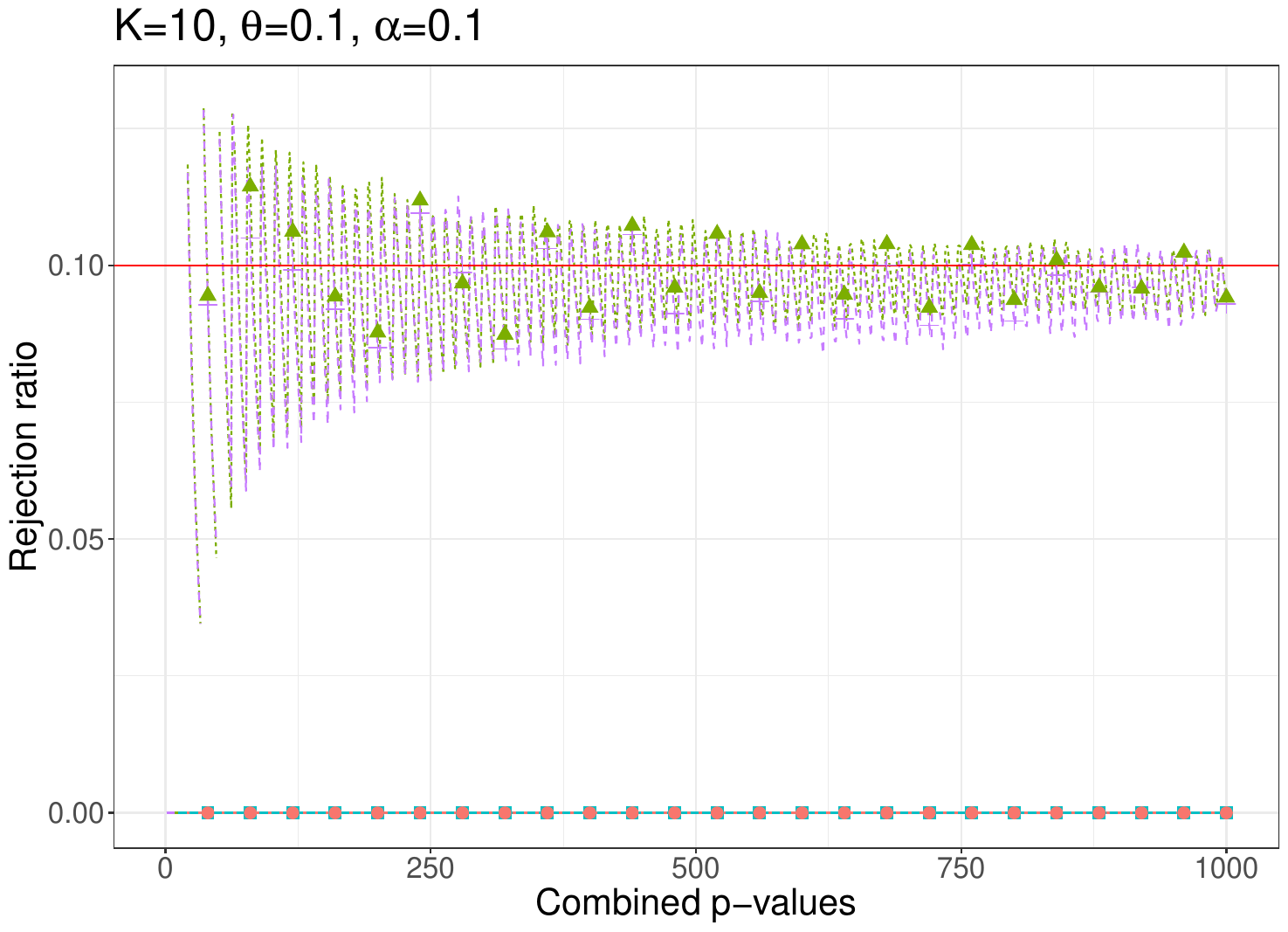}}
			{\includegraphics[width=0.32\linewidth, height=3.5cm]{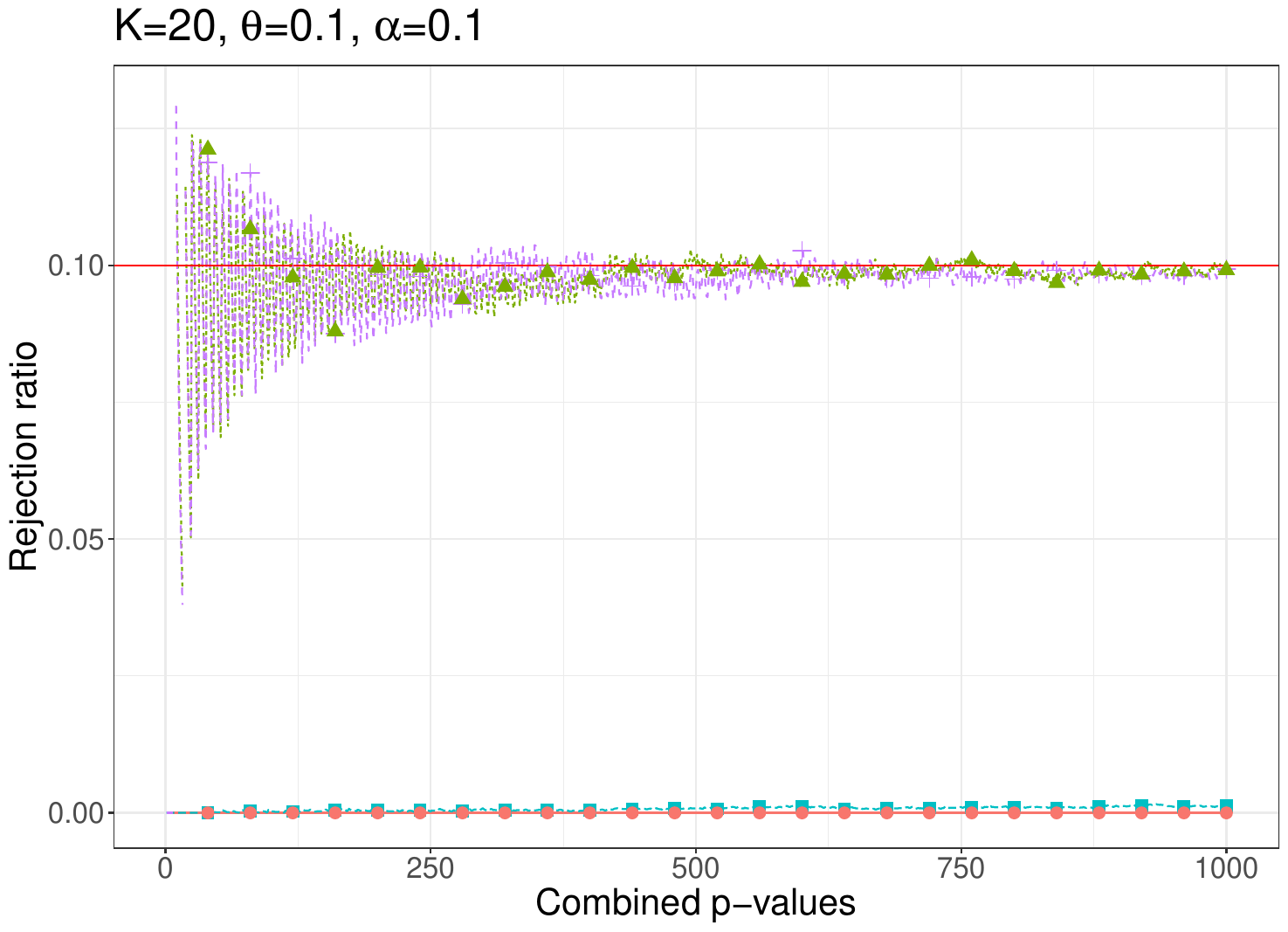}}\\
			{\includegraphics[width=0.32\linewidth, height=3.5cm]{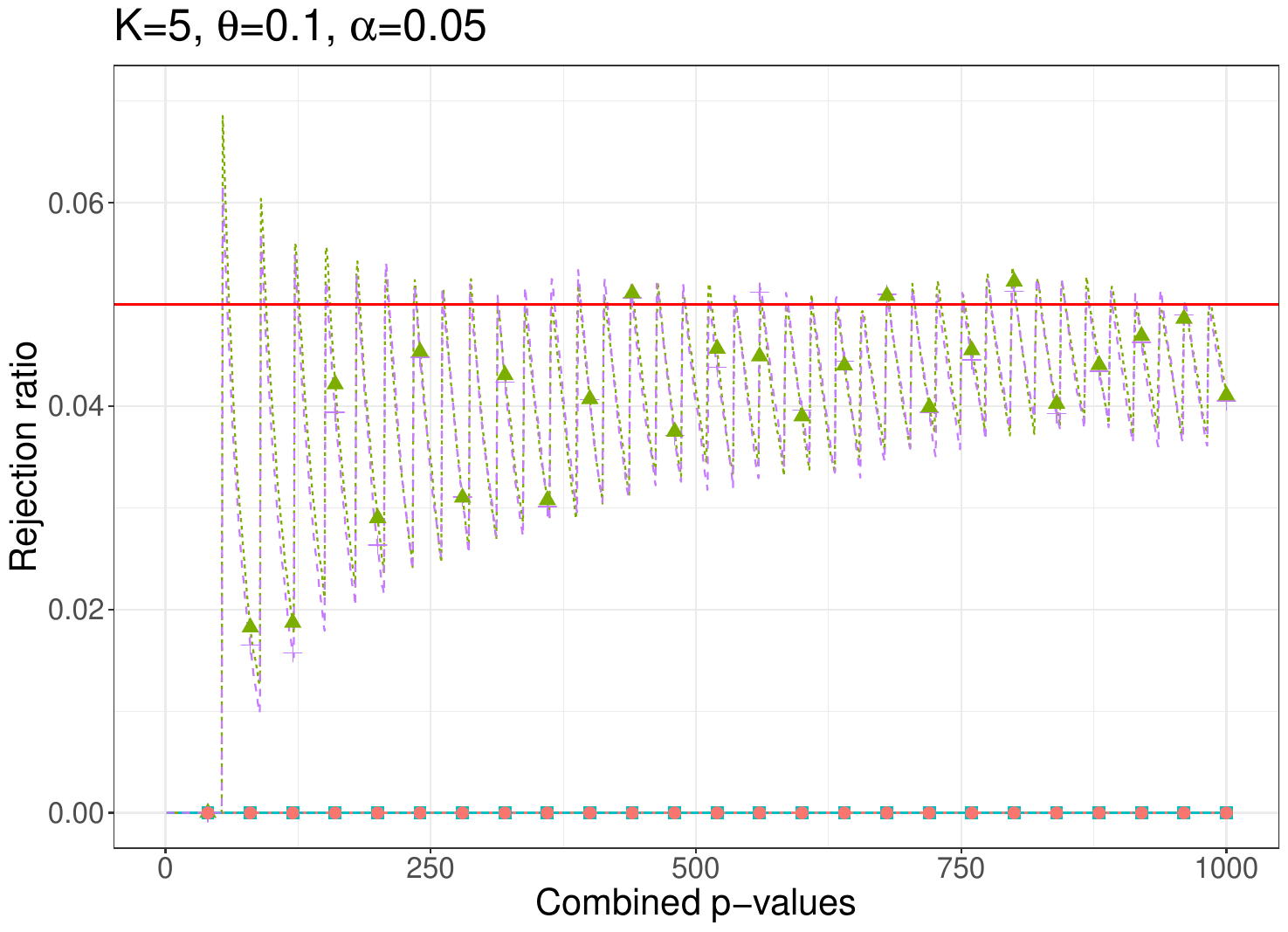}}
			{\includegraphics[width=0.32\linewidth, height=3.5cm]{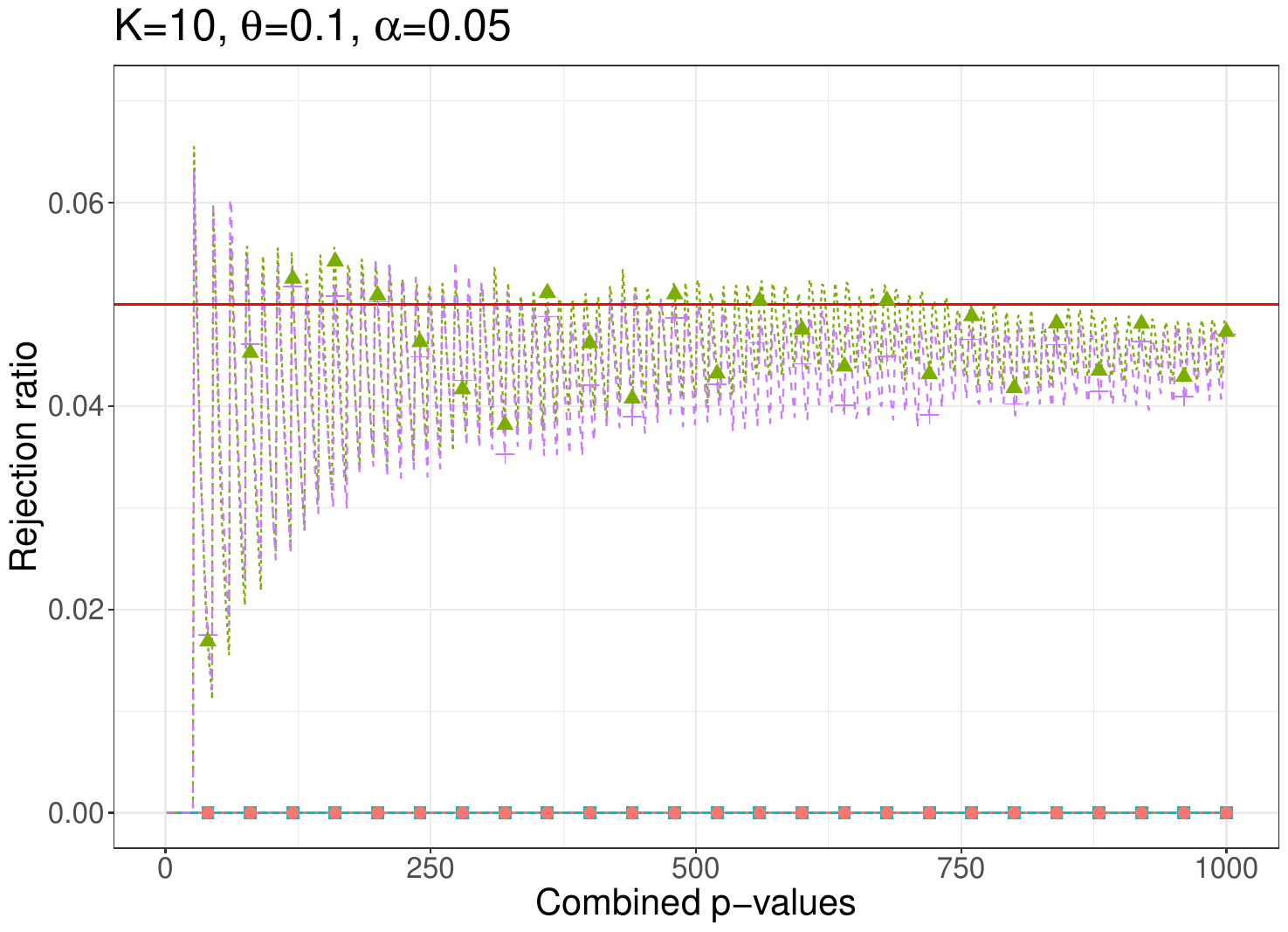}}
			{\includegraphics[width=0.32\linewidth, height=3.5cm]{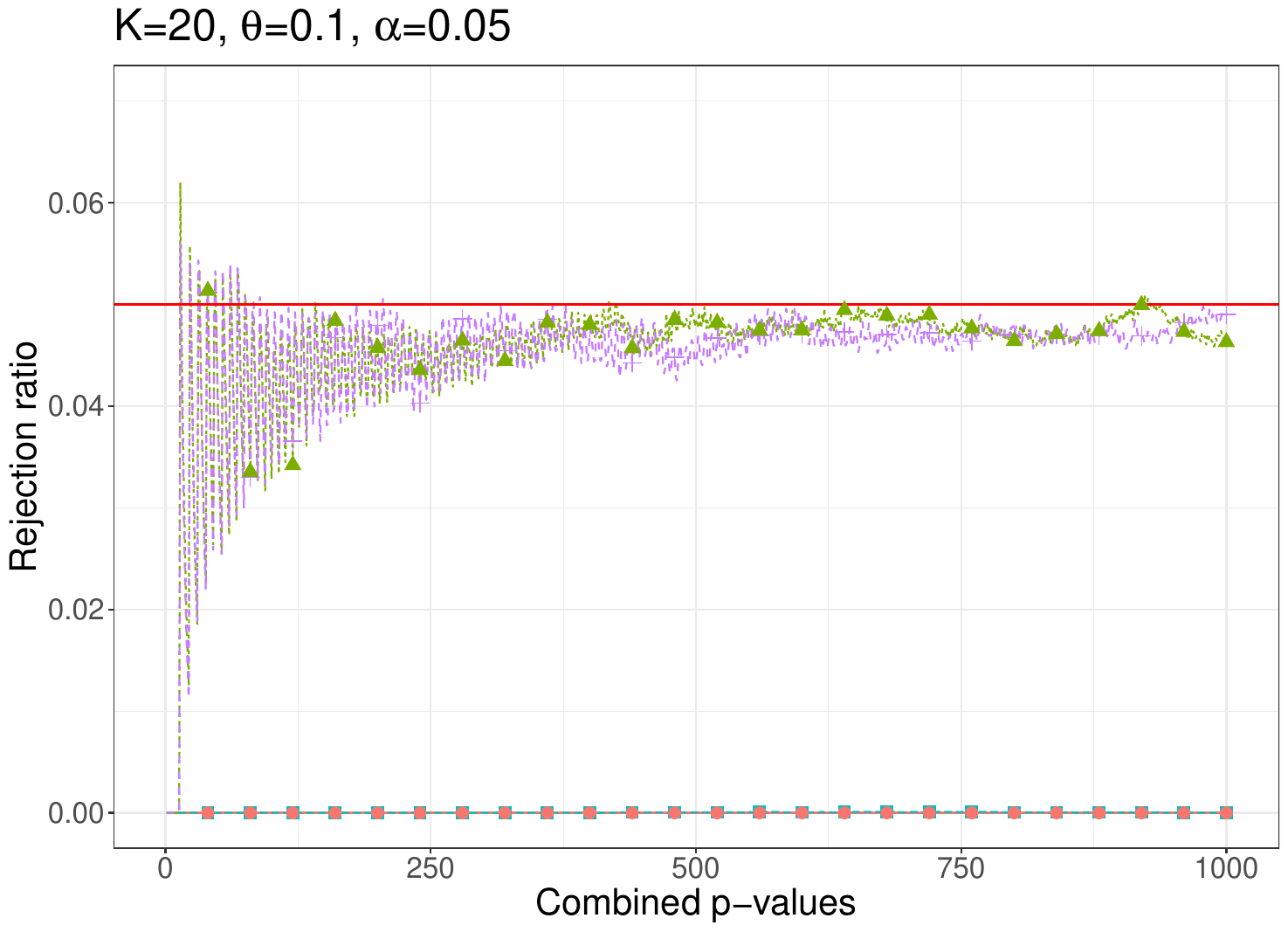}}\\
			{\includegraphics[width=0.32\linewidth, height=3.5cm]{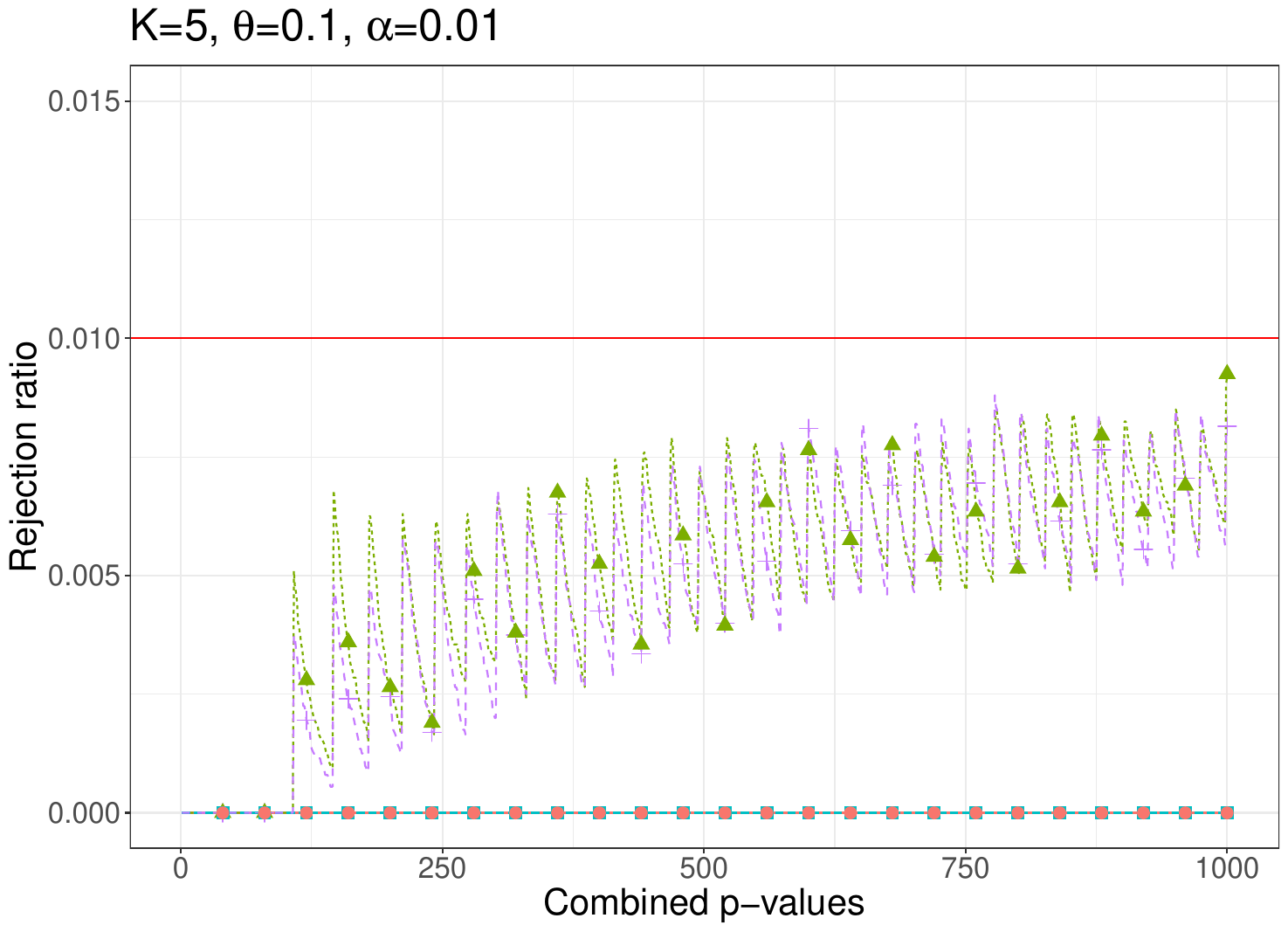}}
			{\includegraphics[width=0.32\linewidth, height=3.5cm]{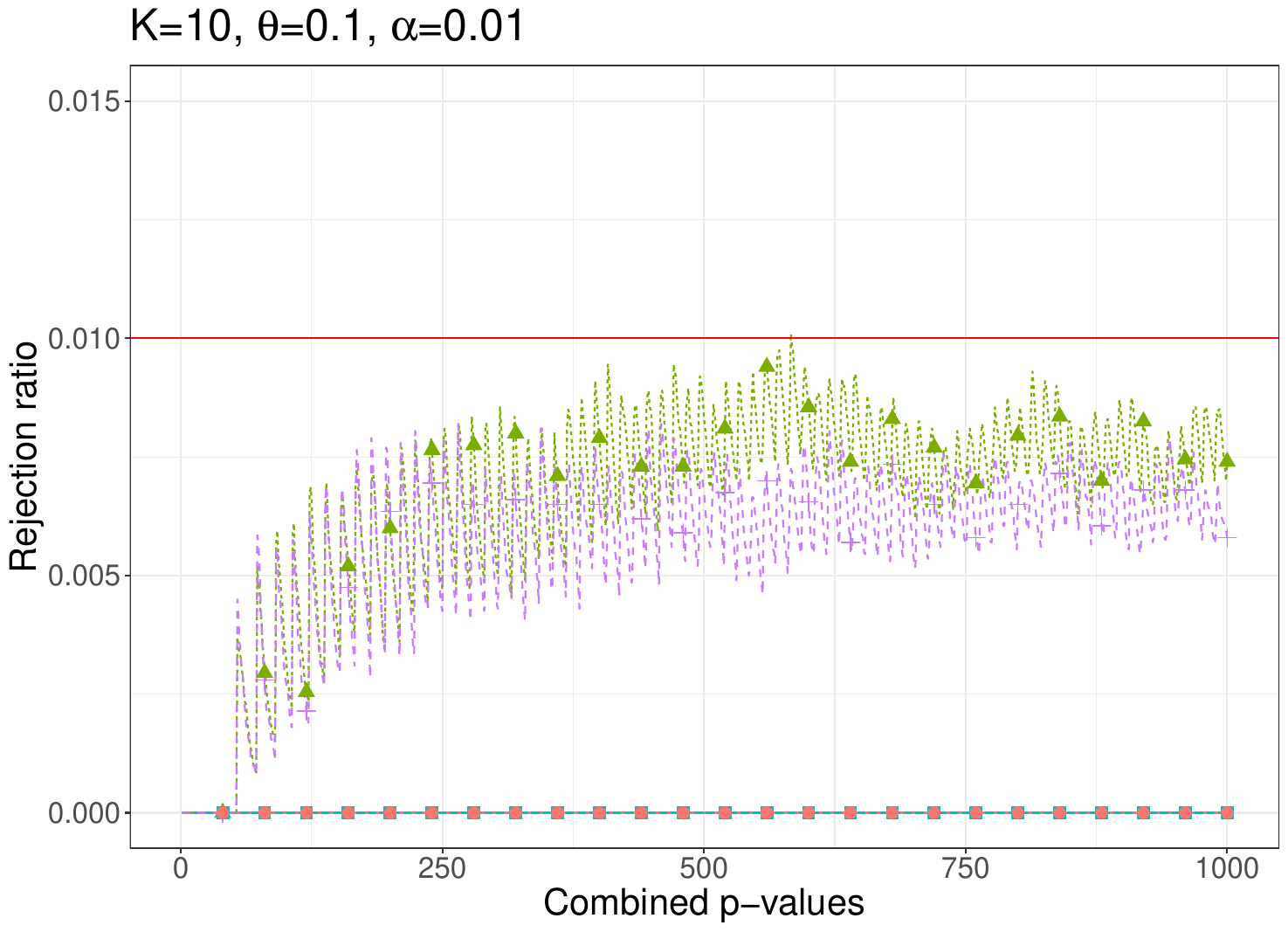}}
			{\includegraphics[width=0.32\linewidth, height=3.5cm]{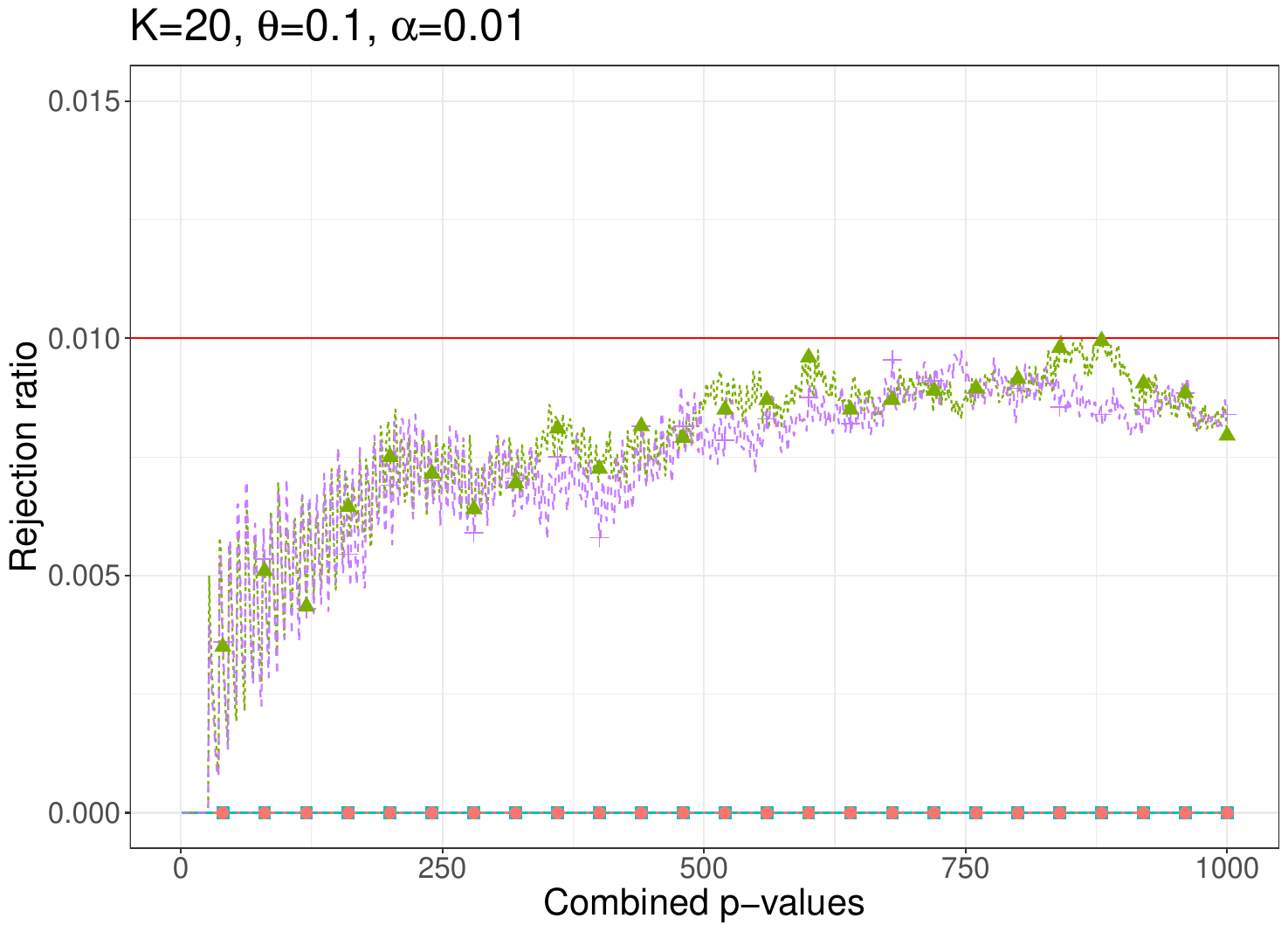}}
			\caption{Empirical type I error rates at various $\theta$, $K$ and $\alpha$ values. 
				Circle dots: $\tilde{S}_n$ with $\chi^2_{2n}$; 
				Triangles: $\tilde{S}_n$ with optimal gamma;
				Squares: $S_n$ with $\chi^2_{2n}$; 
				Crosses: $S_n$ with optimal gamma. 
			}
			\label{fig:pvcontrol05}
		\end{centering}
	\end{figure}

\begin{figure}[h!]
		\begin{tabular}{p{0.42\linewidth} p{0.12\linewidth} p{0.42\linewidth}}
			\vspace{0pt} {\includegraphics[width=\linewidth]{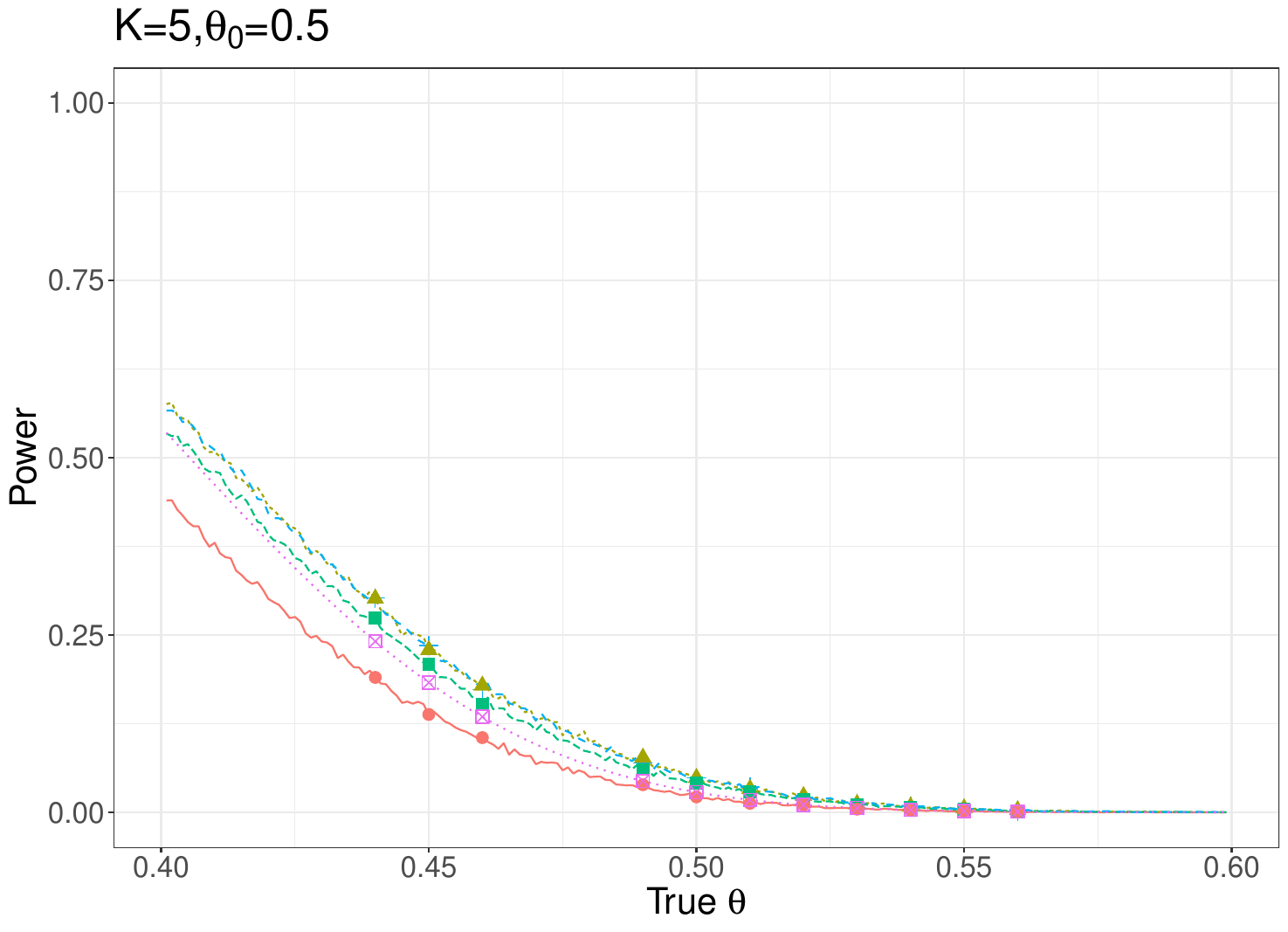}} &
			\vspace{20pt}
			{\includegraphics[width=\linewidth]{legendside.pdf}} &
			\vspace{0pt}
			{\includegraphics[width=\linewidth]{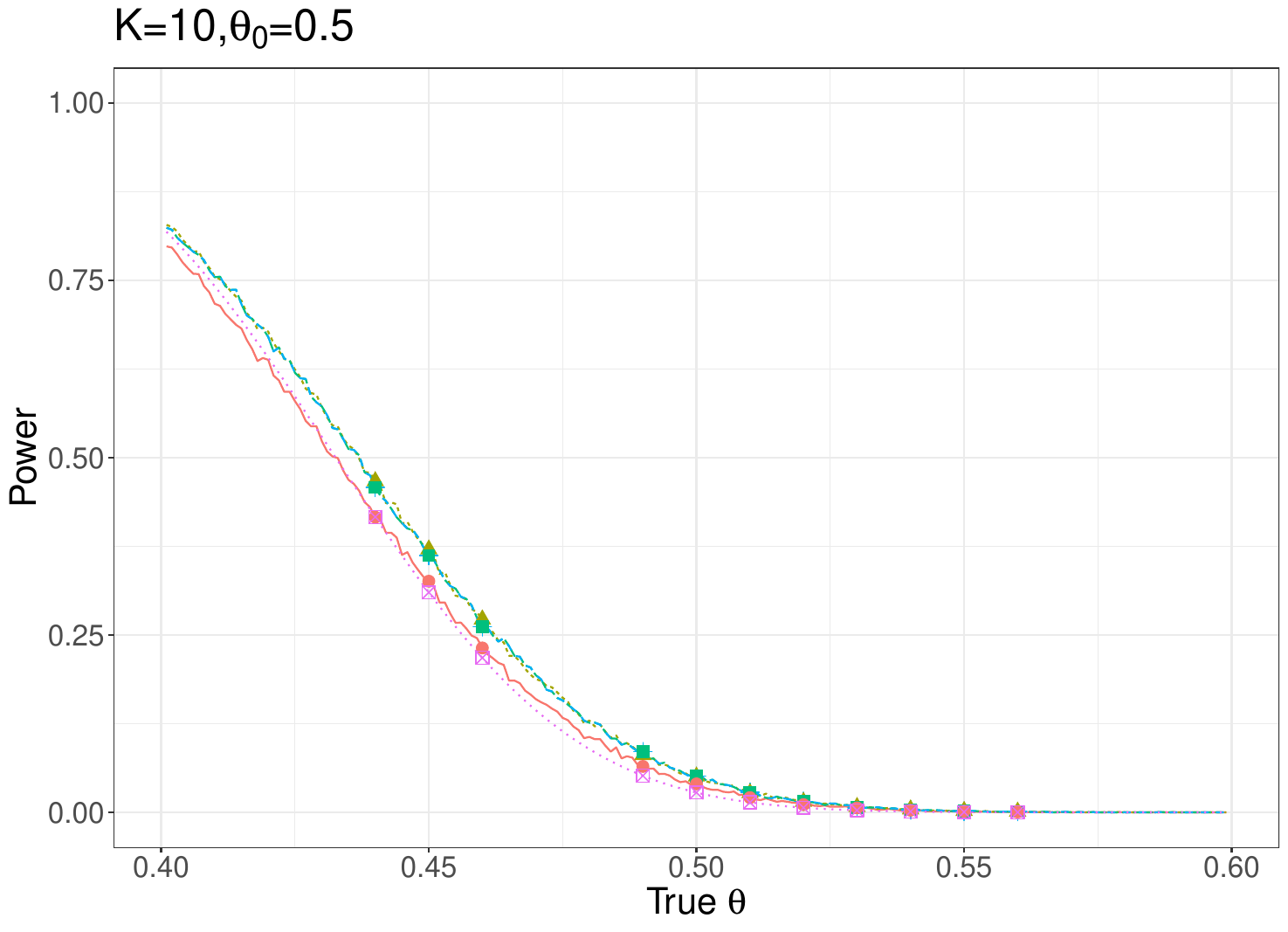}}  
		\end{tabular}
		\caption{Power curves over various actual $\theta$ values (x-axis) at size $\alpha=0.05$. 
			The null parameter is $\theta_0=0.5$ and $K=5$ (left) or $10$ (right). 
			Circle dots: $\tilde{S}_n$ with $\chi^2_{2n}$; 
			Triangles: $\tilde{S}_n$ with optimal gamma;
			Solid Squares: $S_n$ with $\chi^2_{2n}$; 
			Crosses: $S_n$ with optimal gamma;
			Hollow Squares: exact test in \eqref{eq:pvalueExactBinom}.
			The number of combined tests is $n=40$; the number of simulations is $N=10^4$. 
		} 
		\label{fig:power_theta0_0.5}
	\end{figure}

	The numerical and graphical results presented here and their descriptions arise from data simulated as described in section \ref{sec:numerial}  at various binomial parameters. Figure \ref{fig:histgof} presents simulated $S_n$ and $\tilde{S}_n$ compared with gamma and chi-squared approximations; Figure \ref{fig:pvcontrol05} shows the empirical type I error rates for varying nominal $\alpha$; and Figure \ref{fig:power_theta0_0.5} presents power comparisons for the approximations.
	
	\subsection{Combining independent, non-identically distributed p-values}
	
	We simulate the empirical distributions of $S_n$ and $\tilde{S}_n$ with $n=90$, $n=225$ and $n=900$ using $10^7$ simulated p-values obtained with equal frequency among distributions given in Table \ref{table:paramfit}. Optimal gamma null distributions are $S_n\sim Gamma(1.78n,1.12)$ and $\tilde{S}_n\sim Gamma(1.597n,1.08)$ based on mean (null) parameters $\bar{\nu}=2.24$, $\bar{m}=1.72$ and $\bar{v}=1.86$ (rounded to three decimal places).
	
	\begin{figure}[h!]
		\begin{centering}
			\resizebox{1\linewidth}{!}{\includegraphics{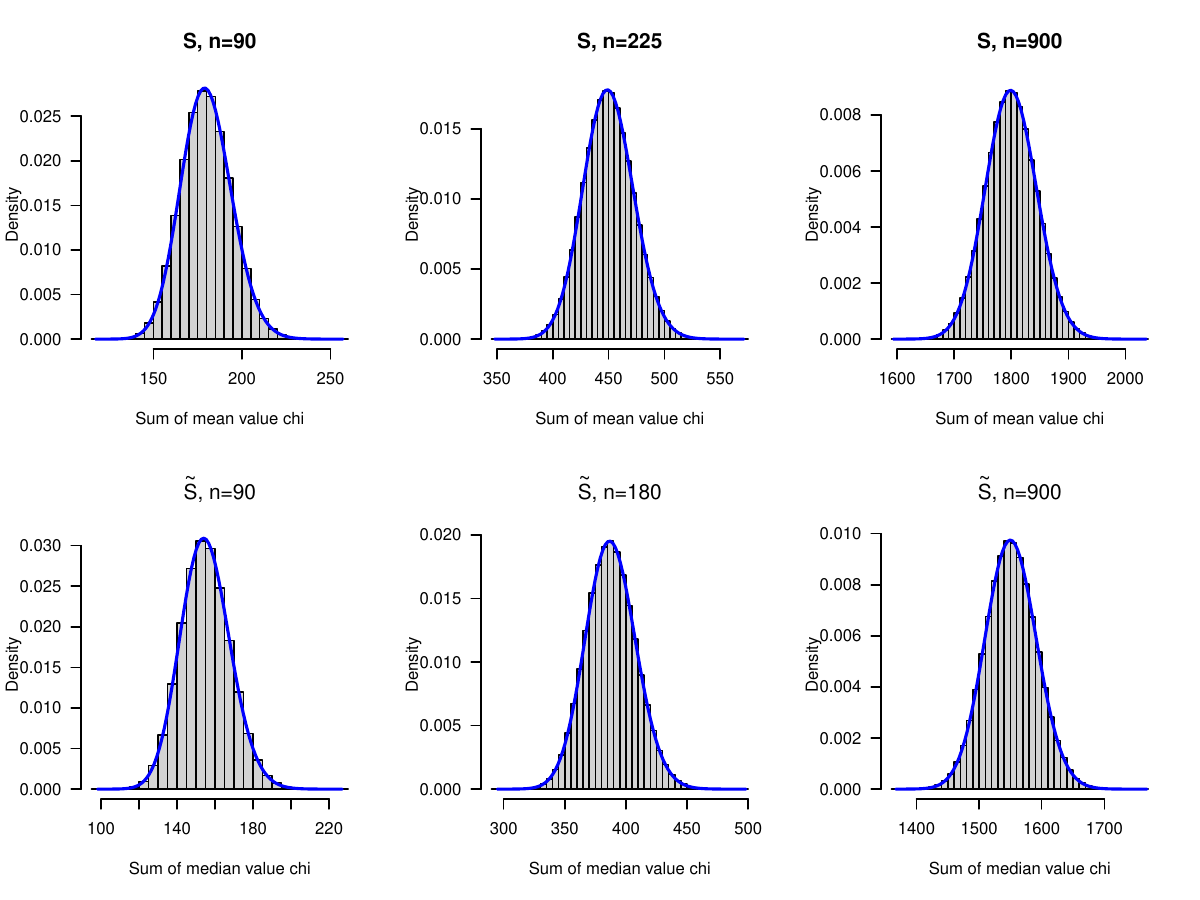}}
			\caption{Density histograms of simulated data and adjusted Gamma for non i.i.d. combinations.
				$10^7$ simulated $S_n$ and $\tilde{S}_n$ with $n=90$, $n=225$ and $n=900$. Null parameters chosen with equal frequency among distributions given in Table \ref{table:paramfit}. }
			\label{fig:convg}
			
		\end{centering}
	\end{figure}
	
	\begin{figure}[h!]
		\begin{centering}
			\resizebox{\linewidth}{!}{\includegraphics{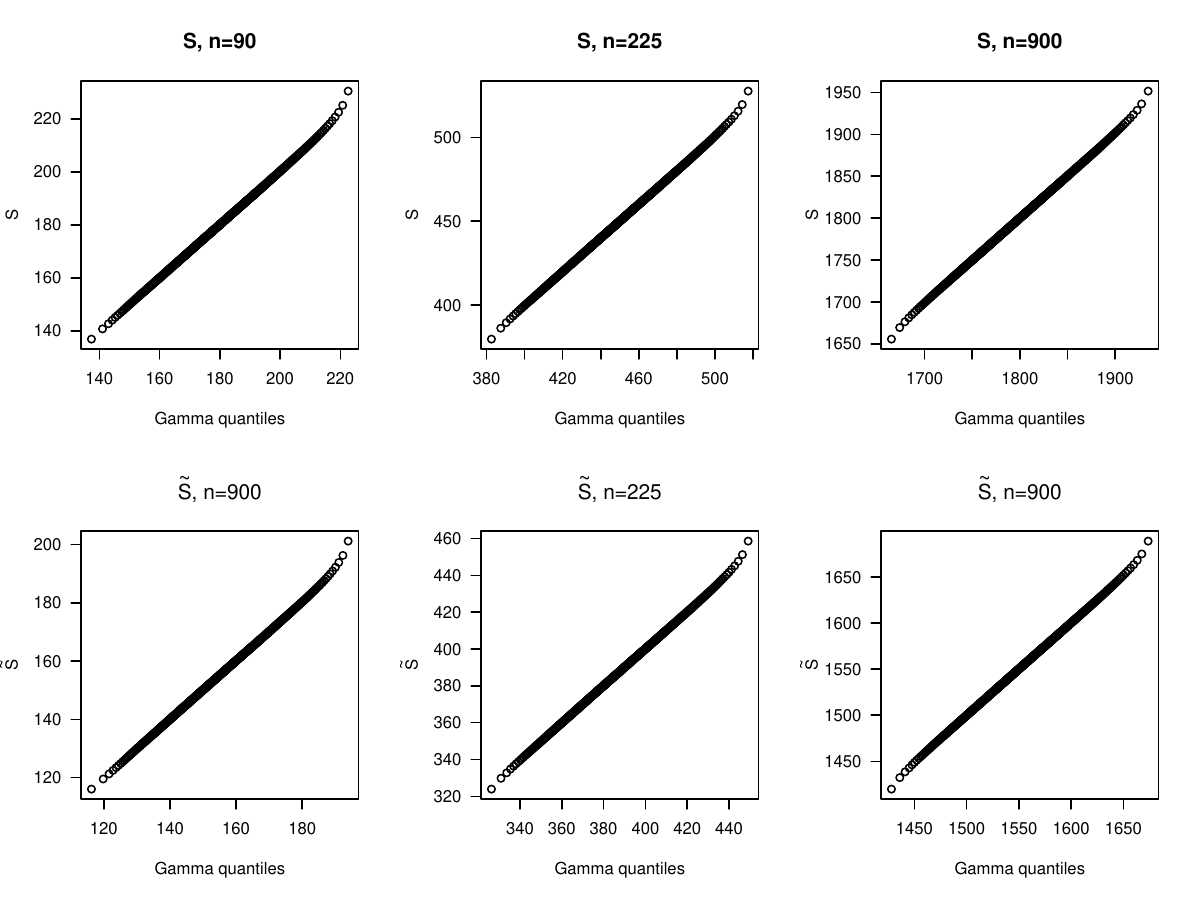}}
			\caption{Quantile plots of simulated data and adjusted Gamma for non i.i.d. combinations. Null parameters chosen with equal frequency from distributions given in Table \ref{table:paramfit}.}
			\label{fig:quantilesconv}
			
		\end{centering}
	\end{figure}

	{\bf Power analysis for non i.i.d. combinations.}
	
	We obtain binomial samples with the combination of parameters $(K,\theta_0)$ chosen with equal frequency among 9 distributions in Table \ref{table:paramfit} and $\kappa \in (0.01, 1.99)$ (so that $0\leq \kappa \theta_0\leq 1$ and $\kappa <1$ provides evidence favoring $\theta\leq \theta_0$ for all $\theta_0$ values in Table \ref{table:paramfit}). 
	The rejection rules for the global null are  those in \eqref{eq:rejection_chisq_gamma}, using the parameters $\bar{\nu}$, $\bar{m}$, and $\bar{v}$ to calculate the gamma quantiles.
	In this particular setting, for combinations of $n=40$ tests, the rejection rules are $J_n = I(S_n \geq 101.88)$, $\tilde{J}_n = I(\tilde{S}_n \geq 101.88)$, $I_n = I( S_n \geq 99.86)$ and $\tilde{I}_n= I(\tilde{S}_n \geq 93.92)$.

	Estimated power curves can be seen in Figure \ref{fig:powerniid}. We observe that, once again, Gamma based tests are uniformly more powerful than their chi-squared counterparts.

	\section{Discussions}\label{sec:sfig}

	\subsection{Discrepancy due to discreteness}\label{sec:discrepancy}
	
	
	Although appropriate continuous distributions can properly approximate adjusted discrete statistics, there is always a discrepancy between them due to the discreteness of the latter. In this section, we discuss this discrepancy through a few lower bounds related to $\Delta = \sup_{i \in \mathbb{N}} \{F_i - F_{i-1}\}$, which is the largest probability mass of the discrete statistic. These lower bounds provide an ``at least" discrepancy with interpretable formulas based on the properties of the discrete and continuous distributions. They help quantify the quality and limitations of the approximation procedure and explain how the proposed gamma approximation in Theorem \ref{thm:optgamma} reduces the overall discrepancy compared to that of $\chi^2_2$ for approximating the distributions of $\ddot{Z}$ and $\tilde{Z}$, the mean-value-$\chi^2$ and median-value-$\chi^2$ in \eqref{eq:lmod_lexplicit} and \eqref{eq:medianchi}, respectively.

	
	\begin{proposition}
		\label{cor:lowerbounds}
		Let $X\in \{x_i; i \in \mathbb{N}\}$ have cumulative distribution function $F_i=\P(X \leq X_i)$ with $0=F_0$, and $Y\sim G$ with $G$ being a strictly increasing cumulative distribution function. 
		It holds that: 
		\begin{align}
			\label{eq:lowerb_general}
			W_2^2(X,Y)  \geq \sup_{i \in \mathbb{N}}\int_{G^{-1}(F_{i-1})}^{G^{-1}(F_i)}(x_i-y)^2g(y)dy =\sup_{i \in \mathbb{N}} \E((X-Y)^2|A_i)(F_i-F_{i-1}). 
		\end{align}
		
		Particularly, if $Y\sim \chi_2^2$ and $\ddot{Z}$ is the mean-value-$\chi^2$, we have
		\begin{equation}
			\label{eq:lowerb_chisq}
			W_2^2(\ddot{Z},Y)\geq 4\max_{i \in \mathbb{N}} \{ F_i-F_{i-1}-\frac{F_iF_{i-1}(\log F_i -\log F_{i-1})^2}{F_i-F_{i-1}} \}.
		\end{equation}
		If $Y\sim \chi_2^2$ and $\tilde{Z}$ is the median-value-$\chi^2$, we have
		\begin{equation}
			\label{eq:lowerb_chisq_midp}
			W_2^2(\tilde{Z},Y)\geq 4\max_{i \in \mathbb{N}} \{ (F_i-F_{i-1})(1+c^2_i)-\frac{F_iF_{i-1}(\log F_i -\log F_{i-1})^2}{F_i-F_{i-1}} \},
		\end{equation}
		where 
		$$
		c_i=1-\frac{F_i\log F_i -F_{i-1}\log F_{i-1}}{F_i-F_{i-1}}+ \log \left(\frac{F_i+F_{i-1}}{2}\right).
		$$
	\end{proposition}

	These lower bounds account for the discrepancies due to the largest probability point mass. Equation \eqref{eq:lowerb_general} provides a general lower bound for the discrepancy between an arbitrary discrete-continuous pair $(X, Y)$ based on the point masses $F_i-F_{i-1}$ of $X$ and the distribution of $Y$. As special cases, when $Y\sim \chi_2^2$, \eqref{eq:lowerb_chisq} and \eqref{eq:lowerb_chisq_midp} correspond to $X=\ddot{Z}$ and $X=\tilde{Z}$, respectively. These lower bounds all shrink to zero as $\Delta \to 0$, which is consistent with Corollary \ref{cor:closetc}.
	
	The lower bounds in \eqref{eq:lowerb_chisq} and \eqref{eq:lowerb_chisq_midp} differ by a term $c^2_i$, which is equal to $(\ddot{Z} -\tilde{Z})^2$ on the event $P=F_i$. 
	Therefore, the lower bound of the discrepancy between $\tilde{Z}$ and $\chi^2_2$ is always larger than that between $\ddot{Z}$ and $\chi^2_2$, which is consistent with Corollary \ref{cor:closetochi1}. Meanwhile, $c_i \to 0$ as $\Delta \to 0$ (because, for some $F_{i-1} \leq \xi_i \leq F_i$, $|c_i|=\left|\log\{(F_i+F_{i-1})/(2\xi_i)\}\right| \leq \left|(F_i+F_{i-1})/(2\xi_i)-1\right| $).  \cite{lancaster1949combination} argued that in this close-to-continuous scenario, $\tilde{Z}$ can be used as an approximation to $\ddot{Z}$. Indeed, since $W_2^2(\tilde{Z},\chi_2^2)=W_2^2(\ddot{Z},\chi_2^2)+\E[(\ddot{Z} -\tilde{Z})^2]$, both statistics are close to each other, and their distributions are close to the target $\chi_2^2$ distribution when $\Delta$ is small.
	
	
	These lower bounds suggest that the discrepancy between a discrete-continuous pair can be substantial when there is a single large probability mass. In such cases, the lower bounds may provide a reasonable estimate of the entire discrepancy. We illustrate this point with an example: consider a discrete statistic $X$ with cumulative distribution function $F_i=i/(10K)$ for $i = 0,1,\ldots ,K$ and $F_{K+1}=1$, where $K \in \mathbb{N}$. The corresponding left-sided p-value is $P=F_i$ with probability $p_i = F_{i}-F_{i-1} = 0.1/K$ for $i \leq K$ and $P=1$ with a relatively large probability $F_{K+1}-F_{K}= 0.9$. For any $K$ value, the lower bound of $W_2^2(\ddot{Z},Y)$ in \eqref{eq:lowerb_chisq} is $4(0.9-0.1\log^2( 0.1)/0.9)\approx 1.2436$. This lower bound value is close to the overall discrepancy at large $K$, e.g., $W_2^2(\ddot{Z},Y) = 4-\nu\approx1.2479$ by Corollary \ref{cor:wdmoments2} when $K=100$. Similarly, the lower bound of $W_2^2(\tilde{Z},Y)$ in \eqref{eq:lowerb_chisq_midp} is $4(0.9(1+(1+\log 0.1/9 + \log 0.55)^2)-0.1\log^2( 0.1)/0.9)\approx 1.32$, which is very close to $W_2^2(\tilde{Z},Y)=4-\nu+\E([\ddot{Z}-\tilde{Z}]^2)\approx 1.3223$, suggesting that the majority of the discrepancy between $\tilde{Z}$ and $\chi^2_2$ is due to the distribution discreteness characterized by the largest probability point mass.

	\begin{figure}[h!]
		\begin{centering}
			\subfloat[$\ddot{Z}$ vs. $\chi_2^2$]{\includegraphics[width=2.2in, height=1.8in]{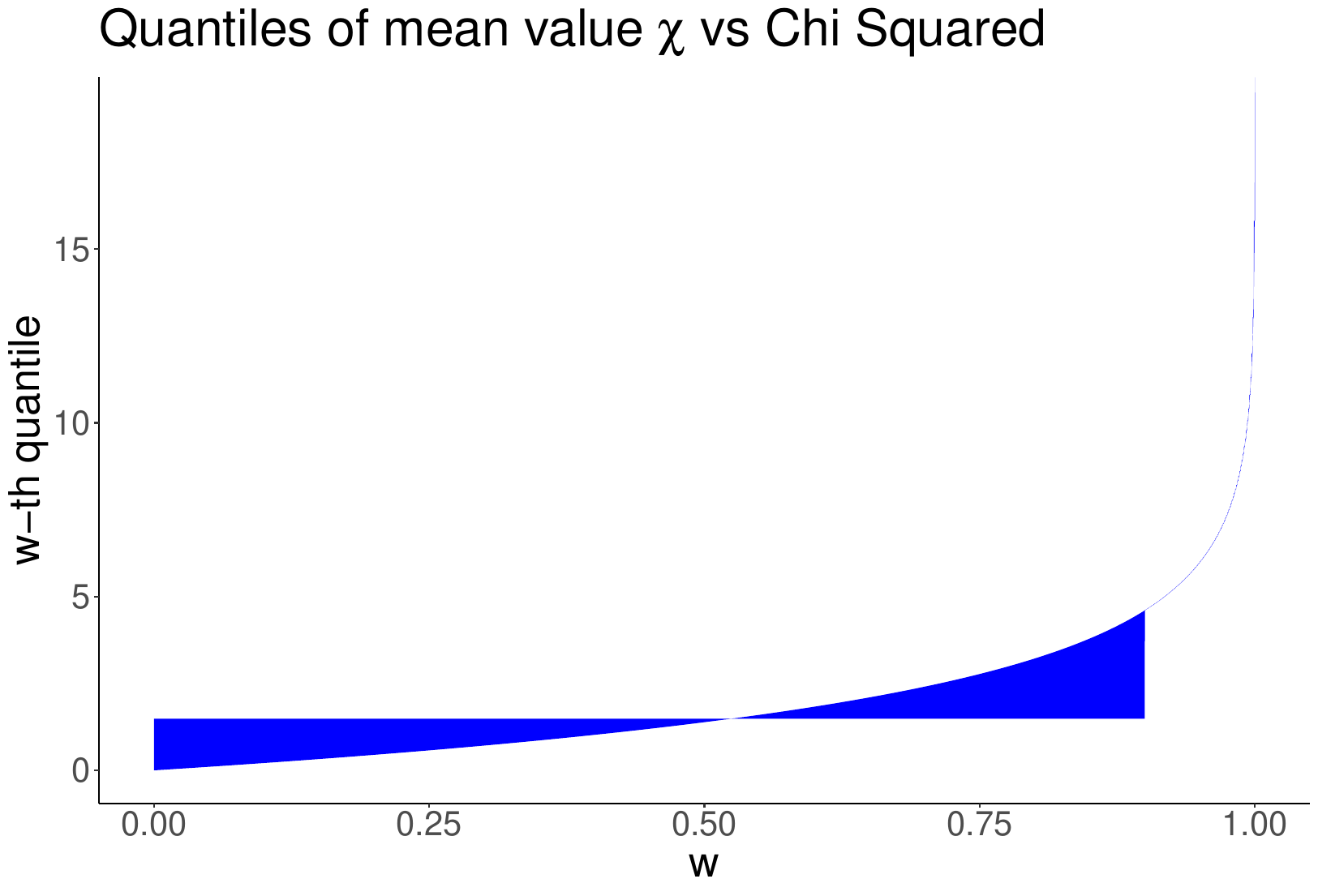}}
			\subfloat[$\ddot{Z}$ vs. gamma]{\includegraphics[width=2.2in, height=1.8in]{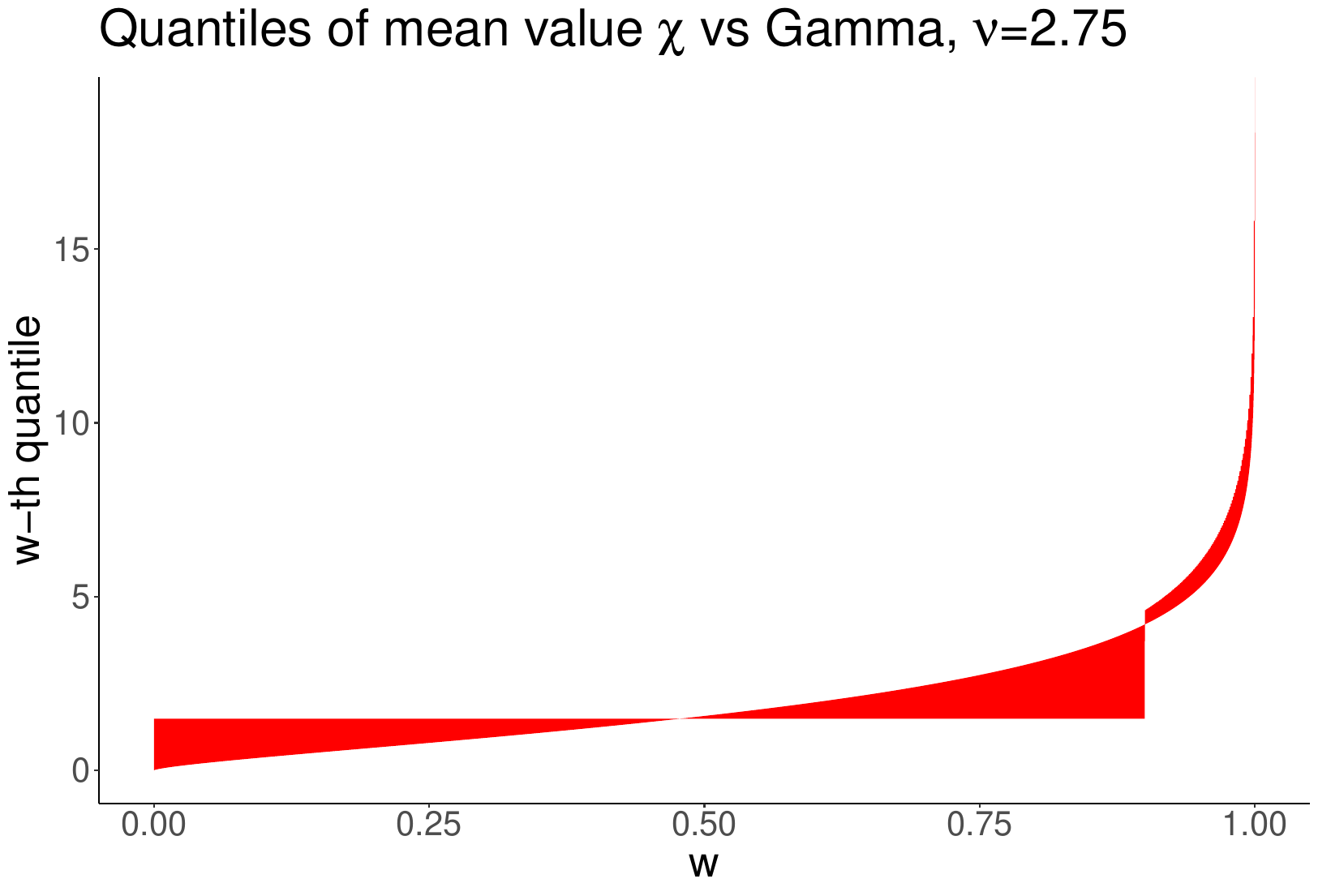}} \\
			\subfloat[$\tilde{Z}$ vs. $\chi_2^2$]{\includegraphics[width=2.2in, height=1.8in]{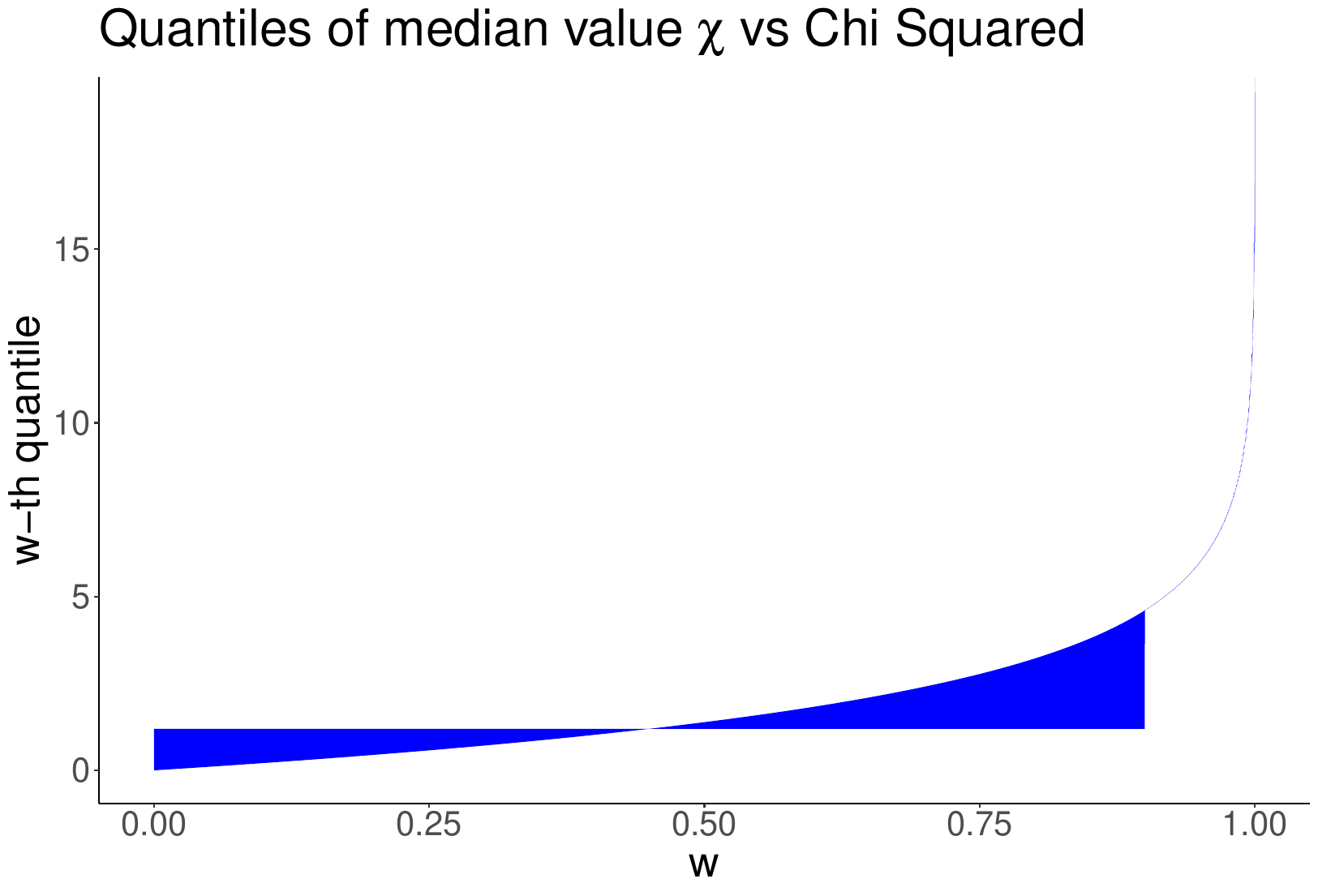}}
			\subfloat[$\tilde{Z}$ vs. gamma]{\includegraphics[width=2.2in, height=1.8in]{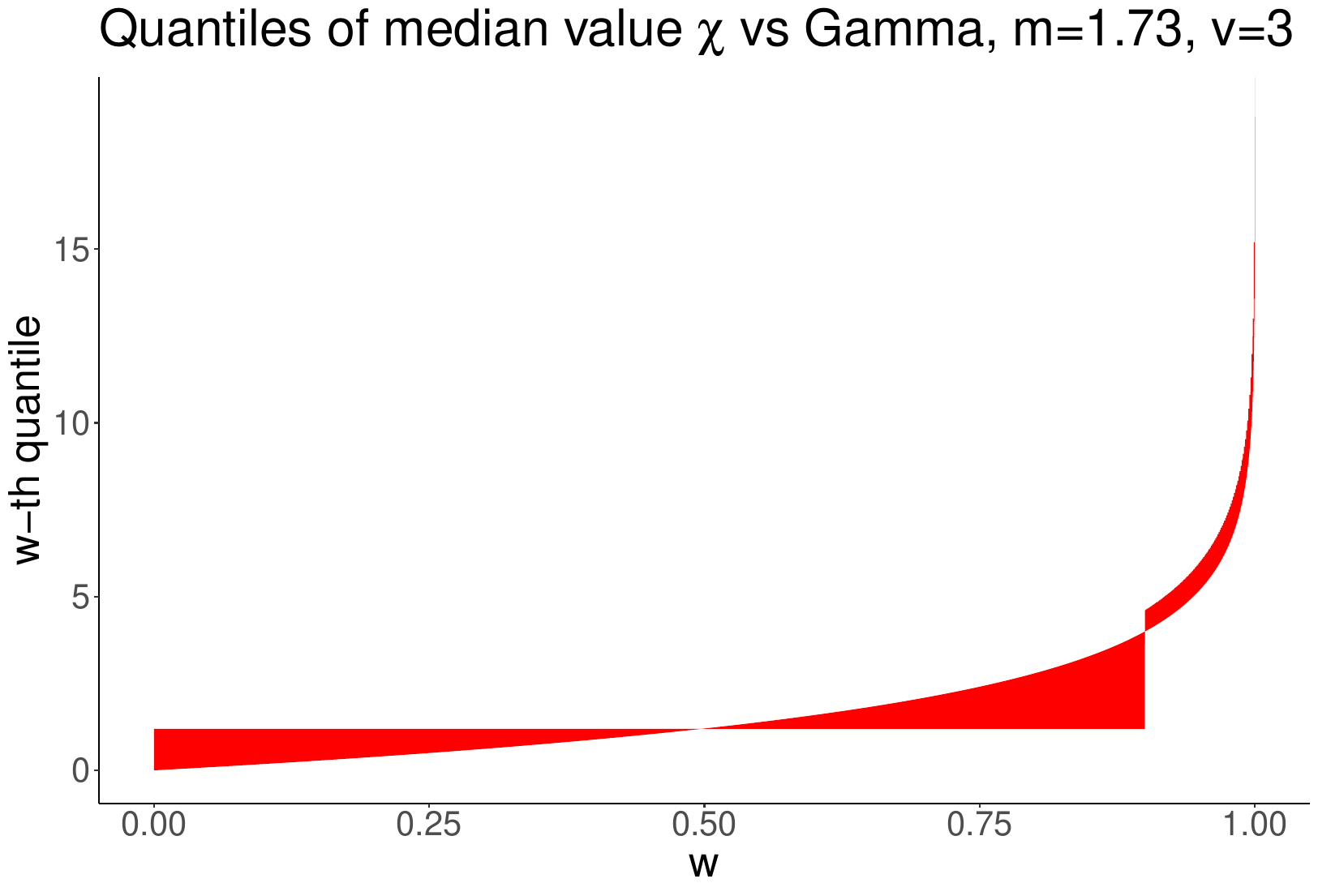}}
			\caption{Shaded area between a pair of discrete and continuous 
				distributions is the Wasserstein distance between them, as per \eqref{eq:wass} with $p=2$.
				Y-axis: distribution quantiles; X-axis: integration variable in \eqref{eq:wass}.
				The discrete statistics are the mean-value-$\chi^2_2$ statistic $\ddot{Z}$ in \eqref{eq:lmod_lexplicit} (upper row) and the median-value-$\chi^2$ statistic $\tilde{Z}$ in \eqref{eq:medianchi} (lower row). 
				The continuous distributions are $\chi_2^2$ (left column) and the optimal gamma distributions by Theorem \ref{thm:optgamma} (right column). 
				The dominant bulks of discrepancy, at the left of each figure, are the lower bounds given by Corollary \ref{cor:lowerbounds}. 
				Compared to $\chi^2_2$ (left column, blue), the gamma adjustments (right column, red) reduces overall discrepancy by substantially reducing its dominant bulk of discrepancy at the largest probability point mass, even if at a cost of slightly increasing the discrepancy at the closer-to-continuous tail of the distribution.
			} 
			\label{fig:quoaexlanc}
		\end{centering}
	\end{figure}

	The lower bounds in Proposition \ref{cor:lowerbounds} reveal one mechanism of how the optimal gamma approximations reduce the discrepancies when compared with the $\chi^2_2$ distribution. Specifically, the gamma distributions decrease the overall discrepancy between the discrete and continuous distributions by significantly reducing the influence of the large probability mass. To see this point, we can compare the ratios between the total discrepancy and the lower bound given in \eqref{eq:lowerb_general}. In the above example, with $X=\ddot{Z}$ and $Y\sim\chi^2_2$, the ratio is $1.2454/1.2479 \approx 0.998$. Changing $Y$ to the optimal gamma distribution reduces the ratio to $1.0694/1.111\approx 0.9626$. That is, the gamma distribution decreases the overall discrepancy from 1.2479 to 1.111 by significantly reducing the contribution of the largest probability mass (calculated by the lower bound) from $1.2436$ to $1.0694$ (even if at the cost of slightly increasing the discrepancy at the other part of the distribution from $1.2479-1.2436 = 0.0043$ to $1.111-1.0694 = 0.0416$). Similar results hold for $\tilde{Z}$. 
	Figure \ref{fig:quoaexlanc} visualizes the overall discrepancies by the areas between the discrete and the continuous distributions, where the largest bulk is contributed by the largest point probability mass. The ratios between the lower bounds and the total discrepancies are  (a) $1.2436/1.247 \approx 0.9975$; (b) $1.0694/1.111\approx 0.9626$; (c) $1.32/1.323\approx 0.997$; (d) $0.944/1.03\approx 0.916$.
	
	\bigskip
	{\bf Proof of Proposition \ref{cor:lowerbounds}}
	
	\begin{proof}
		Since for $w \in [F_{i-1},F_i)$, $F^{-1}(w)=x_i$ one has: 
		\begin{align*}
			W_2^2(X,Y) &= \int_{0}^1 |F^{-1}(w)-G^{-1}(w)|dw\\
			= &\sum_{i \in \mathbb{N}}\int_{G^{-1}(F_{i-1})}^{G^{-1}(F_i)}(x_i-y)^2g(y)dy\\  \geq & \sup_{i \in \mathbb{N}} \int_{G^{-1}(F_{i-1})}^{G^{-1}(F_i)}(x_i-y)^2g(y)dy\\
			= & \sup_{i \in \mathbb{N}} \E((X-Y)^2|A_i)(F_i-F_{i-1}).\\
		\end{align*}
		To get the last equation, note that $X=x_i$ almost surely when constrained to the sets $A_i$ defined in \eqref{eq:partition}, and therefore $\E((X-Y)^2|A_i)= \int_{G^{-1}(F_{i-1})}^{G^{-1}(F_i)}(x_i-y)^2g(y)dy/\P(A_i)$.  
		
		When $Y\sim \chi_2^2$ and $\ddot{Z}$ is the mean-value-$\chi^2$, we can further obtain
		\begin{align*}
			W_2^2(\ddot{Z},Y)&=4-\nu \\
			&=4\sum_{i \in \mathbb{N}} \left[ (F_i-F_{i-1})-\frac{(F_i\log F_i -F_{i-1}\log F_{i-1})^2}{F_i-F_{i-1}}\right]\\
			&\geq 4\max_{i \in \mathbb{N}} \left[ (F_i-F_{i-1})-\frac{(F_i\log F_i -F_{i-1}\log F_{i-1})^2}{F_i-F_{i-1}}\right]\\
			&\geq 4\max_{i \in \mathbb{N}} \left[ F_i-F_{i-1}-\frac{F_iF_{i-1}(\log F_i -\log F_{i-1})^2}{F_i-F_{i-1}} \right].
		\end{align*}
		
		For the case of $Y\sim \chi_2^2$ and $\tilde{Z}$ is the median-value-$\chi^2$, we apply the equality (see the proof of Corollary \ref{cor:closetochi1}):
		\begin{equation}
			\label{eq:midpdistance}
			W_2^2(\tilde{Z},\chi_2^2)=W_2^2(\ddot{Z},\chi_2^2)+\Var(\ddot{Z} -\tilde{Z}) +(m-2)^2=W_2^2(\ddot{Z},\chi_2^2)+\E[(\ddot{Z} -\tilde{Z})^2].
		\end{equation}
		
		Denote $c_i=1-(F_i\log F_i -F_{i-1}\log F_{i-1})/(F_i-F_{i-1})+\log [(F_i+F_{i-1})/2]$. Using \eqref{eq:midpdistance}, we have that for $\tilde{Z}$, 
		
		\begin{align*}
			W_2^2(\tilde{Z},Y) 
			&= 4\sum_{i \in \mathbb{N}}(F_i-F_{i-1}) \biggr[1 -\frac{(F_i\log F_i -F_{i-1}\log F_{i-1})^2}{(F_i-F_{i-1})^2} \\
			&+(1-\frac{(F_i\log F_i -F_{i-1}\log F_{i-1})}{(F_i-F_{i-1})}+\log(F_i/2+F_{i-1}/2))^2\biggr]\\ 
			&= 4\sum_{i \in \mathbb{N}} \left[ (F_i-F_{i-1})(1+c^2_i)-\frac{(F_i\log F_i -F_{i-1}\log F_{i-1})^2}{F_i-F_{i-1}}  \right]\\
			&\geq 4\max_{i \in \mathbb{N}} \left[ (F_i-F_{i-1})(1+c^2_i)-\frac{(F_i\log F_i -F_{i-1}\log F_{i-1})^2}{F_i-F_{i-1}}  \right]\\
			&\geq 4\max_{i \in \mathbb{N}} \left[ (F_i-F_{i-1})(1+c^2_i)-\frac{F_iF_{i-1}(\log F_i -\log F_{i-1})^2}{F_i-F_{i-1}} \right].
		\end{align*}
	\end{proof}

	\subsection{Limitation of the original Fisher's statistic in combining discrete p-values}
	
	
	In this section, we briefly discuss a limitation of the original Fisher's statistic, $T_n=\sum_{j=1}^n -2\log(P_j)$ in \eqref{eq:fisherchi}, when it comes to combining discrete p-values, especially when $P_j=1$ has large probability. 
	As illustrated by an example, although a gamma distribution could be used to approximate the distribution of $T_n$, such a testing procedure is less favored compared to those using the adjusted combination statistics $S_n$ and $\tilde{S}_n$ defined in \eqref{eq:sum12}.

		Consider the combination of $n=40$ independent, identically distributed experiments with discrete p-values $P=F_i$ with probability $p_i = F_{i}-F_{i-1}=0.001$, where $F_i=i/(1000)$ for $i = 0,1,\ldots ,100$ and $F_{101}=1$ with probability $p_{101} = F_{101}-F_{100}=0.9$. When $P=1$, $-2\log(P)=0$, $\ddot{Z}=1.49$ and $\tilde{Z}=1.2$, where $\ddot{Z}$ and $\tilde{Z}$ denote the mean-value-$\chi^2$ and median-value-$\chi^2$ as expressed in \eqref{eq:lmod_lexplicit} and \eqref{eq:medianchi}, respectively. Suppose that $P_i=1$ for all 40 observed significances, an event that occurs with probability $0.9^{40}\approx 0.0148$. In this scenario, $T_n=0$, hence, any Gamma distribution would yield $\P (T_n\leq 0)=0$. In contrast, the adjusted combination statistics $S_n=59.5326$ and $\tilde{S}_n=47.827$. The quantiles of their respective gamma approximations ($\nu=2.7521$, $m=1.7359$, $v=3.01$) are $\P(S_n\leq 59.5326)= 0.0177$ and $\P(\tilde{S}_n\leq47.827)=0.0151$, which approximate $0.0148$, the null probability of all p-values being 1, more accurately. 
		
\begin{figure}
			\begin{centering}
				\resizebox{\linewidth}{!}{\includegraphics{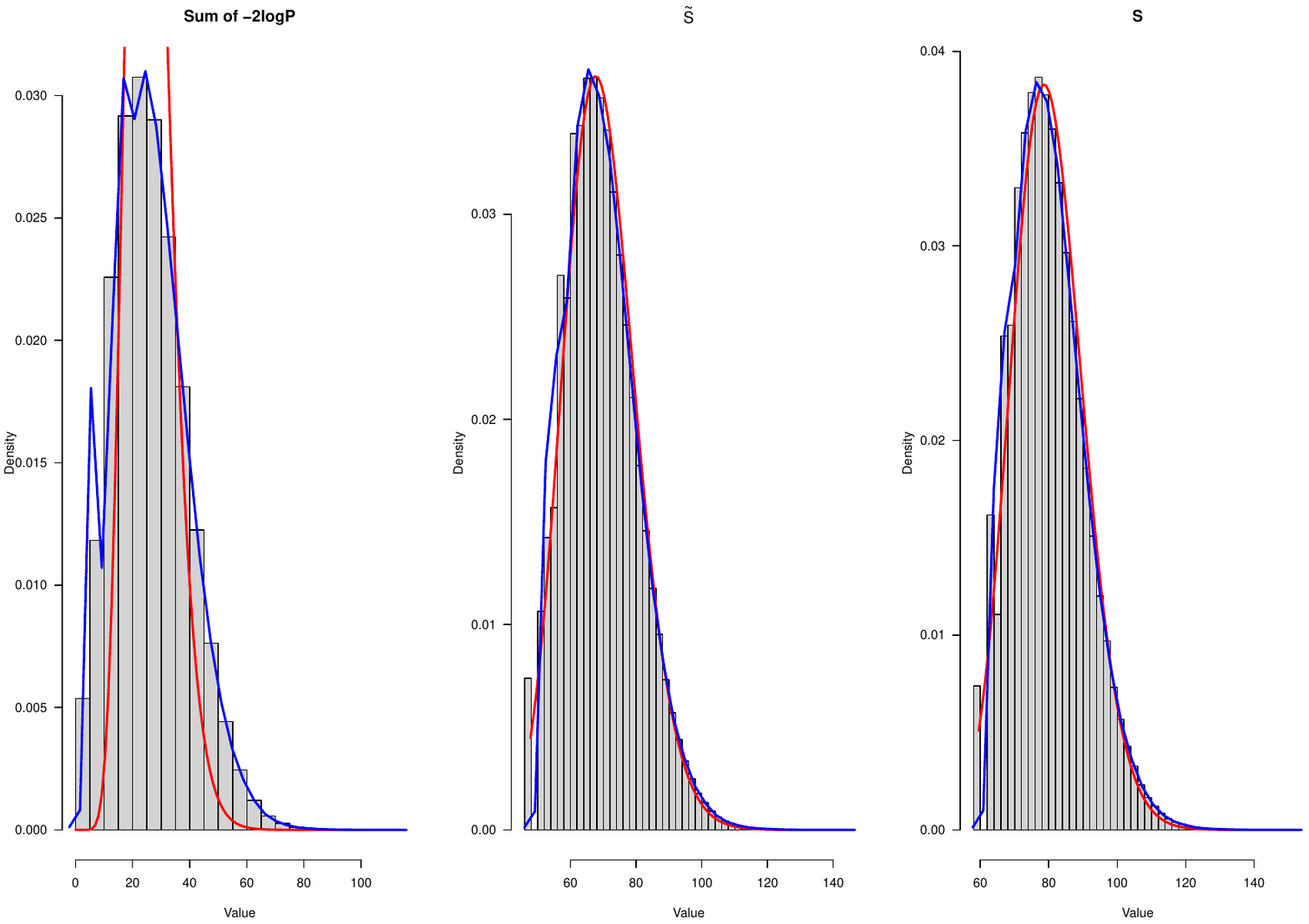}}
				\caption{Density histograms of simulated data in contrast with the optimal gamma density curves in red. The simulations are based on $10^7$ runs of $T_n$ (left), $\tilde{S}_n$ (middle), and $S_n$ (right) combining $n=40$ discrete p-values. The Gaussian kernel density estimator  (using 32 kernels, in blue) is drawn to effectively capture the profile of the histograms. Overall, the optimal gamma distribution for $T_n$ does not align well with its histogram. In contrast, the optimal gamma distributions of $S_n$ and $\tilde{S}_n$ fit well with their respective histograms. 
				}
				\label{fig:notraw}
				
			\end{centering}
		\end{figure}

		
		Figure \ref{fig:notraw} illustrates this issue with $N=10^7$ simulations of these random variables. The distribution of $T_n$ has a few modal peaks close to 0, which any gamma distribution would fail to model. Furthermore, the distribution of $T_n$ exhibits a heavy right tail, which the gamma does not accommodate due to its left-skewed mean. In contrast, the empirical distributions of $S_n$ and $\tilde{S}_n$ are more accurately approximated by the Gamma distribution. 
		

\end{appendix}


	\bibliographystyle{imsart-number} 
	\bibliography{CombDiscreteTests}      

\begin{thebibliography}{}

\bibitem[\protect\citeauthoryear{Barnard}{Barnard}{1989}]{barnard}
Barnard, G. (1989, December).
\newblock On alleged gains in power from lower p-values.
\newblock {\em Statistics in medicine\/}~{\em 8\/}(12), 1469---1477.

\bibitem[\protect\citeauthoryear{Berry and Armitage}{Berry and Armitage}{1995}]{Berry}
Berry, G. and P.~Armitage (1995).
\newblock Mid-p confidence intervals: a brief review.
\newblock {\em The Statistician\/}~{\em 44}, 417--423.

\bibitem[\protect\citeauthoryear{Billingsley}{Billingsley}{1986}]{Bill86}
Billingsley, P. (1986).
\newblock {\em Probability and Measure\/} (Second ed.).
\newblock John Wiley and Sons.

\bibitem[\protect\citeauthoryear{Ca{\~n}as and Rosasco}{Ca{\~n}as and Rosasco}{2012}]{Caas2012LearningPM}
Ca{\~n}as, G.~D. and L.~Rosasco (2012).
\newblock Learning probability measures with respect to optimal transport metrics.
\newblock In {\em Neural Information Processing Systems}.

\bibitem[\protect\citeauthoryear{Casella and Berger}{Casella and Berger}{2002}]{Casella2002Statistical-inf}
Casella, G. and R.~L. Berger (2002).
\newblock {\em Statistical inference\/} (2nd ed.), Volume~2.
\newblock Duxbury Pacific Grove, CA.

\bibitem[\protect\citeauthoryear{Clement and Desch}{Clement and Desch}{2008}]{Clement2008WD}
Clement, P. and W.~Desch (2008).
\newblock An elementary proof of the triangle inequality for the wasserstein metric.
\newblock {\em Proceedings of the American Mathematical Society\/}~{\em 136\/}(1), 333--339.

\bibitem[\protect\citeauthoryear{Csiszar}{Csiszar}{1975}]{Csizar1975metrics}
Csiszar, I. (1975).
\newblock $i$-divergence geometry of probability distributions and minimization problems.
\newblock {\em The Annals of Probability\/}~{\em 3\/}(1), 146--158.

\bibitem[\protect\citeauthoryear{Edgington}{Edgington}{1972}]{edgington1972additive}
Edgington, E.~S. (1972).
\newblock An additive method for combining probability values from independent experiments.
\newblock {\em The Journal of Psychology\/}~{\em 80\/}(2), 351--363.

\bibitem[\protect\citeauthoryear{Fisher}{Fisher}{1925}]{fisher1925statistical}
Fisher, R. (1925).
\newblock {\em Statistical methods for research workers}.
\newblock Edinburgh Oliver \& Boyd.

\bibitem[\protect\citeauthoryear{Fisher}{Fisher}{1922}]{fisher1922contingency}
Fisher, R.~A. (1922).
\newblock On the interpretation of chi-squared from contingency tables, and the calculation of p.
\newblock {\em Journal of the Royal Statistical Society\/}~{\em 85\/}(1), 87--94.

\bibitem[\protect\citeauthoryear{George and Mudholkar}{George and Mudholkar}{1983}]{george1983convolution}
George, E. and G.~Mudholkar (1983).
\newblock On the convolution of logistic random variables.
\newblock {\em Metrika\/}~{\em 30}, 1--13.

\bibitem[\protect\citeauthoryear{Kantorovich}{Kantorovich}{2006}]{Kantorovich2006}
Kantorovich, L.~V. (2006, Mar).
\newblock On the translocation of masses.
\newblock {\em Journal of Mathematical Sciences\/}~{\em 133\/}(4), 1381--1382.

\bibitem[\protect\citeauthoryear{Kincaid}{Kincaid}{1962}]{kinkaid}
Kincaid, W.~M. (1962).
\newblock The combination of tests based on discrete distributions.
\newblock {\em Journal of the American Statistical Association\/}~{\em 57\/}(297), 10--19.

\bibitem[\protect\citeauthoryear{Lancaster}{Lancaster}{1949}]{lancaster1949combination}
Lancaster, H. (1949).
\newblock The combination of probabilities arising from data in discrete distributions.
\newblock {\em Biometrika\/}~{\em 36\/}(3/4), 370--382.

\bibitem[\protect\citeauthoryear{Lin}{Lin}{2014}]{lin2014association}
Lin, W.-Y. (2014).
\newblock Association testing of clustered rare causal variants in case-control studies.
\newblock {\em PloS one\/}~{\em 9\/}(4), e94337.

\bibitem[\protect\citeauthoryear{Lin}{Lin}{2016}]{lin2016beyond}
Lin, W.-Y. (2016).
\newblock Beyond rare-variant association testing: pinpointing rare causal variants in case-control sequencing study.
\newblock {\em Scientific reports\/}~{\em 6\/}(1), 21824.

\bibitem[\protect\citeauthoryear{Mielke~Jr, Johnston, and Berry}{Mielke~Jr et~al.}{2004}]{mielke2004combining}
Mielke~Jr, P.~W., J.~E. Johnston, and K.~J. Berry (2004).
\newblock Combining probability values from independent permutation tests: a discrete analog of fisher's classical method.
\newblock {\em Psychological reports\/}~{\em 95\/}(2), 449--458.

\bibitem[\protect\citeauthoryear{Moschopoulos}{Moschopoulos}{1985}]{Mo1985}
Moschopoulos, P.~G. (1985).
\newblock The distribution of the sum of independent gamma random variables.
\newblock {\em Annals of the Institute of Statistical Mathematics\/}~{\em 37}, 541--544.

\bibitem[\protect\citeauthoryear{Neale, Rivas, Voight, Altshuler, Devlin, Orho-Melander, Kathiresan, Purcell, Roeder, and Daly}{Neale et~al.}{2011}]{Neale2011}
Neale, B.~M., M.~A. Rivas, B.~F. Voight, D.~Altshuler, B.~Devlin, M.~Orho-Melander, S.~Kathiresan, S.~M. Purcell, K.~Roeder, and M.~J. Daly (2011).
\newblock Testing for an unusual distribution of rare variants.
\newblock {\em PLoS genetics\/}~{\em 7\/}(3), e1001322.

\bibitem[\protect\citeauthoryear{Pearson}{Pearson}{1950}]{pearson1950questions}
Pearson, E. (1950).
\newblock On questions raised by the combination of tests based on discontinuous distributions.
\newblock {\em Biometrika\/}~{\em 37\/}(3/4), 383--398.

\bibitem[\protect\citeauthoryear{Rice}{Rice}{1990}]{Rice}
Rice, W.~R. (1990).
\newblock A consensus combined p-value test and the family-wide significance of component tests.
\newblock {\em Biometrics\/}~{\em 46\/}(2), 303--308.

\bibitem[\protect\citeauthoryear{Routledge}{Routledge}{1994}]{routledge1994practicing}
Routledge, R. (1994).
\newblock Practicing safe statistics with the mid-p.
\newblock {\em Canadian Journal of Statistics\/}~{\em 22\/}(1), 103--110.

\bibitem[\protect\citeauthoryear{Seneta}{Seneta}{1999}]{adjmidp}
Seneta, E., B. G. . M.~P. (1999).
\newblock Adjustment to lancaster's mid-p.
\newblock {\em Methodology and Computing in Applied Probability\/}~(1), 229--240.

\bibitem[\protect\citeauthoryear{Shao}{Shao}{2006}]{shao2006}
Shao, J. (2006).
\newblock {\em Mathematical Statistics: Exercises and Solutions}.
\newblock Springer New York.

\bibitem[\protect\citeauthoryear{Stouffer, Suchman, DeVinney, Star, and Williams}{Stouffer et~al.}{1949}]{Stouffer1949}
Stouffer, S.~A., E.~A. Suchman, L.~C. DeVinney, S.~A. Star, and R.~M. Williams (1949).
\newblock {\em The American Soldier: Adjustment during Army Life}, Volume~I.
\newblock New Jersey: Princeton University Press.

\bibitem[\protect\citeauthoryear{van~der Vaart}{van~der Vaart}{2000}]{van2000asymptotic}
van~der Vaart, A. (2000).
\newblock {\em Asymptotic Statistics}.
\newblock Cambridge University Press.

\bibitem[\protect\citeauthoryear{Villani}{Villani}{2003}]{villani2003topics}
Villani, C. (2003).
\newblock {\em Topics in Optimal Transportation}.
\newblock Graduate studies in mathematics. American Mathematical Society.

\bibitem[\protect\citeauthoryear{Wallenius}{Wallenius}{1963}]{wallenius1963biased}
Wallenius, K. (1963).
\newblock {\em Biased Sampling: the Noncentral Hypergeometric Probability Distribution}.
\newblock Stanford University.

\bibitem[\protect\citeauthoryear{Whitlock}{Whitlock}{2005}]{whitlock2005combining}
Whitlock, M.~C. (2005).
\newblock Combining probability from independent tests: the weighted z-method is superior to fisher's approach.
\newblock {\em Journal of evolutionary biology\/}~{\em 18\/}(5), 1368--1373.

\end{thebibliography}


\end{document}